\documentclass{amsart}

\numberwithin{equation}{section}

\usepackage{graphicx,amsthm,enumerate}
\usepackage{verbatim,amsmath,amscd,amssymb,color}
\usepackage[mathscr]{eucal}
\begin{document}
\renewcommand{\labelenumi}{$($\roman{enumi}$)$}
\renewcommand{\labelenumii}{$(${\rm \alph{enumii}}$)$}
\font\germ=eufm10
\newcommand{\cI}{{\mathcal I}}
\newcommand{\cA}{{\mathcal A}}
\newcommand{\cB}{{\mathcal B}}
\newcommand{\cC}{{\mathcal C}}
\newcommand{\cD}{{\mathcal D}}
\newcommand{\cE}{{\mathcal E}}
\newcommand{\cF}{{\mathcal F}}
\newcommand{\cG}{{\mathcal G}}
\newcommand{\cH}{{\mathcal H}}
\newcommand{\cK}{{\mathcal K}}
\newcommand{\cL}{{\mathcal L}}
\newcommand{\cM}{{\mathcal M}}
\newcommand{\cN}{{\mathcal N}}
\newcommand{\cO}{{\mathcal O}}
\newcommand{\cR}{{\mathcal R}}
\newcommand{\cS}{{\mathcal S}}
\newcommand{\cT}{{\mathcal T}}
\newcommand{\cV}{{\mathcal V}}
\newcommand{\cX}{{\mathcal X}}
\newcommand{\cY}{{\mathcal Y}}
\newcommand{\fra}{\mathfrak a}
\newcommand{\frb}{\mathfrak b}
\newcommand{\frc}{\mathfrak c}
\newcommand{\frd}{\mathfrak d}
\newcommand{\fre}{\mathfrak e}
\newcommand{\frf}{\mathfrak f}
\newcommand{\frg}{\mathfrak g}
\newcommand{\frh}{\mathfrak h}
\newcommand{\fri}{\mathfrak i}
\newcommand{\frj}{\mathfrak j}
\newcommand{\frk}{\mathfrak k}
\newcommand{\frI}{\mathfrak I}
\newcommand{\fm}{\mathfrak m}
\newcommand{\frn}{\mathfrak n}
\newcommand{\frp}{\mathfrak p}
\newcommand{\fq}{\mathfrak q}
\newcommand{\frr}{\mathfrak r}
\newcommand{\frs}{\mathfrak s}
\newcommand{\frt}{\mathfrak t}
\newcommand{\fru}{\mathfrak u}
\newcommand{\frA}{\mathfrak A}
\newcommand{\frB}{\mathfrak B}
\newcommand{\frF}{\mathfrak F}
\newcommand{\frG}{\mathfrak G}
\newcommand{\frH}{\mathfrak H}
\newcommand{\frJ}{\mathfrak J}
\newcommand{\frN}{\mathfrak N}
\newcommand{\frP}{\mathfrak P}
\newcommand{\frT}{\mathfrak T}
\newcommand{\frU}{\mathfrak U}
\newcommand{\frV}{\mathfrak V}
\newcommand{\frX}{\mathfrak X}
\newcommand{\frY}{\mathfrak Y}
\newcommand{\fry}{\mathfrak y}
\newcommand{\frZ}{\mathfrak Z}
\newcommand{\frx}{\mathfrak x}
\newcommand{\rA}{\mathrm{A}}
\newcommand{\rC}{\mathrm{C}}
\newcommand{\rd}{\mathrm{d}}
\newcommand{\rB}{\mathrm{B}}
\newcommand{\rD}{\mathrm{D}}
\newcommand{\rE}{\mathrm{E}}
\newcommand{\rH}{\mathrm{H}}
\newcommand{\rK}{\mathrm{K}}
\newcommand{\rL}{\mathrm{L}}
\newcommand{\rM}{\mathrm{M}}
\newcommand{\rN}{\mathrm{N}}
\newcommand{\rR}{\mathrm{R}}
\newcommand{\rT}{\mathrm{T}}
\newcommand{\rZ}{\mathrm{Z}}
\newcommand{\bbA}{\mathbb A}
\newcommand{\bbB}{\mathbb B}
\newcommand{\bbC}{\mathbb C}
\newcommand{\bbG}{\mathbb G}
\newcommand{\bbF}{\mathbb F}
\newcommand{\bbH}{\mathbb H}
\newcommand{\bbP}{\mathbb P}
\newcommand{\bbN}{\mathbb N}
\newcommand{\bbQ}{\mathbb Q}
\newcommand{\bbR}{\mathbb R}
\newcommand{\bbV}{\mathbb V}
\newcommand{\bbZ}{\mathbb Z}
\newcommand{\adj}{\operatorname{adj}}
\newcommand{\Ad}{\mathrm{Ad}}
\newcommand{\Ann}{\mathrm{Ann}}
\newcommand{\rcris}{\mathrm{cris}}
\newcommand{\ch}{\mathrm{ch}}
\newcommand{\coker}{\mathrm{coker}}
\newcommand{\diag}{\mathrm{diag}}
\newcommand{\Diff}{\mathrm{Diff}}
\newcommand{\Dist}{\mathrm{Dist}}
\newcommand{\rDR}{\mathrm{DR}}
\newcommand{\ev}{\mathrm{ev}}
\newcommand{\Ext}{\mathrm{Ext}}
\newcommand{\cExt}{\mathcal{E}xt}
\newcommand{\fin}{\mathrm{fin}}
\newcommand{\Frac}{\mathrm{Frac}}
\newcommand{\GL}{\mathrm{GL}}
\newcommand{\Hom}{\mathrm{Hom}}
\newcommand{\hd}{\mathrm{hd}}
\newcommand{\rht}{\mathrm{ht}}
\newcommand{\id}{\mathrm{id}}
\newcommand{\im}{\mathrm{im}}
\newcommand{\inc}{\mathrm{inc}}
\newcommand{\ind}{\mathrm{ind}}
\newcommand{\coind}{\mathrm{coind}}
\newcommand{\Lie}{\mathrm{Lie}}
\newcommand{\Max}{\mathrm{Max}}
\newcommand{\mult}{\mathrm{mult}}
\newcommand{\op}{\mathrm{op}}
\newcommand{\ord}{\mathrm{ord}}
\newcommand{\pt}{\mathrm{pt}}
\newcommand{\qt}{\mathrm{qt}}
\newcommand{\rad}{\mathrm{rad}}
\newcommand{\res}{\mathrm{res}}
\newcommand{\rgt}{\mathrm{rgt}}
\newcommand{\rk}{\mathrm{rk}}
\newcommand{\SL}{\mathrm{SL}}
\newcommand{\soc}{\mathrm{soc}}
\newcommand{\Spec}{\mathrm{Spec}}
\newcommand{\St}{\mathrm{St}}
\newcommand{\supp}{\mathrm{supp}}
\newcommand{\Tor}{\mathrm{Tor}}
\newcommand{\Tr}{\mathrm{Tr}}
\newcommand{\wt}{\mathrm{wt}}
\newcommand{\Ab}{\mathbf{Ab}}
\newcommand{\Alg}{\mathbf{Alg}}
\newcommand{\Grp}{\mathbf{Grp}}
\newcommand{\Mod}{\mathbf{Mod}}
\newcommand{\Sch}{\mathbf{Sch}}\newcommand{\bfmod}{{\bf mod}}
\newcommand{\Qc}{\mathbf{Qc}}
\newcommand{\Rng}{\mathbf{Rng}}
\newcommand{\Top}{\mathbf{Top}}
\newcommand{\Var}{\mathbf{Var}}
\newcommand{\pmbx}{\pmb x}
\newcommand{\pmby}{\pmb y}
\newcommand{\gromega}{\langle\omega\rangle}
\newcommand{\lbr}{\begin{bmatrix}}
\newcommand{\rbr}{\end{bmatrix}}
\newcommand{\cd}{commutative diagram }
\newcommand{\SpS}{spectral sequence}
\newcommand\C{\mathbb C}
\newcommand\hh{{\hat{H}}}
\newcommand\eh{{\hat{E}}}
\newcommand\F{\mathbb F}
\newcommand\fh{{\hat{F}}}
\newcommand\Z{{\mathbb Z}}
\newcommand\Zn{\Z_{\geq0}}
\newcommand\et[1]{\tilde{e}_{#1}}
\newcommand\ft[1]{\tilde{f}_{#1}}

\def\ge{\frg}
\def\AA{{\mathcal A}}
\def\al{\alpha}
\def\bq{B_q(\ge)}
\def\bqm{B_q^-(\ge)}
\def\bqz{B_q^0(\ge)}
\def\bqp{B_q^+(\ge)}
\def\beneme{\begin{enumerate}}
\def\beq{\begin{equation}}
\def\beqn{\begin{eqnarray}}
\def\beqnn{\begin{eqnarray*}}
\def\bfi{{\mathbf i}}
\def\bfii0{{\bf i_0}}
\def\bigsl{{\hbox{\fontD \char'54}}}
\def\bbra#1,#2,#3{\left\{\begin{array}{c}\hspace{-5pt}
#1;#2\\ \hspace{-5pt}#3\end{array}\hspace{-5pt}\right\}}
\def\cd{\cdots}
\def\ci(#1,#2){c_{#1}^{(#2)}}
\def\Ci(#1,#2){C_{#1}^{(#2)}}
\def\mpp(#1,#2,#3){#1^{(#2)}_{#3}}
\def\bCi(#1,#2){\ovl C_{#1}^{(#2)}}
\def\ch(#1,#2){c_{#2,#1}^{-h_{#1}}}
\def\cc(#1,#2){c_{#2,#1}}
\def\CC{\mathbb{C}}
\def\CBL{\cB_L(\TY(B,1,n+1))}
\def\CBM{\cB_M(\TY(B,1,n+1))}
\def\CVL{\cV_L(\TY(D,1,n+1))}
\def\CVM{\cV_M(\TY(D,1,n+1))}
\def\ddd{\hbox{\germ D}}
\def\del{\delta}
\def\Del{\Delta}
\def\Delr{\Delta^{(r)}}
\def\Dell{\Delta^{(l)}}
\def\Delb{\Delta^{(b)}}
\def\Deli{\Delta^{(i)}}
\def\Delre{\Delta^{\rm re}}
\def\di(#1,#2){D_{#1}^{(#2)}}
\def\dbi(#1,#2){\ovl D_{#1}^{(#2)}}
\def\ei{e_i}
\def\eit{\tilde{e}_i}
\def\eneme{\end{enumerate}}
\def\ep{\epsilon}
\def\eeq{\end{equation}}
\def\eeqn{\end{eqnarray}}
\def\eeqnn{\end{eqnarray*}}
\def\fit{\tilde{f}_i}
\def\FF{{\rm F}}
\def\ft{\tilde{f}_}
\def\gau#1,#2{\left[\begin{array}{c}\hspace{-5pt}#1\\
\hspace{-5pt}#2\end{array}\hspace{-5pt}\right]}
\def\gl{\hbox{\germ gl}}
\def\hom{{\hbox{Hom}}}
\def\ify{\infty}
\def\io{\iota}
\def\ji(#1,#2){j_{#1}^{(#2)}}
\def\kp{k^{(+)}}
\def\km{k^{(-)}}
\def\llra{\relbar\joinrel\relbar\joinrel\relbar\joinrel\rightarrow}
\def\lan{\langle}
\def\lar{\longrightarrow}
\def\lm{\lambda}
\def\Lm{\Lambda}
\def\mapright#1{\smash{\mathop{\longrightarrow}\limits^{#1}}}
\def\Mapright#1{\smash{\mathop{\Longrightarrow}\limits^{#1}}}
\def\mm{{\bf{\rm m}}}
\def\nd{\noindent}
\def\nn{\nonumber}
\def\nnn{\hbox{\germ n}}
\def\catob{{\mathcal O}(B)}
\def\oint{{\mathcal O}_{\rm int}(\ge)}
\def\ot{\otimes}
\def\op{\oplus}
\def\opi{\ovl\pi_{\lm}}
\def\osigma{\ovl\sigma}
\def\ovl{\overline}
\def\plm{\Psi^{(\lm)}_{\io}}
\def\qq{\qquad}
\def\q{\quad}
\def\qed{\hfill\framebox[2mm]{}}
\def\QQ{\mathbb Q}
\def\qi{q_i}
\def\qii{q_i^{-1}}
\def\ra{\rightarrow}
\def\ran{\rangle}
\def\rlm{r_{\lm}}
\def\ssl{\hbox{\germ sl}}
\def\slh{\widehat{\ssl_2}}
\def\syl{\scriptstyle}
\def\ti{t_i}
\def\tii{t_i^{-1}}
\def\til{\tilde}
\def\tm{\times}
\def\tri{\bigtriangleup}
\def\tt{\frt}
\def\TY(#1,#2,#3){#1^{(#2)}_{#3}}
\def\ua{U_{\AA}}
\def\ue{U_{\vep}}
\def\uq{U_q(\ge)}
\def\uqp{U'_q(\ge)}
\def\ufin{U^{\rm fin}_{\vep}}
\def\ufinp{(U^{\rm fin}_{\vep})^+}
\def\ufinm{(U^{\rm fin}_{\vep})^-}
\def\ufinz{(U^{\rm fin}_{\vep})^0}
\def\uqm{U^-_q(\ge)}
\def\uqmq{{U^-_q(\ge)}_{\bf Q}}
\def\uqpm{U^{\pm}_q(\ge)}
\def\uqq{U_{\bf Q}^-(\ge)}
\def\uqz{U^-_{\bf Z}(\ge)}
\def\ures{U^{\rm res}_{\AA}}
\def\urese{U^{\rm res}_{\vep}}
\def\uresez{U^{\rm res}_{\vep,\ZZ}}
\def\util{\widetilde\uq}
\def\uup{U^{\geq}}
\def\ulow{U^{\leq}}
\def\bup{B^{\geq}}
\def\blow{\ovl B^{\leq}}
\def\vep{\varepsilon}
\def\vp{\varphi}
\def\vpi{\varphi^{-1}}
\def\VV{{\mathcal V}}
\def\xii{\xi^{(i)}}
\def\Xiioi{\Xi_{\io}^{(i)}}
\def\xxi(#1,#2,#3){\displaystyle {}^{#1}\Xi^{(#2)}_{#3}}
\def\xsi(#1,#2,#3){\displaystyle {}^{#1}\Sigma^{(#2)}_{#3}}
\def\xE(#1,#2,#3){\displaystyle {}^{#1}E_{#2}[#3]}
\def\xF(#1,#2){\displaystyle {}^{#1}F_{#2}}
\def\xx(#1,#2){\displaystyle {}^{#1}\Xi_{#2}}
\def\W1{W(\varpi_1)}
\def\WW{{\mathcal W}}
\def\wt{{\rm wt}}
\def\wtil{\widetilde}
\def\what{\widehat}
\def\wpi{\widehat\pi_{\lm}}
\def\ZZ{\mathbb Z}
\def\RR{\mathbb R}

\def\m@th{\mathsurround=0pt}
\def\fsquare(#1,#2){
\hbox{\vrule$\hskip-0.4pt\vcenter to #1{\normalbaselines\m@th
\hrule\vfil\hbox to #1{\hfill$\scriptstyle #2$\hfill}\vfil\hrule}$\hskip-0.4pt
\vrule}}

\newtheorem{pro}{Proposition}[section]
\newtheorem{thm}[pro]{Theorem}
\newtheorem{lem}[pro]{Lemma}
\newtheorem{ex}[pro]{Example}
\newtheorem{cor}[pro]{Corollary}
\newtheorem{conj}[pro]{Conjecture}
\theoremstyle{definition}
\newtheorem{df}[pro]{Definition}

\newcommand{\cmt}{\marginpar}
\newcommand{\seteq}{\mathbin{:=}}
\newcommand{\cl}{\colon}
\newcommand{\be}{\begin{enumerate}}
\newcommand{\ee}{\end{enumerate}}
\newcommand{\bnum}{\be[{\rm (i)}]}
\newcommand{\enum}{\ee}
\newcommand{\ro}{{\rm(}}
\newcommand{\rf}{{\rm)}}
\newcommand{\set}[2]{\left\{#1\,\vert\,#2\right\}}
\newcommand{\sbigoplus}{{\mbox{\small{$\bigoplus$}}}}
\newcommand{\ba}{\begin{array}}
\newcommand{\ea}{\end{array}}
\newcommand{\on}{\operatorname}
\newcommand{\eq}{\begin{eqnarray}}
\newcommand{\eneq}{\end{eqnarray}}
\newcommand{\hs}{\hspace*}

\title[Decorations on Geometric Crystals and 
Monomial Realizations of Crystal Bases ]
{Decorations on Geometric Crystals and 
Monomial Realizations of Crystal Bases for Classical Groups}

\author{Toshiki N\textsc{akashima}}
\address{Department of Mathematics, 
Sophia University, Kioicho 7-1, Chiyoda-ku, Tokyo 102-8554,
Japan}
\email{toshiki@sophia.ac.jp}
\thanks{supported in part by JSPS Grants 
in Aid for Scientific Research $\sharp 22540031$.}

%    General info
\date{}

%\dedicatory{Dedicated to Professor Michio Jimbo on the occasion of 
%his 60th birthday}


\begin{abstract}
We shall describe explicitly 
the decoration functions for certain 
decorated geometric crystals of classical groups
and we shall show that they are
represented in terms of monomial realizations of crystal bases.
\end{abstract}

\maketitle

\tableofcontents
%%%%%%% Sect 1. %%%%%%%%%%%%%%
\renewcommand{\thesection}{\arabic{section}}
\section{Introduction}
\setcounter{equation}{0}
\renewcommand{\theequation}{\thesection.\arabic{equation}}

Since the theory of crystal bases was invented by
Kashiwara(\cite{K0,K1}),
there have been several kinds of realizations for crystal bases, 
{\it e.g.,} tableaux,
paths, polytopes, monomials, etc.
In the article, the monomial realization of crystal bases, which 
is introduced by Nakajima \cite{Nj} and refined by Kashiwara \cite{K7},
will be treated and used to describe the decoration functions for the
decorated geometric crystals (\cite{BK2}, see also below).
Let $\cY$ be the set of Laurant monomials in doubly-indexed variables
$\{Y_{m,i}\}$ $(m\in\bbZ,\,\,i\in I:=\{1,2,\cd,n\})$ as follows:
\[
  \cY:=\{Y=\prod_{m\in\bbZ,i\in I}
Y_{m,i}^{l_{m,i}}|l_{m,i}\in \bbZ\setminus\{0\}\text{
 except for finitely many }(m,i)\}.
\]
We can define the crystal structure on $\cY$ 
by the way as in \ref{mono-subsec} and it is shown that certain
connected component 
$B(Y)\subset \cY$ including the highest monomial $Y$ 
is isomorphic to the crystal $B(\lm)$
where $\lm=\wt(Y)$ is a dominant integral weight.
For example, for type $A_n$ and any integer $m$ we have
\[
 B(\Lm_1)\cong\{Y_{m,1},\frac{Y_{m,2}}{Y_{m+1,1}},\frac{Y_{m,3}}{Y_{m+1,2}},
\cd,\frac{Y_{m,n}}{Y_{m+1,n-1}},\frac{1}{Y_{m,n}}\},
\]
where $\Lm_1$ is the first fundamental weight and $Y_{m,1}$ is the 
highest monomial.

A geometric crystal is a sort of geometric lifting of 
Kashiwara's crystal bases (\cite{BK}), which is generalized to 
the affine/Kac-Moody settings (\cite{KNO,KNO2,N}). In this paper 
we do not treat such general settings and then we shall consider the
simple classical settings below.
Let $\ge$ be a simple complex Lie algebra, $G$ the corresponding
complex algebraic group, $B^{\pm}\subset G$ Borel subgroups and 
$U^{\pm}$ maximal unipotent subgroups such that $U^{\pm}\subset B^{\pm}$.
The notion of decorated 
geometric crystals has been initiated by Berenstein and Kazhdan(\cite{BK2}).
Let $I$ be the index set of the simple roots. 
Associated with 
the Cartan matrix $A=(a_{i,j})_{i,j\in I}$, 
the decorated geometric crystal $\cX=(\chi,f)$ is defined as
a pair of geometric crystal 
$\chi=(X,\{e_i\}_i,\{\gamma_i\}_i,\{\vep_i\}_i)$ 
and  a certain special rational function $f$ on the algebraic variety
$X$ satisfying the condition
\[
 f(e_i^c(x))=f(x)+(c-1)\vp_i(x)+(c^{-1}-1)\vep_i(x),
\]
for any $i\in I$, where $e_i^c$ is the unital rational $\bbC^\times$ action on
$X$, and $\vp_i:=\vep_i\cdot \gamma_i$ 
is the rational functions on $X$. The function $f$ is called the
decoration (function) of the decorated geometric crystal $\cX$.
In \cite{BK2} the ultra-discretization of the 
decoration is used to describe the Kashiwara's crystal
base $B(\lm)$, which is quite similar to the polyhedral realizations 
of crystal bases(\cite{N2,NZ}). 
This similarity let us conceive some link between crystal bases and 
the decorations. 
Indeed, in \cite{N4} we made this link clear for type $A_n$.
Here we shall consider another link between them, which
is the main purpose of this article. 
The purpose here is to present the explicit form of the
decorations for the decorated geometric crystals on $\bbB_{w_0}$ (see
\ref{BBB}) and describe the decorations in terms of the
monomial realizations of crystal bases(\cite{K7,Nj}).
In \cite{N4}, we presented the conjecture 
(see also Conjecture \ref{conj1} below) and 
gave the positive answer for type $A_n$. 
To be more precise, we shall introduce certain part of the 
results in \cite{N4}. First, we consider the geometric crystal structure
on the variety $\bbB_{w_0}:=TB^-_{w_0}$ where $T\subset G$ is the
maximal torus, $w_0$ is the longest element of the Weyl group and 
$B^-_{w_0}:=B^-\cap U\ovl w_0U$. Let $\chi_i:U\to\bbC$ be the elementary
character and $\eta:G\to G$ the positive inverse (see \ref{char}).
Then the decoration $f_B$ on $\bbB_{w_0}$ is defined by the formula
\[
 f_B(g):=\sum_i\chi_i(\pi^+(w_0^{-1}g))+\chi_i(\pi^+(w_0^{-1}\eta(g))),
\]
where $\pi^+:B^-U\to U$ is the projection.
By the definition of $f_B$, it suffices to get the explicit form of 
$\chi_i(\pi^+(w_0^{-1}g))$ and $\chi_i(\pi^+(w_0^{-1}\eta(g)))$ for our 
purpose. 
Furthermore, the elementary characters are expressed by the ``generalized
minors $\Del_{\gamma,\del}$''
(\cite{BFZ,BZ,BZ2})  and then as in (\ref{chi-min}) we have  
\[
f_B(g)= \sum_i\frac{\Del_{w_0\Lm_i,s_i\Lm_i}(g)
+\Del_{w_0s_i\Lm_i,\Lm_i}(g)}
{\Del_{w_0\Lm_i,\Lm_i}(g)},
\]
where $\Lm_i$ is the $i$-th fundamental weight. 
We shall see the explicit forms of $\Del_{w_0\Lm_i,s_i\Lm_i}(g)$
and $\Del_{w_0s_i\Lm_i,\Lm_i}(g)$ as in \eqref{del1} and \eqref{del2}. 
In most cases except for the spin
representations of type $B_n$ and $D_n$, it is performed by 
direct calculations. For the cases of the spin representations, 
we prepare the ``triangles'',
which has some interesting combinatorial properties and is useful to calculate 
the above generalized minors.
Then, we can find their relations to the monomial realizations of crystals
(Sect.\ref{sect-mono}).

In \cite{N4}, we also describe the relations to the polyhedral
realizations explicitly for type $A_n$ though we do not treat 
that part herein. However, we strongly believe that there exist the 
relations similar to the ones for type $A_n$. 
As for the relations to the polyhedral realizations for 
other classical cases, we shall discuss in forthcoming papers.

The organization of the article is as follows:
After the introduction in this section and the preliminaries in Sect.2, 
we review the explicit descriptions for the 
fundamental representation of the classical Lie algebras in Sect.3.
In Sect.4, first we introduce the theory of decorated geometric 
crystals following \cite{BK2}.  Next, we define the decoration by using the 
elementary characters and certain special positive decorated geometric
crystal on ${\bbB}_w=TB^-_w$. Finally, the ultra-discretization of $TB^-_w$ is
described explicitly. 
In Sect.5, 
the theory of monomial realizations would be introduced and we shall see
some duality on monomial realizations.
In Sect.6, we review the generalized minors and their relations to 
our elementary characters and certain bilinear forms.
The main conjecture will be presented at the end of the section.
In the last three sections, we describe the explicit form of decorations
and express them in terms of the monomial realization of crystal bases, 
which means that the conjecture is positively resolved for the classical
groups.

The author would like to acknowledge Masaki Kashiwara for 
discussions and his helpful suggestions.

%%%%%%% Sect 2. %%%%%%%%%%%%%%
\renewcommand{\thesection}{\arabic{section}}
\section{Preliminaries and Notations}
\setcounter{equation}{0}
\renewcommand{\theequation}{\thesection.\arabic{equation}}

We list the notations used in this paper.
Let  
 $A=({\bf a}_{ij})_{i,j\in I}$ be 
an indecomposable Cartan matrix
with a finite index set $I$
(though we can consider more general Kac-Moody setting.).
Let $(\frt,\{\al_i\}_{i\in I},\{h_i\}_{i\in I})$ 
be the associated
root data 
satisfying $\al_j(h_i)={\bf a}_{ij}$ where 
$\al_i\in \tt^*$ is a simple root and 
$h_i\in \tt$ is a simple coroot.
Let $\ge=\ge(A)=\lan \frt,e_i,f_i(i\in I)\ran$ 
be the simple Lie algebra associated with $A$
over $\bbC$ and $\Del=\Del_+\sqcup\Del_-$
be the root system associated with $\ge$, where $\Del_{\pm}$ is 
the set of positive/negative roots.
Let $P \subset \tt^*$ be the weight lattice, 
$\lan h,\lm\ran=\lm(h)$ the pairing between $\tt$ and $\tt^*$,
and $(\al, \beta)$ be an inner product on
$\tt^*$ such that $(\al_i,\al_i)\in 2\bbZ_{\geq 0}$ and
$\lan h_i,\lm\ran={{2(\al_i,\lm)}\over{(\al_i,\al_i)}}$
for $\lm\in\tt^*$.
Let $P^*=\{h\in \tt: \lan h,P\ran\subset\ZZ\}$ and
$P_+:=\{\lm\in P:\lan h_i,\lm\ran\in\ZZ_{\geq 0}\}$.
We call an element in $P_+$ a {\it dominant integral weight}.
Let $\{\Lm_i|i\in I\}$ be the set of the fundamental weights satisfying
$\lan h_i,\Lm_j\ran=\del_{i,j}$ is a $\bbZ$-basis of $P$.

The quantum algebra $\uq$
is an associative
$\QQ(q)$-algebra generated by the $e_i$, $f_i \,\, (i\in I)$,
and $q^h \,\, (h\in P^*)$
satisfying the usual relations, where we use the same notations for the 
generators $e_i$ and $f_i$ as the ones for $\ge$.
The algebra $\uqm$ is the subalgebra of $\uq$ generated 
by the $f_i$ $(i\in I)$.

For the irreducible highest weight module of $\uq$
with the highest weight $\lm\in P_+$, we denote it by $V(\lm)$
and we denote its {\it crystal base}  by $(L(\lm),B(\lm))$.
Similarly, for the crystal base of the algebra $\uqm$ we denote 
$(L(\ify),B(\ify))$ (see \cite{K0,K1}).
Let $\pi_{\lm}:\uqm\longrightarrow V(\lm)\cong \uqm/{\sum_i\uqm
f_i^{1+\lan h_i,\lm\ran}}$ be the canonical projection and 
$\widehat \pi_{\lm}:L(\ify)/qL(\ify)\longrightarrow L(\lm)/qL(\lm)$
be the induced map from $\pi_{\lm}$. Here note that 
$\widehat \pi_{\lm}(B(\ify))=B(\lm)\sqcup\{0\}$.

By the terminology {\it crystal } we mean some combinatorial object 
obtained by abstracting the properties of crystal bases.
Indeed, crystal constitutes a set $B$ and the maps
$wt:B\longrightarrow P$, $\vep_i,\vp_i:B\longrightarrow \ZZ\sqcup\{-\ify\}$
and $\eit,\fit:B\sqcup\{0\}\longrightarrow B\sqcup\{0\}$
($i\in I$) satisfying several axioms (see \cite{K3},\cite{NZ},\cite{N2}).
In fact, $B(\ify)$ and $B(\lm)$ are the typical examples 
of crystals.

%%%%%%% Sect.3 %%%%%%%%%%%%%%%%%
\renewcommand{\thesection}{\arabic{section}}
\section{Fundamental Representations}
\setcounter{equation}{0}
\renewcommand{\theequation}{\thesection.\arabic{equation}}

%%%%%%%%%%%%%%%%%%%%%%%%%%%%%%%%%%%
\subsection{Type $A_n$}%3.1

Let $V_1:=V(\Lm_1)$ be the vector representation of 
$\ssl_{n+1}(\bbC)$ with
the standard basis $\{v_1,\cd,v_{n+1}\}$, and 
$\{e_i,f_i,h_i\}_{i=1,\cd,n}$ the 
Chevalley generators of $\ssl_{n+1}(\bbC)$.
Their actions on the basis vectors are as follows:
\begin{equation}
e_iv_j=\begin{cases}v_{i}&\text{ if }j=i+1,\\
0&\text{otherwise},
\end{cases}
\q
f_iv_j=\begin{cases}v_{i+1}&\text{ if }j=i,\\
0&\text{otherwise},
\end{cases}
\q
h_iv_j=\begin{cases}v_{i}&\text{ if }j=i,\\
-v_{i+1}&\text{ if }j=i+1,\\
0&\text{otherwise},
\end{cases}
\end{equation}

%%%%%%%%%%%%%%%%%%%%%%%%%%%%%%%%%%%%
\subsection{Type $C_n$}%3.2

Let $I:=\{1,2,\cd ,n\}$ be the index set of the 
simple roots of type $C_n$. The Cartan matrix 
$A=({\bf a}_{i,j})_{i,j\in I}$ of type 
$C_n$ is given by 
\[
 {\bf a}_{i,j}=\begin{cases}
2&\text{if }i=j,\\
-1&\text{if }|i-j|=1\text{ and }(i,j)\ne(n-1,n)\\
-2&\text{if }(i,j)=(n-1,n),\\
0&\text{otherwise}.
\end{cases}
\]
Here $\al_i$ ($i\ne n$) is a short root and 
$\al_n$ is the long root.
Let $\{h_i\}_{i\in I}$ be the set of the simple co-roots
and $\{\Lm_i\}_{i\in I}$ be the set of the fundamental 
weights satisfying $\al_j(h_i)={\bf a}_{i,j}$ 
and $\Lm_i(h_j)=\del_{i,j}$.

First, let us describe the vector representation 
$V(\Lm_1)$. Set ${\mathbf B}^{(n)}:=
\{v_i,v_{\ovl i}|i=1,2,\cd,n.\}$ and define 
$V(\Lm_1):=\bigoplus_{v\in{\mathbf B}^{(n)}}\bbC v$. 
The weight of $v_i$ is as follows:
\[
 {\rm wt}(v_i)=\begin{cases}
\Lm_i-\Lm_{i-1}&\text{ if }i=1,\cd,n,\\
\Lm_{i-1}-\Lm_{i}&\text{ if }i=\ovl 1,\cd,\ovl n,
\end{cases}
\]
where $\Lm_0=0$.
The actions of $e_i$ and $f_i$ 
are given by:
\begin{eqnarray}
&&f_iv_i=v_{i+1},\q f_iv_{\ovl i+1}=v_{\ovl i},\q
e_iv_{i+1}=v_i,\q e_iv_{\ovl i}=v_{\ovl{i+1}}
\q(1\leq i<n),\label{c-f1}\\
&&f_nv_n=v_{\ovl n},\qq 
e_nv_{\ovl n}=v_n,\label{c-f2}
\end{eqnarray}
and the other actions are trivial.

Let $\Lm_i$ be the $i$-th fundamental weight of type $C_n$.
As is well-known that the fundamental representation 
$V(\Lm_i)$ $(1\leq i\leq n)$
is embedded in $V(\Lm_1)^{\ot i}$
with multiplicity free.
The explicit form of the highest(resp. lowest) weight 
vector $u_{\Lm_i}$ (resp. $v_{\Lm_i}$)
of $V(\Lm_i)$ is realized in 
$V(\Lm_1)^{\ot i}$ as follows:
\begin{equation}
\begin{array}{ccc}\displaystyle
u_{\Lm_i}&=&\displaystyle
\sum_{\sigma\in{\mathfrak S}_i}{\rm sgn}(\sigma)
v_{\sigma(1)}\ot\cd\ot v_{\sigma(i)},\\
v_{\Lm_i}&=&\displaystyle
\sum_{\sigma\in{\mathfrak S}_i}{\rm sgn}(\sigma)
v_{\ovl{\sigma(i)}}\ot\cd\ot v_{\ovl{\sigma(1)}}, 
\end{array}
\label{h-l}
\end{equation}
where ${\mathfrak S}_i$ is the $i$-th symmetric group.

We review the crystal $B(\Lm_k)$ following \cite{KN,N0}.
Set the order  on the set $J=\{i,\ovl i|1\leq i\leq n\}$ by 
\[
 1< 2<\cd< n-1< n
< \ovl n< \ovl{n-1}< \cd< \ovl
 2< \ovl 1
\]
Then, the crystal $B(\Lm_k)$ is described:
\begin{equation}
B(\Lm_k)=\left\{[j_1,\cd,j_k]\,\,\begin{array}{|l}
1\leq  j_1<\cd< j_k\leq\ovl 1,\\
\text{if }\ovl{j_a}=j_b=\ovl i,\text{ then }a+b\leq i
\end{array}
\right\}
\label{Ytab-c}
\end{equation}
Note that in \cite{KN} the vector $[j_1,\cd,j_k]$ is represented 
as the column Young tableau.
Note also that 
the highest(resp. lowest) weight vector in $B(\Lm_k)$ is 
$[1,2,\cd,k]$ (resp. $[\ovl k,\ovl{k-1},\cd, \ovl 2,\ovl 1]$).
%%%%%%%%%%%  %%%%%%%%%%
\subsection{Type $B_n$}\label{Bn}%3.3

Let $I:=\{1,2,\cd,n\}$ be the index set of the 
simple roots of type $B_n$. The Cartan matrix 
$A=({\bf a}_{i,j})_{i,j\in I}$ of type 
$B_n$ is given by 
\[
 {\bf a}_{i,j}=\begin{cases}
2&\text{if }i=j,\\
-1&\text{if }|i-j|=1\text{ and }(i,j)\ne(n,n-1)\\
-2&\text{if }(i,j)=(n,n-1),\\
0&\text{otherwise}.
\end{cases}
\]
Here $\al_i$ ($i\ne n$) is a long root and 
$\al_n$ is the short root.
Let $\{h_i\}_{i\in I}$ be the set of the simple co-roots
and $\{\Lm_i\}_{i\in I}$ be the set of the fundamental 
weights satisfying $\al_j(h_i)={\bf a}_{i,j}$ 
and $\Lm_i(h_j)=\del_{i,j}$.

First, let us describe the vector representation $V(\Lm_1)$ for 
$B_n$.
Set ${\mathbf B}^{(n)}:=\{v_i,\,v_{\ovl i}|i=1,2,\cd,n\}\cup\{v_0\}$ and 
$V(\Lm_1):=\bigoplus_{v\in {\mathbf B}^{(n)}}\bbC v$.
The weight of $v_i$ is as follows:
\begin{eqnarray*}
&& {\rm wt}(v_i)=
\Lm_i-\Lm_{i-1},\q
{\rm wt}(v_{\ovl i})=\Lm_{i-1}-\Lm_{i}\q
(i=1,\cd,n-1),\\
&&{\rm wt}(v_n)=
2\Lm_n-\Lm_{n-1},\q
{\rm wt}(v_{\ovl n})=
\Lm_{n-1}-2\Lm_n,\q {\rm wt}(v_0)=0,
\end{eqnarray*}
where $\Lm_0=0$.
The actions of $e_i$ and $f_i$ 
are given by:
\begin{eqnarray}
&&f_iv_i=v_{i+1},\q f_iv_{\ovl{i+1}}=v_{\ovl i},\q
e_iv_{i+1}=v_i,\q e_iv_{\ovl i}=v_{\ovl{i+1}}
\q(1\leq i<n),\label{b-f1}\\
&&f_nv_n=v_0,\q f_nv_0=2v_{\ovl n},\q 
e_nv_0=2v_n,\q e_nv_{\ovl n}=v_0,\label{b-f2}
\end{eqnarray}
and the other actions are trivial.
For $i=1,2,\cd, n-1$, 
the $i$-th fundamental representation $V(\Lm_i)$ is realized in
$V(\Lm_1)^{\ot i}$ as the case $C_n$ and their highest (resp. lowest)
weight vector $u_{\Lm_i}$ (resp. $v_{\Lm_i}$) is given by the formula 
\eqref{h-l}.

The last fundamental representation $V(\Lm_n)$ 
is called the ``spin representation'' whose dimension
is $2^n$. It is realized as follows:
Set $V^{(n)}_{sp}:=\bigoplus_{\ep\in B^{(n)}_{sp}}\bbC\ep$ where
\[
 B^{(n)}_{sp}:=\{(\ep_1,\cd,\ep_n)|\ep_i\in\{+,-\}\,(i=1,2,\cd,n)\}.
\]
Define the explicit actions of $h_i$, $e_i$ and $f_i$ on 
$V^{(n)}_{sp}$ by 
\begin{eqnarray}
h_i(\ep_1,\cd,\ep_n)&=&\begin{cases}
\frac{\ep_i\cdot1-\ep_{i+1}\cdot1}{2}
(\ep_1,\cd,\ep_n),&\text{ if }i\ne n,\\
\ep_n(\ep_1,\cd,\ep_n)&\text{ if }i=n,
\end{cases}\\
f_i(\ep_1,\cd,\ep_n)&=&\begin{cases}
(\cd,\buildrel{i}\over-,\buildrel{i+1}\over+,\cd)
&\text{ if }\ep_i=+,\,\,\ep_{i+1}=-,\,\,
i\ne n,\\
(\cd\cd\cd,\buildrel{n}\over-)
&\text{ if }\ep_n=+,\,\,i=n,\\
\qq0&\text{otherwise}
\end{cases}
\label{B-fi}
\\
e_i(\ep_1,\cd,\ep_n)&=&\begin{cases}
(\cd,\buildrel{i}\over+,\buildrel{i+1}\over-,\cd)
&\text{ if }\ep_i=-,\,\,\ep_{i+1}=+,\,\,
i\ne n,\\
(\cd\cd\cd,\buildrel{n}\over+)
&\text{ if }\ep_n=-,\,\,i=n,\\
\qq0&\text{otherwise.}
\end{cases}
\label{B-ei}
\end{eqnarray}
Then the module $V^{(n)}_{sp}$ is isomorphic to 
$V(\Lm_n)$ as a $B_n$-module.

{\sl Remark.}
We can associate the crystal structure on the set $B^{(n)}_{sp}$ by 
setting $\fit=f_i$ and $\eit=e_i$ in (\ref{B-fi}) and (\ref{B-ei})
respectively, which is also denoted by $B^{(n)}_{sp}$ and is isomorphic
to $B(\Lm_n)$.

%%%%%%%%%%%%%%%%%%%%%%%%%%%%%%
\subsection{Type $D_n$}\label{Dn}%3.4

Let $I:=\{1,2,\cd,n\}$ be the index set of the 
simple roots of type $D_n$. The Cartan matrix 
$A=({\bf a}_{i,j})_{i,j\in I}$ of type 
$D_n$ is as follows:
\[
 {\bf a}_{i,j}=\begin{cases}
2&\text{if }i=j,\\
-1&\text{if }|i-j|=1
\text{ and }(i,j)\ne(n,n-1),\,(n-1,n), 
\text{ or }(i,j)=(n-2,n),\,\,(n,n-2)\\
0&\text{otherwise}.
\end{cases}
\]
Let $\{h_i\}_{i\in I}$ be the set of the simple co-roots
and $\{\Lm_i\}_{i\in I}$ be the set of the fundamental 
weights satisfying $\al_j(h_i)={\bf a}_{i,j}$ 
and $\Lm_i(h_j)=\del_{i,j}$.

First, let us describe the vector representation $V(\Lm_1)$ for 
$D_n$.
Set ${\mathbf B}^{(n)}:=\{v_i, v_{\ovl i}|i=1,2,\cd,n\}$. 
The weight of $v_i$ is as follows:
\begin{eqnarray*}
&& {\rm wt}(v_i)=
\Lm_i-\Lm_{i-1},\q
{\rm wt}(v_{\ovl i})=\Lm_{i-1}-\Lm_{i}\q
(i=1,\cd,n-1),\\
&&{\rm wt}(v_n)=
\Lm_{n-1}+\Lm_n-\Lm_{n-2},\q
{\rm wt}(v_{\ovl n})=
\Lm_{n-2}-\Lm_{n-1}+\Lm_n,
\end{eqnarray*}
where $\Lm_0=0$.
The actions of $e_i$ and $f_i$ 
are given by:
\begin{eqnarray}
&&f_iv_i=v_{i+1},\q f_iv_{\ovl{i+1}}=v_{\ovl i},\q
e_iv_{i+1}=v_i,\q e_iv_{\ovl i}=v_{\ovl{i+1}}
\q(1\leq i<n),\label{d-fi}\\
&&f_nv_n=v_{\ovl{n-1}},\q f_{n-1}v_{\ovl n}=v_{\ovl{n-1}},\q 
e_{n-1}v_{\ovl {n-1}}=v_{\ovl n},\q e_nv_{\ovl{n-1}}=v_n,
\label{d-fn}
\end{eqnarray}
and the other actions are trivial.
For $i=1,2,\cd, n-2$, 
the $i$-th fundamental representation $V(\Lm_i)$ is realized in
$V(\Lm_1)^{\ot i}$ as the cases $B_n$ and $C_n$ and their highest (resp.lowest)
weight vector $u_{\Lm_i}$ (resp. $v_{\Lm_i}$) is given by the formula 
\eqref{h-l}.

The last two fundamental representations $V(\Lm_{n-1})$ and $V(\Lm_n)$ 
are also called the ``spin representations'' whose dimensions
are $2^{n-1}$. They are realized as follows:
Set $V^{(+,n)}_{sp}$(resp. $V^{(-,n)}_{sp}$)
$:=\bigoplus_{\ep\in
B^{(+,n)}_{sp}(\text{resp. }B^{(-,n)}_{sp})}\bbC\ep$ 
where
\[
 B^{(+,n)}_{sp}(\text{resp. }B^{(-,n)}_{sp})
:=\{(\ep_1,\cd,\ep_n)|\ep_i\in\{+,-\},\ep_1\cd\ep_n=+(\text{resp. }-)\}.
\]
Define the explicit actions of $h_i$, $e_i$ and $f_i$ on 
$V^{(\pm,n)}_{sp}$ by 
\begin{eqnarray}
h_i(\ep_1,\cd,\ep_n)&=&\begin{cases}
\frac{\ep_i\cdot1-\ep_{i+1}\cdot1}{2}
(\ep_1,\cd,\ep_n),&\text{ if }i\ne n,\\
\frac{\ep_{n-1}\cdot1+\ep_n\cdot1}{2}(\ep_1,\cd,\ep_n)&\text{ if }i=n,
\end{cases}\\
f_i(\ep_1,\cd,\ep_n)&=&\begin{cases}
(\cd,\buildrel{i}\over-,\buildrel{i+1}\over+,\cd)
&\text{ if }\ep_i=+,\,\,\ep_{i+1}=-,\,\,
i\ne n,\\
(\cd\cd\cd,\buildrel{n-1}\over-,\buildrel{n}\over-)
&\text{ if }\ep_{n-1}=+,\ep_n=+,\,\,i=n,\\
\qq0&\text{otherwise}
\end{cases}
\label{D-fi}
\\
e_i(\ep_1,\cd,\ep_n)&=&\begin{cases}
(\cd,\buildrel{i}\over+,\buildrel{i+1}\over-,\cd)
&\text{ if }\ep_i=-,\,\,\ep_{i+1}=+,\,\,
i\ne n,\\
(\cd\cd\cd,\buildrel{n-1}\over+,\buildrel{n}\over+)
&\text{ if }\ep_{n-1}=-,\ep_n=-,\,\,i=n,\\
\qq0&\text{otherwise.}
\end{cases}
\label{D-ei}
\end{eqnarray}
Then the module $V^{(+,n)}_{sp}$ (resp. $V^{(-,n)}_{sp}$) is isomorphic to 
$V(\Lm_n)$ (resp. $V(\Lm_{n-1})$) as a $D_n$-module.

{\sl Remark.}
Similar to the case $B_n$, in this case 
we can associate the crystal structure on the set $B^{(+,n)}_{sp}$ 
(resp. $B^{(-,n)}_{sp}$) by 
setting $\fit=f_i$ and $\eit=e_i$ in (\ref{D-fi}) and (\ref{D-ei})
respectively, which is also denoted by $B^{(\pm,n)}_{sp}$ and is isomorphic
to $B(\Lm_n)$ (resp. $B(\Lm_{n-1})$).

%%%%%%% Sect 4. %%%%%%%%%%%%%%
\renewcommand{\thesection}{\arabic{section}}
\section{Decorated geometric crystals}
\setcounter{equation}{0}
\renewcommand{\theequation}{\thesection.\arabic{equation}}

The basic reference for this section is \cite{BK,BK2,N}.
%%%%%%%%%%%%%%%%%%%%%%%%%%%%%%
\subsection{Definitions}%4.1

Let  
 $A=({\bf a}_{ij})_{i,j\in I}$ be 
an indecomposable Cartan matrix.
Let $\ge=\ge(A)=\lan \frt,e_i,f_i(i\in I)\ran$ 
be the simple Lie algebra associated with $A$
over $\bbC$ as above and $\Del=\Del_+\sqcup\Del_-$
be the root system associated with $\ge$.
Define the simple reflections $s_i\in{\rm Aut}(\frt)$ $(i\in I)$ by
$s_i(h)\seteq h-\al_i(h)h_i$, which generate the Weyl group $W$.
Let $G$ be the simply connected simple algebraic group 
over $\bbC$ whose Lie algebra is $\ge=\frn_+\oplus \frt\oplus \frn_-$, 
which is the usual triangular decomposition.
Let $U_{\al}\seteq\exp\ge_{\al}$ $(\al\in \Del)$
be the one-parameter subgroup of $G$.
The group $U^\pm$ are generated by 
$\{U_\al|\al\in\Del_{\pm}\}$.
Here $U^\pm$ is a unipotent radical of $G$ and 
${\rm Lie}(U^\pm)=\frn_{\pm}$.
For any $i\in I$, there exists 
a unique group homomorphism
$\phi_i\cl SL_2(\bbC)\rightarrow G$ such that
\[
\phi_i\left(
\left(
\begin{array}{cc}
1&t\\
0&1
\end{array}
\right)\right)=\exp(t e_i),\,\,
 \phi_i\left(
\left(
\begin{array}{cc}
1&0\\
t&1
\end{array}
\right)\right)=\exp(t f_i)\qquad(t\in\bbC).
\]
Set $\al^\vee_i(c)\seteq
\phi_i\left(\left(
\begin{smallmatrix}
c&0\\
0&c^{-1}\end{smallmatrix}\right)\right)$,
$x_i(t)\seteq\exp{(t e_i)}$, $y_i(t)\seteq\exp{(t f_i)}$, 
$G_i\seteq\phi_i(SL_2(\bbC))$,
$T_i\seteq \alpha_i^\vee(\bbC^\times)$ 
and 
$N_i\seteq N_{G_i}(T_i)$. Let
$T$ be a maximal torus of $G$ 
which has $P$ as its weight lattice and 
Lie$(T)=\frt$.
Let 
$B^{\pm}(\supset T,U^{\pm})$ be the Borel subgroup of $G$.
We have the isomorphism
$\phi:W\mapright{\sim}N/T$ defined by $\phi(s_i)=N_iT/T$.
An element $\ovl s_i:=x_i(-1)y_i(1)x_i(-1)$ is in 
$N_G(T)$, which is a representative of 
$s_i\in W=N_G(T)/T$.

\begin{df}
\label{def-gc}
Let $X$ be an affine algebraic variety over $\bbC$, 
$\gamma_i$, $\vep_i, f$ 
$(i\in I)$ rational functions on $X$, and 
$e_i:\bbC^\times\times X\to X$ a unital rational $\bbC^\times$-action
($i\in I$).
A 5-tuple $\chi=(X,\{e_i\}_{i\in I},\{\gamma_i,\}_{i\in I},
\{\vep_i\}_{i\in I},f)$ is a 
$G$ (or $\ge$)-{\it decorated geometric crystal} 
if
\begin{enumerate}
%\bnum
\item
$(\{1\}\times X)\cap dom(e_i)$ 
is open dense in $\{1\}\times X$ for any $i\in I$, where
$dom(e_i)$ is the domain of definition of
$e_i\cl\C^\times\times X\to X$.
%\cmt{changed}
\item
The rational functions  $\{\gamma_i\}_{i\in I}$ satisfy
$\gamma_j(e^c_i(x))=c^{{\bf a}_{ij}}\gamma_j(x)$ for any $i,j\in I$.
\item
The function $f$ satisfies
\begin{equation}
f(e_i^c(x))=f(x)+{(c-1)\vp_i(x)}+{(c^{-1}-1)\vep_i(x)},
\label{f}
\end{equation}
for any $i\in I$ and $x\in X$,  where $\vp_i:=\vep_i\cdot\gamma_i$.
\item
$e_i$ and $e_j$ satisfy the following relations:
\[
 \begin{array}{lll}
&\hspace{-20pt}e^{c_1}_{i}e^{c_2}_{j}
=e^{c_2}_{j}e^{c_1}_{i}&
{\rm if }\,\,{\bf a}_{ij}={\bf a}_{ji}=0,\\
&\hspace{-20pt} e^{c_1}_{i}e^{c_1c_2}_{j}e^{c_2}_{i}
=e^{c_2}_{j}e^{c_1c_2}_{i}e^{c_1}_{j}&
{\rm if }\,\,{\bf a}_{ij}={\bf a}_{ji}=-1,\\
&\hspace{-20pt}
e^{c_1}_{i}e^{c^2_1c_2}_{j}e^{c_1c_2}_{i}e^{c_2}_{j}
=e^{c_2}_{j}e^{c_1c_2}_{i}e^{c^2_1c_2}_{j}e^{c_1}_{i}&
{\rm if }\,\,{\bf a}_{ij}=-2,\,
{\bf a}_{ji}=-1,\\
&\hspace{-20pt}
e^{c_1}_{i}e^{c^3_1c_2}_{j}e^{c^2_1c_2}_{i}
e^{c^3_1c^2_2}_{j}e^{c_1c_2}_{i}e^{c_2}_{j}
=e^{c_2}_{j}e^{c_1c_2}_{i}e^{c^3_1c^2_2}_{j}e^{c^2_1c_2}_{i}
e^{c^3_1c_2}_je^{c_1}_i&
{\rm if }\,\,{\bf a}_{ij}=-3,\,
{\bf a}_{ji}=-1.
\end{array}
\]
\item
The rational functions $\{\vep_i\}_{i\in I}$ satisfy
$\vep_i(e_i^c(x))=c^{-1}\vep_i(x)$ and 
$\vep_i(e_j^c(x))=\vep_i(x)$ if ${\bf a}_{i,j}={\bf a}_{j,i}=0$.
\end{enumerate}
%\ee
\end{df}
We call the function $f$ in (iii) the {\it decoration} of $\chi$ and 
the relations in (iv) are called 
{\it Verma relations}.
If $\chi=(X,\{e_i\},\,\{\gamma_i\},\{\vep_i\})$
satisfies the conditions (i), (ii), (iv) and (v), 
we call $\chi$ a {\it geometric crystal}.
{\sl Remark.}
The definitions of $\vep_i$ and $\vp_i$ are different from the ones in 
e.g., \cite{BK2} since we adopt the definitions following
\cite{KNO,KNO2}. Indeed, if we flip $\vep_i\to \vep^{-1}$ and 
$\vp_i\to \vp^{-1}$, they coincide with ours.
%%%%%%%%%%%%%%%%%%%%%%%%%%%%%%%%%%
\subsection{Characters}\label{char}%4.2

Let $\what U:={\rm Hom}(U,\bbC)$ be the set of additive characters 
of $U$.
The {\it elementary character }$\chi_i\in \what U$ and  
the {\it standard regular character} $\chi^{\rm st}\in \what U$ are  
defined as follows:
\[
\chi_i(x_j(c))=\del_{i,j}\cdot c \q(c\in \bbC,\,\, i\in I),\qq
\chi^{st}=\sum_{i\in I}\chi_i.
\]
We also define the anti-automorphism $\eta:G\to G$ 
by 
\[
 \eta(x_i(c))=x_i(c),\q  \eta(y_i(c))=y_i(c),\q \eta(t)=t^{-1}\q
 (c\in\bbC,\,\, t\in T),
\]
which is called the {\it positive inverse}(\cite{BK2}).

The rational function $f_B$ on $G$ is defined by 
\begin{equation}
f_B(g)=\chi^{st}(\pi^+(w_0^{-1}g))+\chi^{st}(\pi^+(w_0^{-1}\eta(g))),
\label{f_B}
\end{equation}
for $g\in B\ovl w_0 B$, where $\pi^+:B^-U\to U$ is the projection
defined by 
$\pi^+(bu)=u$.

For a split algebraic torus $T$ over $\bbC$, let us denote 
its lattice of (multiplicative )characters(resp. co-characters) by $X^*(T)$ 
(resp. $X_*(T)$). By the usual way, we identify $X^*(T)$
(resp. $X_*(T)$) with the weight lattice $P$ (resp. the dual weight
lattice $P^*$). 

%%%%%%%%%%%%%%%%%%%%%%%%%%%%
\subsection{Positive structure and ultra-discretization}%3.3
\label{subsec-posi}

Let us review the notion positive structure and 
the ultra-discretization.

\begin{df}
Let $T,T'$ be  split algebraic tori over $\bbC$.
\begin{enumerate}
\item
A regular function $f=\sum_{\mu\in X^*(T)}c_\mu\cdot\mu$ on $T$ 
is {\it positive} if all coefficients $c_\mu$ are non-negative 
numbers. A rational function on $T$ is said to be {\it positive} if 
there exist positive regular functions $g,h$ such that 
$f=\frac{g}{h}$ ($h\ne0$).
\item
Let $f:T\to T'$ be a rational map between 
$T$ and $T'$. Then we say that $f$ is {\it positive} if 
for any $\xi\in X^*(T')$ we have that $\xi\circ f$ is positive in the
above sense.
\end{enumerate}
\end{df}
Note that if $f,g$ are positive rational functions on $T$, then 
$f\cdot g$, $f/g$ and $f+g$ are all positive.

\begin{df}
Let $\chi=(X,\{e_i\}_{i\in I},\{{\rm wt}_i\}_{i\in I},
\{\vep_i\}_{i\in I},f)$ be a decorated
geometric crystal, $T'$ an algebraic torus
and $\theta:T'\rightarrow X$ 
a birational map.
The birational map $\theta$ is called 
{\it positive structure} on
$\chi$ if it satisfies:
\begin{enumerate}
\item For any $i\in I$ the rational functions
$\gamma_i\circ \theta, \vep_i\circ \theta, f\circ\theta:T'\rightarrow \bbC$ 
are all positive in the above sense.
\item
For any $i\in I$, the rational map
$e_{i,\theta}:\bbC^\tm \tm T'\rightarrow T'$ defined by
$e_{i,\theta}(c,t)
:=\theta^{-1}\circ e_i^c\circ \theta(t)$
is positive.
\end{enumerate}
\end{df}

Let $v:\bbC(c)\setminus\to\bbZ$ be a map defined by 
$v(f(c)):=\deg(f(c^{-1}))$, which is different from that in e.g.,
\cite{KNO,KNO2,N,N3}. 

Let $f\cl T\rightarrow T'$ be a 
positive rational mapping
of algebraic tori $T$ and 
$T'$.
We define the map $\what f\cl X_*(T)\rightarrow X_*(T')$ by 
\[
\langle\chi,\what f(\xi)\rangle
=v(\chi\circ f\circ \xi),
\]
where $\chi\in X^*(T')$ and $\xi\in X_*(T)$.

Let $\cT_+$ be the category whose objects are algebraic tori over
$\bbC$ and whose morphisms are positive rational maps.
Then, we obtain the ``{\it ultra-discretization}'' functor
\[
\begin{array}{cccc}
{\mathcal UD}:& \cT_+&\longrightarrow &{\mathfrak Set}\\
&T& \mapsto &X_*(T)\\
&(f:T\to T') &\mapsto &(\what f:X_*(T)\to X_*(T')).
\end{array}
\]
Note that this definition of the functor $\mathcal UD$ is
called tropicalization in \cite{BK} and much simpler than the one 
in \cite{BK2}.

Let $\theta:T'\rightarrow X$ be a positive structure on 
a decorated geometric crystal $\chi=(X,\{e_i\}_{i\in I},
\{{\rm wt}_i\}_{i\in I},
\{\vep_i\}_{i\in I},f)$.
Applying the functor ${\mathcal UD}$ 
to positive rational morphisms
$e_{i,\theta}:\bbC^\tm \tm T'\rightarrow T'$ and
$f\circ\Theta,\gamma_i\circ \theta,\vep_i\circ\Theta:T'\ra \bbC$, 
we obtain
\begin{eqnarray*}
\til e_i&:=&{\mathcal UD}(e_{i,\theta}):
\ZZ\tm X_*(T') \rightarrow X_*(T')\\
{\rm wt}_i&:=&{\mathcal UD}(\gamma_i\circ\theta):
X_*(T')\rightarrow \bbZ,\\
\wtil\vep_i&:=&{\mathcal UD}(\vep_i\circ\theta):
X_*(T')\rightarrow \bbZ,\\
\wtil f&:=& {\mathcal UD}(f \circ\theta):
X_*(T')\rightarrow \bbZ.
\end{eqnarray*}
Now, for given positive structure $\theta:T'\rightarrow X$
on a geometric crystal 
$\chi=(X,\{e_i\}_{i\in I},\{{\rm wt}_i\}_{i\in I},
\{\vep_i\}_{i\in I})$, we associate 
the quadruple $(X_*(T'),\{\til e_i\}_{i\in I},
\{{\rm wt}_i\}_{i\in I},\{\wtil\vep_i\}_{i\in I})$
with a free pre-crystal structure (see \cite[2.2]{BK}) 
and denote it by ${\mathcal UD}_{\theta,T'}(\chi)$.
We have the following theorem:

\begin{thm}[\cite{BK,BK2,N}]
For any geometric crystal 
$\chi=(X,\{e_i\}_{i\in I},\{\gamma_i\}_{i\in I},
\{\vep_i\}_{i\in I})$ and positive structure
$\theta:T'\rightarrow X$, the associated pre-crystal 
${\mathcal UD}_{\theta,T'}(\chi)=
(X_*(T'),\{e_i\}_{i\in I},\{{\rm wt}_i\}_{i\in I},
\{\wtil\vep_i\}_{i\in I})$ 
is a Langlands dual Kashiwara's crystal.
\end{thm}
{\sl Remark.}
The definition of $\wtil\vep_i$ is different from the one in 
\cite[6.1.]{BK2} since our definition of $\vep_i$ corresponds to 
$\vep_i^{-1}$ in \cite{BK2}.

For a positive decorated geometric crystal 
$\cX=((X,\{e_i\}_{i\in I},\{\gamma_i\}_{i\in I},
\{\vep_i\}_{i\in I},f),\theta,T')$, set 
\begin{equation}
 \wtil B_{\wtil f}:=\{\wtil x\in X_*(T')(=\bbZ^{\dim(T')})
|\wtil f(\wtil x)\geq0\},
\label{btil}
\end{equation}
and define 
$B_{f,\theta}:=(\wtil B_{\wtil f},\wt_i|_{\wtil B_{\wtil f}},
\vep_i|_{\wtil B_{\wtil f}},e_i|_{\wtil B_{\wtil f}})_{i\in I}$.
\begin{pro}[\cite{BK2}]
For a positive decorated geometric crystal 
$\cX=((X,\{e_i\}_{i\in I},\{\gamma_i\}_{i\in I},
\{\vep_i\}_{i\in I},f),\theta,T')$, the
 quadruple $B_{f,\theta}$ is a normal crystal.
\end{pro}

%%%%%%%%%%%%%%%%%%%%%%%%%%%%%%%%%%%%%%%%%%%
\subsection{Decorated geometric crystal on $\bbB_w$}
\label{BBB}%4.4
For a Weyl group element $w\in W$, define 
$B^-_w:=B^-\cap U\ovl w U$
and set $\bbB_w:=TB^-_w$. 
Let $\gamma_i:\bbB_w\to\bbC$ be the rational function defined by 
\begin{equation}
\gamma_i:\bbB_w\,\,\hookrightarrow \,\,\,
B^-\,\,\mapright{\sim}\,\,T\times U^-\,\,
\mapright{\rm proj}\,\,\, T\,\,\,\mapright{\al_i^\vee}\,\,\,\bbC.
\label{gammai}
\end{equation}

For any $i\in I$, there exists the natural projection 
$pr_i:B^-\to B^-\cap \phi(SL_2)$. Hence, 
for any $x\in \bbB_w$ there exists unique
       $v=\begin{pmatrix}b_{11}&0\\b_{21}&b_{22}\end{pmatrix}
\in SL_2$ such that 
$pr_i(x)=\phi_i(v)$. Using this fact, we define 
the rational function $\vep_i$ on $\bbB_w$ as in \cite{N4}:
\begin{equation}
\vep_i(x)=\frac{b_{22}}{b_{21}}\q(x\in\bbB_w).
\label{vepi}
\end{equation}
The rational $\bbC^\times$-action $e_i$ on $\bbB_w$ is defined by
\begin{equation}
e_i^c(x):=x_i\left((c-1)\vp_i(x)\right)\cdot x\cdot
x_i\left((c^{-1}-1)\vep_i(x)\right)\qq
(c\in\bbC^\times,\,\,x\in \bbB_w),
\label{ei-action}
\end{equation}
if $\vep_i(x)$ is well-defined, that is, $b_{21}\ne0$, 
and define $e_i^c(x)=x$ if $b_{21}=0$.\\
{\sl Remark.} The definition (\ref{vepi}) is different from the one in 
\cite{BK2}. Indeed, if in \eqref{vepi} we take $\vep_i(x)=b_{21}/b_{22}$, 
then it coincides with
the one in \cite{BK2}.
\begin{pro}[\cite{BK2}]
For any $w\in W$,
the 5-tuple $\chi:=(\bbB_w,\{e_i\}_i,\{\gamma_i\}_i,\{\vep_i\}_i,f_B)$
is a decorated geometric crystal, where 
$f_B$ is in (\ref{f_B}), $\gamma_i$ is in (\ref{gammai}), $\vep_i$ is in 
(\ref{vepi}) and $e_i$ is in (\ref{ei-action}).
\end{pro}

\def\ld{\ldots}
For the  longest Weyl group element $w_0\in W$, let 
$\bfii0=i_1\ld i_N$ be one of its reduced expressions and 
define the positive structure on $B^-_{w_0}$ 
$\Theta^-_\bfii0:(\bbC^\times)^N\longrightarrow B^-_{w_0}$ by 
\[
 \Theta^-_\bfii0(c_1,\cd,c_N):=\pmby_{i_1}(c_1)\cd \pmby_{i_N}(c_N),
\]
where $\pmby_i(c)=y_i(c)\al^\vee(c^{-1})$, 
which is different from $Y_i(c)$ in 
\cite{N,N2,KNO,KNO2}. Indeed, $Y_i(c)=\pmby_i(c^{-1})$.
We also define the positive structure on $\bbB_{w_0}$ as
$T\Theta^-_\bfii0:T\times(\bbC^\times)^N\,\,\longrightarrow\,\,\bbB_{w_0}$  
by $T\Theta^-_\bfii0(t,c_1,\cd,c_N)
=t\Theta^-_\bfii0(c_1,\cd,c_N)$.

Now, for this positive structure, we describe the geometric crystal
structure on $\bbB_{w_0}=TB^-_{w_0}$ explicitly.
% See NOTE 16 pp86.
\begin{pro}[\cite{N4}]
The action $e^c_i$ on 
$t\Theta^-_{\bfii0}(c_1,\cd,c_N)$ is given by
\[
e_i^c(t\Theta^-_{\bfii0}(c_1,\cd,c_N))
=t\Theta^-_{\bfii0}(c'_1,\cd,c'_N)
\]
where
\begin{equation}
c'_j\seteq 
c_j\cdot \frac{\displaystyle \sum_{1\leq m< j,\,i_m=i}
{c\cdot c_1^{{\bf a}_{i_1,i}}\cd c_{m-1}^{{\bf a}_{i_{m-1},i}}c_m}
+\sum_{j\leq m\leq N,\,i_m=i} 
{c_1^{{\bf a}_{i_1,i}}\cd c_{m-1}^{{\bf a}_{i_{m-1},i}}c_m}}
{\displaystyle \sum_{1\leq m\leq j,\,i_m=i} 
{c\cdot c_1^{{\bf a}_{i_1,i}}\cd c_{m-1}^{{\bf a}_{i_{m-1},i}}c_m}+
\mathop\sum_{j< m\leq N,\,i_m=i}  
{c_1^{{\bf a}_{i_1,i}}\cd c_{m-1}^{{\bf a}_{i_{m-1},i}}c_m}}.
\label{eici}
\end{equation}
The explicit forms of 
rational functions $\vep_i$ and $\gamma_i$ are:
\begin{equation}
 \vep_i(t\Theta^-_{\bfii0}(c))=
\left(\sum_{1\leq m\leq N,\,i_m=i} \frac{1}
{c_mc_{m+1}^{{\bf a}_{i_{m+1},i}}\cd c_{N}^{{\bf a}_{i_{N},i}}}\right)^{-1},\q
\gamma_i(t\Theta^-_{\bfii0}(c))
=\frac{\al_i(t)}{c_1^{{\bf a}_{i_1,i}}\cd c_N^{{\bf a}_{i_N,i}}}.
\label{th-vep-gamma}
\end{equation}
\end{pro}

%%%%%%% Sect.5 %%%%%%%%%%%%%%%%
\renewcommand{\thesection}{\arabic{section}}
\section{Monomial Realization of Crystals}
\label{sect-mono}
\setcounter{equation}{0}
\renewcommand{\theequation}{\thesection.\arabic{equation}}

\subsection{Definitions of Monomial Realization of Crystals}
\label{mono-subsec}
Following \cite{K7,Nj}, 
we shall introduce the monomial realization of crystals.
For doubly-indexed variables $\{Y_{m,i}|i\in I, m\in\bbZ.\}$, 
define the set of monomials
\[
 \cY:=\{Y=\prod_{m\in\bbZ,i\in I}
Y_{m,i}^{l_{m,i}}|l_{m,i}\in \bbZ\setminus\{0\}\text{
 except for finitely many }(m,i)\}.
\]
Fix a set of integers $p=(p_{i,j})_{i,j\in I,i\ne j}$ such that 
$p_{i,j}+p_{j,i}=1$, which we call a {\it sign}. 
Take a sign $p:=(p_{i,j})_{i,j\in I,i\ne j}$ and a
Cartan matrix $({\bf a}_{i,j})_{i,j\in I}$. For $m\in\bbZ,\,i\in I$ 
define the monomial 
\[
 A_{m,i}=Y_{m,i}Y_{m+1,i}\prod_{j\ne i}Y_{m+p_{j,i},j}^{{\bf a}_{j,i}}.
\]
Here, when we emphasize the sign $p$, we shall denote the
monomial
$A_{m,i}$ by $A^{(p)}_{m,i}$.
For any cyclic sequence of the indices
$\io=\cd (i_1i_2\cd i_n)(i_1i_2\cd i_n)\cd$
s.t. $\{i_1,\cd,i_n\}=I$,
we can associate the following sign $(p_{i,j})$ by:
\begin{equation}
 p_{i_a,i_b}=\begin{cases}
1&a<b,\\0&a>b.\end{cases}
\label{iaib}
\end{equation}
For example, if we take $\io=\cd (213)(213)\cd$, then we have
$p_{2,1}=p_{1,3}=p_{2,3}=1$ and $p_{1,2}=p_{3,1}=p_{3,2}=0$.
Thus, we can identify a cyclic sequence $\cd(i_1\cd i_n)(i_1\cd i_n)\cd$ 
with such $(p_{i,j})$.

For a monomial $Y=\prod_{m,i}Y_{m,i}^{l_{m,i}}$, 
set 
\begin{eqnarray*}
&&\hspace{-20pt}wt(Y)=\sum_{i,m} l_{m,i}\Lm_i, \,\,
\vp_i(Y)=\operatorname{max}_{m\in\bbZ}\{\sum_{k\leq m}l_{k,i}\},\,\,
\vep_i(Y)=\vp_i(Y)-wt(Y)(h_i)
=\max_{m\in\bbZ}\{-\sum_{k>m}l_{i,k}\},\\
&&\hspace{-20pt}\fit(Y)=\begin{cases}
A_{n_f^{(i)}(Y),i}^{-1}\cdot Y&\text{ if }\vp_i(Y)>0,\\
0&\text{ if }\vp_i(Y)=0,
\end{cases}\q\q
\eit(Y)=\begin{cases}
A_{n_e^{(i)}(Y),i}\cdot Y&\text{ if }\vep_i(Y)>0,\\
0&\text{ if }\vep_i(Y)=0,
\end{cases}\\
&&\text{where }n_f^{(i)}(Y)=\min\{n|\vp_i(Y)=\sum_{k\leq n}l_{k,i}\},\q
n_e^{(i)}(Y)=\max\{n|\vp_i(Y)=\sum_{k\leq n}l_{k,i}\}.
\end{eqnarray*}

\begin{thm}[\cite{K7,Nj}]
\begin{enumerate}
\item
In the above setting, $\cY$ is a crystal, which is denoted by $\cY(p)$.
\item
If $Y\in\cY(p)$ satisfies $\vep_i(Y)=0$ (resp. $\vp_i(Y)=0$) for any $i\in I$, 
then the connected component containing $Y$ is isomorphic to
$B(wt(Y))$ (resp. $B(w_0\wt(Y))$), 
where we call such monomial $Y$ a highest (resp. lowest) monomial.
\end{enumerate}
\end{thm}
By the above crystal structure of monomials, we know that 
for any $k\in \bbZ,\,\,i\in I$
the monomial $Y_{k,i}$ is a highest monomial with the weight $\Lm_i$. 
Thus, we can define the embedding of crystal 
$\mpp(m,k,i)$ $(i\in I,k\in\bbZ)$ by 
\begin{eqnarray}
 \mpp(m,k,i):B(\Lm_i)&\hookrightarrow& \cY(p)
\label{embed-mki}\\
u_{\Lm_i}&\mapsto& Y_{k,i} \nn
\end{eqnarray}

\subsection{Duality on Monomial Realizations}

Let $i_1,\cd,i_n$ be the indices satisfying $\{i_1,\cd,i_n\}=I$.
For the cyclic sequence 
${\bf i}=\cd(i_1i_2\cd i_n)(i_1i_2\cd i_n)\cd$
let $p=(p_{i,j})$ be the associated sign as in\eqref{iaib}.

And for the opposite cyclic sequence
${\bf i}^{-1}=\cd(i_ni_{n-1}\cd i_1)(i_ni_{n-1}\cd i_1)\cd$, 
let ${}^tp=p':=(p'_{i,j})=(p_{j,i})$  be the associated sign.
For a monomial $Y=\prod_{m,i}Y_{i,m}^{l_{i,m}}\in\cY(p)$, define the 
map ${\,}^{a-}:\cY(p)\to \cY({}^tp)$ is defined by 
$Y_{i,m}\mapsto Y_{i,a-m}^{-1}$ for $a\in \bbZ$.

Indeed, the following lemma is derived by direct calculations:
\begin{pro}
For a monomial $Y\in\cY(p)$, we have ${}^a\ovl Y\in \cY({^tp})$ and 
\begin{eqnarray}
&&\vp_i({}^a{\ovl Y})=\vep_i(Y),\q
\vep_i({}^a\ovl Y)=\vp_i(Y),\q wt({}^a\ovl Y)=-wt(Y),
\label{bar-ep}\\
&&{}^a\ovl{\fit(Y)}=\eit({}^a\ovl Y), \qq {}^a\ovl{\eit(Y)}=\fit({}^a\ovl Y).
\label{bar-eit}
\end{eqnarray}
\end{pro}
{\sl Proof.}
The last formula in \eqref{bar-ep} is trivial by the definition.
For an arbitrary monomial $Y=\prod_{i,m}Y_{i,m}^{l_{i,m}}\in\cY(p)$, 
we have 
${}^a\ovl Y=\prod_{i,m}Y_{i,a-m}^{-l_{i,m}}=\prod_{i,m}Y_{i,m}^{-l_{i,a-m}}.$
Thus, one gets 
\[
\vp_i({}^a\ovl Y)=
\max_{m\in\bbZ}\{\sum_{k\leq m}-l_{i,a-k}\}
=\max_{m\in\bbZ}\{\sum_{a-k\leq m}-l_{i,k}\}
=\max_{m\in\bbZ}\{\sum_{k\geq m}-l_{i,k}\}
=\vep_i(Y).
\]
The second one is obtained similarly.

As for \eqref{bar-eit},  we shall show the first one.
Let us see the following lemma:
\begin{lem}
We have 
\begin{eqnarray}
&&{}^a{\ovl{A^{(p)}_{i,m}(Y)}}=A^{({^tp})}_{i,a-m-1}(Y),
\label{aim-dual}\\
&&n_e^{(i)}({}^a\ovl Y)=a-n_f^{(i)}(Y)-1,\q 
n_f^{(i)}({}^a\ovl Y)=a-n_e^{(i)}(Y)-1,
\label{ne-nf-dual}
\end{eqnarray}
\end{lem}
{\sl Proof.}
The first formula is trivial from the definition of $A_{i,m}$. 
As for \eqref{ne-nf-dual}, we can see 
\begin{eqnarray*}
n_e^{(i)}({}^a\ovl Y)&&=
\max_{n}\{\vp_i({}^a\ovl Y)=\sum_{k\leq n}-l_{i,a-k}\}
=\max_{n}\{\vep_i(Y)=\sum_{k\leq n}-l_{i,a-k}\}\\
&&=\max_{n}\{\wt_i(Y)+\vep_i(Y)=\wt_i(Y)+\sum_{k\leq n}-l_{i,a-k}\}
=\max_n\{\vp_i(Y)=\sum_{k>n}l_{i,a-k}\}\\
&&=\max_n\{\vp_i(Y)=\sum_{k<a-n}l_{i,k}\}
=a-\min_n\{\vp_i(Y)=\sum_{k<n}l_{i,k}\}=a-(n_f^{(i)}(Y)+1).
\end{eqnarray*}
The second one is shown similarly.\qed\\
By this lemma, we can easily get \eqref{bar-eit}.\qed

For $Y\in \cY(p)$, let us denote the connected component containing $Y$
by $B(Y)$. Then, by the above proposition we obtain:
\begin{thm}\label{thm-dual}
For any $Y\in\cY(p)$ and any $a\in\bbZ$, 
the set ${}^a\ovl{B(Y)}$ is equipped with the crystal 
structure associated with ${}^tp$ and there exists the isomorphism
of crystals :
\begin{equation}
{}^a\ovl{B(Y)}\mapright{\sim}B({}^a\ovl Y)(\subset \cY({^tp}))
\qq({}^a\ovl Y\mapsto {}^a\ovl Y).
\end{equation}
\end{thm}
Indeed, we find that if $Y$ is a highest(resp. lowest) monomial in 
$\cY(p)$, then ${}^a\ovl Y$ is a lowest(resp. highest) monomial 
in $\cY({^tp})$.

%%%%%%% Sect 6. %%%%%%%%%%%%%%
\renewcommand{\thesection}{\arabic{section}}
\section{Explicit form of the decoration $f_B$ for Classical Groups}
\setcounter{equation}{0}
\renewcommand{\theequation}{\thesection.\arabic{equation}}

%%%%%%%%%%%%%%%%%%%%%%%%%%%%%%%%%%%%%%%%%%%%%%%%%%%%%%
\subsection{Generalized Minors and the function $f_B$}%6.1

For this subsection, see \cite{BFZ,BZ,BZ2}.
Let $G$ be a simply connected simple algebraic groups over $\bbC$ and 
$T\subset G$ a maximal torus. 
Let  $X^*(T):=\Hom(T,\bbC^\times)$ and $X_*(T):=\Hom(\bbC^\times,T)$ be
the lattice of characters and co-characters respectively.
We identify $P$ (resp. $P^*$) with $X^*(T)$ 
(resp. $X_*(T)$) as above.
\begin{df}
For a dominant weight $\mu\in P_+$, the
{\it principal minor} $\Del_\mu:G\to\bbC$ is defined as
\[
 \Del_\mu(u^-tu^+):=\mu(t)\q(u^\pm\in U^\pm,\,\,t\in T).
\]
Let $\gamma,\del\in P$ be extremal weights such that 
$\gamma=u\mu$ and $\del=v\mu$ for some $u,v\in W$. 
Then the {\it generalized minor} $\Del_{\gamma,\del}$ is defined
by
\[
 \Del_{\gamma,\del}(g):=\Del_\mu(\ovl u^{-1}g\ovl v)
\q(g\in G),
\]
which is a regular function on $G$.
\end{df}
\begin{lem}[\cite{BK2}]
Suppose that $G$ is simply connected.
\begin{enumerate}
\item
For $u\in U$ and $i\in I$, we have 
$\Del_{\mu,\mu}(u)=1$ and $\chi_i(u)=\Del_{\Lm_i,s_i\Lm_i}(u)$,
where $\Lm_i$ be the $i$th fundamental weight.
\item
Define the map $\pi^+:B^-\cdot U\to U$ by $\pi^+(bu)=u$ for 
$b\in B^-$ and $u\in U$. For any $g\in G$, we have 
\begin{equation}
 \chi_i(\pi^+(g))=\frac{\Del_{\Lm_i,s_i\Lm_i}(g)}
{\Del_{\Lm_i,\Lm_i}(g)}.
\end{equation}
\end{enumerate}
\end{lem}
\begin{pro}[\cite{BK2}]
The function $f_B$ in (\ref{f_B}) is described as follows:
\begin{equation}
 f_B(g)=\sum_i\frac{\Del_{w_0\Lm_i,s_i\Lm_i}(g)
+\Del_{w_0s_i\Lm_i,\Lm_i}(g)}
{\Del_{w_0\Lm_i,\Lm_i}(g)}.
\label{chi-min}
\end{equation}
\end{pro}
%%%%%%%%%%%%%%%%%%%%%%%%%%%%%%%%%%%%%%
\subsection{Bilinear Forms}%4.2

Let $\omega:\ge\to\ge$ be the anti involution 
\[
\omega(e_i)=f_i,\q
\omega(f_i)=e_i\,\q\omega(h)=h,
\] and extend it to $G$ by setting
$\omega(x_i(c))=y_i(c)$, $\omega(y_i(c))=x_i(c)$ and $\omega(t)=t$
$(t\in T)$.

There exists a $\ge$(or $G$)-invariant bilinear form on the
finite-dimensional  irreducible
$\ge$-module $V(\lm)$ such that 
\[
 \lan au,v\ran=\lan u,\omega(a)v\ran,
\q\q(u,v\in V(\lm),\,\, a\in \ge(\text{or }G)).
\]
For $g\in G$, 
we have the following simple fact:
\[
 \Del_{\Lm_i}(g)=\lan gu_{\Lm_i},u_{\Lm_i}\ran,
\]
where $u_{\Lm_i}$ is a properly normalized highest weight vector in
$V(\Lm_i)$. Hence, for $w,w'\in W$ we have
\begin{equation}
 \Del_{w\Lm_i,w'\Lm_i}(g)=
\Del_{\Lm_i}({\ovl w}^{-1}g\ovl w')=
\lan {\ovl w}^{-1}g\ovl w'\cdot u_{\Lm_i},u_{\Lm_i}\ran
=\lan g\ovl w'\cdot u_{\Lm_i}\, ,\, \ovl{w}\cdot u_{\Lm_i}\ran,
\label{minor-bilin}
\end{equation}
where note that $\omega(\ovl s_i^{\pm})=\ovl s_i^{\mp}$.

\begin{pro}\label{form-alt}
Let ${\bf i}=i_1\cd i_N$ be a reduced word for the longest Weyl
group element $w_0$. 
For $t\Theta_{\bf i}^-(c)\in \bbB_{w_0}=T\cdot B^-_{w_0}$, 
we get the following formula.
\begin{equation}
 f_B(t\Theta_{\bf i}^-(c))
=\sum_i\Del_{w_0\Lm_i,s_i\Lm_i}(\Theta_{\bf i}^-(c))
+\al_i(t)\Del_{w_0s_i\Lm_i,\Lm_i}(\Theta_{\bf i}^-(c)).
\label{fb-th}
\end{equation}
\end{pro}
%Note17 pp
{\sl Proof.}
We shall show 
\begin{equation}
\Del_{\ovl w_0\Lm_i,\Lm_i}(\Theta_{\bf i}^{-}(c))=1.
\label{del-1}
\end{equation}
Since $\Theta_{\bf i}^{-}(c)\in U\ovl w_0 U$, we have 
$\ovl w_0^{-1}\Theta_{\bf i}^{-}(c)\in \ovl w_0^{-1}U\ovl w_0U=U^-\cdot
U$.
So, there exist  $u_1\in U^-$ and $u_2\in U$ such that 
$\ovl w_0^{-1}\Theta_{\bf i}^{-}(c)=u_1u_2$.
Thus, it follows from \eqref{minor-bilin} that 
\[
\Del_{\ovl w_0\Lm_i,\Lm_i}(\Theta_{\bf i}^{-}(c))=
\lan \ovl w_0^{-1}\Theta_{\bf i}^-(c)u_{\Lm_i},u_{\Lm_i}\ran
=\lan u_1u_2u_{\Lm_i},u_{\Lm_i}\ran
=\lan u_2u_{\Lm_i},\omega(u_1)u_{\Lm_i}\ran=\lan
u_{\Lm_i},u_{\Lm_i}\ran=1,
\]
since $\omega(u_1)\in U$.
The following is evident from \eqref{minor-bilin}
\begin{eqnarray}
\Del_{w_0s_i\Lm_i,\Lm_i}(t\Theta_{\bf i}^{-}(c))&=&
\lan t\Theta_{\bf i}^{-}(c)u_{\Lm_i},\ovl w_0\ovl s_iu_{\Lm_i\ran}
=\lan \Theta_{\bf i}^{-}(c)u_{\Lm_i},t\cdot\ovl w_0\ovl
 s_iu_{\Lm_i\ran}
\label{deldel}\\
&=&\Del_{w_0s_i\Lm_i,\Lm_i}(\Theta_{\bf i}^{-}(c))\cdot 
w_0s_i\Lm_i(t).\nn
\end{eqnarray}
Since $w_0s_i\Lm_i/w_0\Lm_i=\al_i$ and 
$\Del_{w_0\Lm_i,\Lm_i}(t\Theta_{\bf i}^{-}(c))=
\Del_{w_0\Lm_i,\Lm_i}(\Theta_{\bf i}^{-}(c))\cdot w_0\Lm_i(t)$, 
we obtain the desired result.\qed
%See note 17 pp32-33.

{\sl Remark.}
Note that by virtue of Proposition \ref{form-alt} to get the explicit
form of $f_B(t\Theta_{\bf i}^{-}(c))$ it is sufficient to 
know those of $\Del_{w_0\Lm_i,s_i\Lm_i}(\Theta_{\bf i}^{-}(c))$ and 
 $\Del_{w_0s_i\Lm_i,\Lm_i}(\Theta_{\bf i}^{-}(c))$. 
%%%%%%%%%%%%%%%%%%%%%%%%%%%%%%%%%%%
\subsection{Explicit form of 
$f_B(t\Theta_{\bf i}^-(c))$ for $A_n$}
\label{fb-An}
%4.3

For all classical cases $A_n,B_n,C_n$ and $D_n$,
we fix the reduced longest word $\bfii0$ as follows:
\begin{equation}
\bfii0
=\begin{cases}\underbrace{1,2,\cd,n}_{},
\underbrace{1,2,\cd,n-1}_{},\cd,\underbrace{1,2,3}_{},1,2,1&\text{ for
  }A_n,\\
(1,2,\cd,n-1,n)^n&\text{ for }B_n,\,\,C_n,\\
(1,2,\cd,n-1,n)^{n-1}&\text{ for }D_n.
\end{cases}
\end{equation}

To obtain the explicit form of $f_B(t\Theta_{\bf i_0}^-(c))$ for type
$A_n$, 
by the above remark it suffices to know 
$\Del_{w_0\Lm_i,s_i\Lm_i}(\Theta_{\bf i_0}^-(c))$ and 
$\Del_{w_0s_i\Lm_i,\Lm_i}(\Theta_{\bf i_0}^-(c))$ for 
\[
c=(c^{(i)}_j|i+j\leq n+1)=
(c^{(1)}_1,c^{(1)}_2,\cd,c^{(1)}_n,c^{(2)}_1,c^{(2)}_2,\cd,c^{(2)}_{n-1},\cd
c^{(n-1)}_1,c^{(n-1)}_2,c^{(n)}_1)\in (\bbC^\times)^N.
\]
The following result for type $A_n$ is given in \cite{N4}:
\begin{thm}[\cite{N4}]
\label{thm-a}
For $c\in (\bbC^\times)^N$ as above, 
we have the following explicit forms:
\begin{eqnarray}
&&\Del_{w_0\Lm_j,s_j\Lm_j}(\Theta_\bfii0^-(c))=
c^{(n-j+1)}_1+\frac{c^{(n-j+1)}_2}{c^{(n-j+2)}_1}
+\frac{c^{(n-j+1)}_3}{c^{(n-j+2)}_2}
+\cd+\frac{c^{(n-j+1)}_j}{c^{(n-j+2)}_{j-1}}, 
\label{del1}\\
&&\Del_{w_0s_j\Lm_j,\Lm_j}(\Theta_\bfii0^-(c))=
\frac{1}{c^{(j)}_1}+
\frac{c^{(j-1)}_1}{c^{(j-1)}_2}+\frac{c^{(j-2)}_2}{c^{(j-2)}_3}
+\cd+\frac{c^{(1)}_{j-1}}{c^{(1)}_j},\,\,
(j\in I).
\label{del2}
\end{eqnarray}
\end{thm}

As we mentioned in the introduction, we shall see
the relations between the decoration and the 
monomial realizations of crystals of 
type $A_n$ explicitly following \cite{N4}. 
For type $A_n$ take $(p_{i,j})_{i,j\in I,i\ne j}$
such that $p_{i,j}=1$ for $i<j$, $p_{i,j}=0$
for $i>j$, which corresponds to the cyclic sequence
${\mathbf i}=(12\cd n)(12\cd n)\cd$.
Then we obtain 
\begin{pro}[\cite{N4}]
\label{mono-cry}
The crystal containing the monomial $Y_{n-i+1,1}$ (resp. $Y_{i,1}^{-1}$)
is isomorphic to $B(\Lm_1)$ (resp. $B(\Lm_n)$)
and all basis vectors are given by
\begin{eqnarray*}
&&
\til f_k\cd \til f_2\til
f_1(Y_{n-i+1,1})=\frac{Y_{n-i+1,k+1}}{Y_{n-i+2,k}}\in B(\Lm_1),\\
&&\til e_k\cd\til e_2\til e_1(Y_{i,1}^{-1})
=\frac{Y_{i-k,k}}{Y_{i-k,k+1}}\in B(\Lm_n)\q
(k=1,\cd,n).
\end{eqnarray*}
\end{pro}

Applying this results to Theorem \ref{thm-a} and 
changing the variable $Y_{m,l}$ to $c^{(m)}_l$, we find:
\begin{pro}[\cite{N4}]
\label{del-mono}
For $j=1,\cd,n$  we have
\begin{eqnarray*}
&&\chi_j(\pi^+(w_0^{-1}t\Theta_\bfii0^-(c)))=
\Del_{w_0\Lm_j,s_j\Lm_j}(\Theta_\bfii0^-(c))=
\sum_{k=0}^{j-1}\til f_k\cd \til f_2\til f_1(c^{(n-j+1)}_1),\\
&&\chi_j(\pi^+(w_0^{-1}\eta(t\Theta_\bfii0^-(c))))=
\al_j(t)\Del_{w_0s_j\Lm_j,\Lm_j}(\Theta_\bfii0^-(c))=
\al_j(t)
\sum_{k=0}^{j-1}\til e_k\cd\til e_2\til e_1({c^{(j)}_1}^{-1}).
\end{eqnarray*}
\end{pro}
Note that 
$\{\til f_k\cd \til f_2\til f_1(c_{n-i+1,1})|0\leq k<i\}
=B(\Lm_1){\tiny s_{k-1}\cd s_2s_1}$
is the Demazure crystal associated with the Weyl group element 
$s_{k-1}\cd s_2s_1$ (\cite{K3}).

Observing Proposition \ref{del-mono}, we present the following
conjecture:
\begin{conj}[\cite{N4}]
\label{conj1}
There exists certain reduced longest word 
$\bfi=(i_1,\cd, i_N)$ and a  sign $p=(p_{i,j})_{i\ne j}$ 
such that for any $i\in I$, 
there exist Demazure crystal $B^-_{w}(i)\subset B(\Lm_k)$,
Demazure crystal $B^+_{w'}(i)\subset B(\Lm_j)$ and 
positive integers $\{a_b,a_{b'}|b\in B^-_w,b'\in B^+_{w'}\}$ satisfying
\begin{eqnarray*}
&&\chi_i(\pi^+(w_0^{-1}t\Theta_\bfi^-(c)))=
\Del_{w_0\Lm_i,s_i\Lm_i}(\Theta_\bfi^-(c))=
\sum_{b\in B^-_w(i)}a_b
m_b(c),\\
&&\chi_i(\pi^+(w_0^{-1}\eta(t\Theta_\bfi^-(c))))=
\al_i(t)\Del_{w_0s_i\Lm_i,\Lm_i}(\Theta_\bfi^-(c))=
\al_i(t)\sum_{b'\in B^+_{w'}(i)}a_{b'}m_{b'}(c),
\end{eqnarray*}
where $m_b(c)\in\cY(p)$ is the monomial corresponding to $b\in B(\Lm_k)$
associated with $p=(p_{i,j})_{i\ne j}$.
\end{conj}

We would see the answers to this
conjecture for other type of Lie algebras in the 
subsequent sections.

%%%%%%%% Sect. 7%%%%%%%%%%%%%%%%%%%%%%%%%%%%%%%%%%%
\renewcommand{\thesection}{\arabic{section}}
\section{Explicit form of 
$f_B(t\Theta_{\bfii0}^-(c))$ for $C_n$}\label{cn}
\setcounter{equation}{0}
\renewcommand{\theequation}{\thesection.\arabic{equation}}

\subsection{Main theorems}
In this section we see the results for type $C_n$.
\begin{thm}\label{thm-c-1}
In the case $C_n$,  for $k=1,\cd,n$, $\bfii0=(12\cd n)^n$ and 
$c=(\ci(i,j))_{1\leq i,j\leq n}\\
=(\ci(1,1),\ci(2,1),\cd,\ci(n-1,n),\ci(n,n))\in
(\bbC^\times)^{n^2}$
we have
\begin{equation}
\Del_{w_0\Lm_k,s_k\Lm_k}(\Theta^-_\bfii0(c))
=\ci(1,k)+\frac{\ci(2,k)}{\ci(1,k+1)}+\cd+
\frac{\ci(n,k)}{\ci(n-1,k+1)}+
\frac{\ci(n-1,k+1)}{\ci(n,k+1)}+\frac{\ci(n-2,k+2)}{\ci(n-1,k+2)}+\cd+
\frac{\ci(k,n)}{\ci(k+1,n)},
\end{equation}
where note that $\Del_{w_0\Lm_n,s_n\Lm_n}(\Theta^-_\bfii0(c))=\ci(n,n)$.
\end{thm}
We also get the following theorem.
\begin{thm}\label{thm-c-2}
\begin{enumerate}
\item
Let $k$ be in $\{1,2,\cd,n-1\}$. Then we have
\begin{eqnarray}
\Del_{w_0s_k\Lm_k,\Lm_k}(\Theta^-_\bfii0(c))
=\frac{1}{\ci(1,k)}+\sum_{j=1}^{k-1}\frac{\ci(k-j,j)}{\ci(k-j+1,j)}.
\label{del-c-2}
\end{eqnarray}
\item
For variables $(\ci(i,j))_{1\leq i,j\leq n}$, set 
$\Ci(i,j)=\frac{\ci(i,n-j)}{\ci(i-1,n-j+1)}$ and 
$\bCi(i,j)=\frac{\ci(i-1,n-j)}{\ci(i,n-j)}$. Then we have
\begin{equation}
\Del_{w_0s_n\Lm_n,\Lm_n}(\Theta^-_\bfii0(c))
=\sum_{(*)}
{\Ci(u_1,1)}{\Ci(u_2,2)}\cd {\Ci(u_m,m)}
{\bCi(l_1,l_1-1)}\cdot{\bCi(l_2,l_2-2)}\cd 
{\bCi(l_{k},l_{k}-k)},
\label{del-c-3}
\end{equation}
%note 16 pp 37 red pen's
where 
%$\di(i,j)=\frac{\ci(i,j)}{\ci(i-1,j)}$, 
%$\dbi(i,j)=\frac{\ci(i-1,j-1)}{\ci(i,j)}$ and 
$(*)$ is the conditions:
$k+m=n$, $0\leq m<n$,\\
$1\leq l_1< l_2<\cd<l_k\leq n,$ and 
$1\leq u_1< u_2<\cd<u_m\leq n.$
\end{enumerate}
\end{thm}

%%%%%%%%%%%%%%%%%%%%%%%%%%%%%%%%%%
\subsection{Proof of Theorem \ref{thm-c-1}}
\label{proof-c-1}

On the module $V(\Lm_1)$ we can write 
${\pmb x}_i(c):=\al_i^\vee(c^{-1})x_i(c)=c^{-h_i}(1+c\cdot e_i)$ 
and ${\pmb y}_i(c):=y_i(c)\al_i^\vee(c^{-1})=(1+c\cdot f_i)c^{-h_i}$ 
since $f_i^2=e_i^2=0$ on $V(\Lm_1)$.

We also have $\omega({\pmb y}_i(c))=\al_i^\vee(c^{-1})x_i(c)=
{\pmb x}_i(c)$ and define $\xxi(\ovl i,p,j)=\xxi(\ovl i,p,j)(c^{[1:p]})$ 
and $\xxi(i,p,j)=\xxi(i,p,j)(c^{[1:p]})$ for 
$p,j \in I$, $c=(\ci(i,j))_{1\leq i,j\leq n}\in(\bbC^{\times})^{n^2}$ by 
\begin{eqnarray*}
&&X^{(p)}X^{({p-1})}\cd X^{(1)}v_{i}
=\sum_{j=1}^n{\xxi(i,p,j)}v_j
+\sum_{j=1}^n{\xxi(i,p,\ovl j)}v_{\ovl j}\in V(\Lm_1)
\q(i=1,2,\cd,n),\\
&&X^{(p)}X^{({p-1})}\cd X^{(1)}v_{\ovl i}
=\sum_{j=1}^n{\xxi(\ovl i,p,j)}v_j
+\sum_{j=1}^n{\xxi(\ovl i,p,\ovl j)}v_{\ovl j}\in V(\Lm_1)
\q(i=1,2,\cd,n),
\end{eqnarray*}
where $c^{[1:p]}=(\ci(1,1),\ci(2,1),\cd,\ci(n-1,p),\ci(n,p))$ and 
$X^{(p)}={\pmb x}_n(\ci(n,p)){\pmb x}_{n-1}(\ci(n-1,p))\cd 
{\pmb x}_1(\ci(1,p))$.

By (\ref{minor-bilin}) and $\omega(\Theta_{{\bf i}_0}(c))=X^{(n)}\cd X^{(1)}$, 
we have 
\begin{equation}
 \Del_{w_0\Lm_i,s_i\Lm_i}(\Theta_\bfii0(c))=\lan \ovl s_i\cdot
 u_{\Lm_i}\,,\,
X^{(n)}\cd X^{(1)}v_{\Lm_i}\ran, 
\label{del-x}
\end{equation}
and then e.g.,$\xxi(\ovl 1,n,2)=\Del_{w_0\Lm_1,s_1\Lm_1}(\Theta_\bfii0(c))$.

Here to describe $\xxi(\ovl i,p,j)$ explicitly
let us introduce some combinatorial object {\it segments}
as follows.
For $1\leq p,k\leq n$ 
set $L:=p-n+k$ and $S:=n-k+1$. 
For  $r=0,1,\cd,n-p$, set 
\[
 {\mathcal M}^{(p)}_k[r]
:=\{M=\{m_2,m_3,\cd,m_L\}|2+r\leq m_2<\cd m_L\leq p+r\}.
\]
We usually denote ${\mathcal M}^{(p)}_k[0]$ by 
${\mathcal M}^{(p)}_k$.
Define the {\it segments} of $M\in {\mathcal M}^{(p)}_k[r]$ as 
$M=M_1\sqcup\cd\sqcup M_S$ where 
each segment $M_j$ is a consecutive subsequence of $M$ or an empty set
such that $\min(M_b)=\max(M_a)+(b-a+1)$ for $a<b$ and non-empty $M_a,M_b$.
Note that for $M=\{m_2,\cd,m_L\}\in {\mathcal M}^{(p)}_k[r]$ and
$s\in\{-r,-r+1,\cd,n-p-r\}$
we find that $M[s]:=\{m_2+s,\cd,m_L+r\}$ is an element in 
${\mathcal M}^{(p)}_k[r+s]$.
% See Note17, pp24
\begin{ex}
For $n=p=6$, $k=4$ and $r=0$ we have $L=4$ and $S=3$. 

$M=\{2,3,5\}$ $\Longrightarrow$ 
$M_1=\{2,3\}$, $M_2=\{5\}$, $M_3=\emptyset$.

$M=\{2,3,6\}$ $\Longrightarrow$ 
$M_1=\{2,3\}$, $M_2=\emptyset$, $M_3=\{6\}$, 

$M=\{2,4,6\}$ $\Longrightarrow$ 
$M_1=\{2\}$, $M_2=\{4\}$, $M_3=\{6\}$, 

$M=\{3,4,6\}$ $\Longrightarrow$ 
$M_1=\emptyset$, $M_2=\{3,4\}$, $M_3=\{6\}$, 
\end{ex}
For $m\in M=M_1\sqcup\cd\sqcup M_S\in{\mathcal M}^{(p)}_k[r]$, define 
$n(m):=n-j+1$ if $m\in M_j$.
For $M=M_1\sqcup\cd\sqcup M_S\in{\mathcal M}^{(p)}_k[r]$, write
$M_1=\{2+r,3+r,\cd,a\}$ where note that non-empty $M_1$ has to include $2+r$.
For $i_{2+r},\cd,i_a$ satisfying 
$i-1\leq i_{2+r}\leq i_{3+r}\leq\cd\leq i_a\leq n$, define
\[
C_{i_{2+r},i_{3+r},\cd,i_a}^M:=\frac{\ci(i_{2+r}+1-2\ep_{i_{2+r}},
1+r+\ep_{i_{2+r}})\cd
\ci(i_a+1-2\ep_{i_a},a+\ep_{i_a}-1)}
{\ci(i_{2+r},2+r)\cd \ci(i_a,a)},\qq
D^M:=\prod_{m\in M\setminus M_1}\frac{\ci(n(m)-1,m)}{\ci(n(m),m)},
\]
where $\ep_i=\del_{i,n}$ and $C_{i_{2+r},i_{3+r},\cd,i_a}^M=1$ 
(resp. $D^M=1$) if $M_1=\emptyset$ (resp. $M\setminus M_1=\emptyset$).
Here, 
for $(\ci(k,l))$ we set $c^{(l)}:=(\ci(1,l),\cd,\ci(n,l))$ and 
$c^{[a:b]}:=(c^{(a)},c^{(a+1)},\cd,c^{(b)})$. Then, $c=c^{[1:n]}$.
Indeed, for $M\in\cM^{(p)}_k[r]$ the monomial $C^M\cdot D^M$ 
depends on $c^{[2+r:p+r]}$ and then for any $s=1,2,\cd,n-p+2$
and $q=0,1,\cd,r$ we have 
\begin{equation}
C_{i_{2+r},i_{3+r},\cd,i_a}^M\cdot D^M=
C_{i_{2+r},i_{3+r},\cd,i_a}^M\cdot D^M(c^{[s:s+p-2]})
=C_{i_{2+r},i_{3+r},\cd,i_a}^{M[-q]}\cdot D^{M[-q]}(c^{[s:s+p-2]}),
\label{shift-cd}
\end{equation}
where $M[-q]\in\cM^{(p)}_k[r-q]$.

\begin{pro}
In the setting above, we have
\begin{eqnarray}
&& \xxi(\ovl i,p,k)
=\frac{1}{\ci(i-1,1)}
\sum_{\tiny\begin{array}{l}i-1\leq i_2\leq\cd\leq i_a\leq n\\
M=M_1\sqcup\cd\sqcup M_S\in{\mathcal M}^{(p)}_k
\end{array}}
C_{i_2,i_3,\cd,i_a}^M\cdot D^M,
\label{c-xi-1}\\ 
%See note 15 pp20
&&\xxi(\ovl i,p,\ovl k)=\sum_{i=i_1\leq i_2\leq \cd\leq i_p\leq k}
(\ci(i_1-1,1)\ci(i_2-1,2)\cd \ci(i_p-1,p))^{-1}(\ci(i_2,1)\ci(i_3,2)\cd
\ci(i_p,p-1)\ci(k,p)).
%See note 15 pp17
\label{c-xi-2}
\end{eqnarray}
\end{pro}
{\sl Proof.}
Set $\cX:={\pmb x}_n(c_n)\cd{\pmb x}_1(c_1)$. By calculating directly
we have the formula:
\begin{eqnarray}
&&\cX v_i=\begin{cases}c_1^{-1}v_1&\text{ if }i=1,\\
c_{i-1}c_i^{-1}v_i+v_{i-1}&\text{ if }i=2,\cd,n,
\end{cases}\\
&&\cX v_{\ovl i}=
c_{i-1}^{-1}(c_iv_{\ovl i}+c_{i+1}v_{\ovl{i+1}}+\cd +c_{n-1}v_{\ovl{n-1}}+
c_nv_{\ovl n}+v_n),
\end{eqnarray}
where we understand $c_0=1$. Using these, for $k=1,2,\cd,n$ and
$p=2,\cd,n$
we get 
\begin{eqnarray}
&&\xxi(\ovl i,p,\ovl k)=\sum_{j=i}^k\xxi(\ovl i,p-1,\ovl j)
\frac{\ci(k,p)}{\ci(j-1,p)},\label{recur-c-f}\\
&&\xxi(\ovl i,p,k)=\xxi(\ovl i,p-1,k+1)+\xxi(\ovl i,p-1,k)
\frac{\ci(k-1,p)}{\ci(k,p)}.
%See note 15 pp21
\label{recur-c-l}
\end{eqnarray}
Indeed, the formula (\ref{c-xi-2}) is easily shown by the induction
on $p$ using (\ref{recur-c-f}).

To obtain (\ref{c-xi-1}) we see the segments of elements in $\cM_{k}^{(p)}$, 
$\cM_{k+1}^{(p-1)}$ and $\cM_{k}^{(p-1)}$.
$M=M_1\sqcup\cd\sqcup M_{n-k}\in{\mathcal M}^{(p-1)}_{k+1}$ 
can be seen as 
an element in $\cM_{k}^{(p)}$ by setting $M_{n-k+1}=\emptyset$.
For any $M'=M'_1\sqcup\cd\sqcup M'_{n-k+1}\in{\mathcal
M}^{(p-1)}_{k}$, the last segment $M'_{n-k+1}$ is empty or includes
$p-1$. 
Indeed, the following lemma insures this fact:
\begin{lem}
\label{lem-M}
%See Note 17 pp11-12.
For any $M=M_1\sqcup\cd\sqcup M_S\in\cM^{(p)}_k$ and any $i=1,\cd,S:=n-k+1$
we have 
\begin{equation}
\min M_i\geq i+1,\qq \max M_i\leq L+i-1.
\end{equation}
\end{lem}
The proof of this lemma is done by the induction on $p$.\\
Thus, in any case if we set $M''_{n-k+1}=M'_{n-k+1}\cup\{p\}$, then
$M'':=M'_1\sqcup\cd\sqcup M'_{n-k}\sqcup M''_{n-k+1}$ turns out to be
 an element in $\cM_{k}^{(p)}$ and we have
\[
 C^{M''}_{i_2,\cd,i_a}\cdot D^{M''}=
 C^{M'}_{i_2,\cd,i_a}\cdot D^{M'}\frac{\ci(k-1,p)}{\ci(k,p)}.
\]
Using this formula and \eqref{recur-c-l}, we obtain \eqref{c-xi-1}.\qed

Thus, for example, we have 
$\Del_{w_0\Lm_1,s_1\Lm_1}(\Theta^-_\bfii0(c))=\xxi(\ovl 1,n,2)
=\sum_{j=1}^n\frac{\ci(j,1)}{\ci(j-1,2)}+
\sum_{j=2}^n\frac{\ci(n-j+1,j)}{\ci(n-j+2,j)}$.

To get the explicit form of 
$\Del_{w_0\Lm_k,s_k\Lm_k}(\Theta^-_\bfii0(c))$, we show the following
lemma:
\begin{lem}\label{lem-det-c}
For $k=1,\cd, n-1$ set 
\begin{equation}
W_k:=
\begin{pmatrix}
\xxi(\ovl 1,n,k+1)&\xxi(\ovl 1,n,k-1)&\cd&\xxi(\ovl 1,n,2)&
\xxi(\ovl 1,n,1)\\
\xxi(\ovl 2,n,k+1)&\xxi(\ovl 2,n,k-1)&\cd&\xxi(\ovl 2,n,2)&
\xxi(\ovl 2,n,1)\\
\vdots&\vdots&\cd&\vdots&\vdots\\
\xxi(\ovl k,n,k+1)&\xxi(\ovl k,n,k-1)&\cd&\xxi(\ovl k,n,2)&
\xxi(\ovl k,n,1)
\end{pmatrix}.
\label{det-c}
\end{equation}
Then, we have $\Del_{w_0\Lm_k,s_k\Lm_k}(\Theta^-_\bfii0(c))=\det W_k$.
\end{lem}
{\sl Proof.}
As has been given in (\ref{h-l}), for the lowest weight vector
$v_{\Lm_k}\in V(\Lm_k)$ 
and $X:=X^{(n)}\cd X^{(1)}\in G$ we have 
\[
X v_{\Lm_k}=\displaystyle
\sum_{\sigma\in{\mathfrak S}_k}{\rm sgn}(\sigma)
X v_{\ovl{\sigma(k)}}\ot\cd\ot X v_{\ovl{\sigma(1)}}.
\]
Here note that for the simple reflection $s_k\in W$, 
\[
\ovl s_k(v_1\ot\cd\ot v_{k-1}\ot v_{k})=
v_1\ot\cd\ot v_{k-1}\ot v_{k+1}.
\]
Thus, it follows from the formula (\ref{del-x}) that 
$\Del_{w_0\Lm_k,s_k\Lm_k}(\Theta^-_\bfii0(c))$ coincides with 
the coefficient of the vector $v_1\ot\cd\ot v_{k-1}\ot v_{k+1}$
in $X v_{\Lm_k}$, which completes the proof.
\qed

The last column of $W_k$ is just 
\[
{}^t(\xxi(\ovl 1,n,1),\xxi(\ovl 2,n,1),\cd,\xxi(\ovl k,n,1)) 
={}^t(1,{\ci(1,1)}^{-1},{\ci(2,1)}^{-1},\cd,{\ci(k-1,1)}^{-1}).
\]
Considering the elementary transformations on $W_k$ by 
($i$-th row) $-\frac{\ci(i,1)}{\ci(i-1,1)}\times$($i+1$-th row)
for $i=1,\cd,k-1$, the $(i,j)$-entry of 
the transformed matrix $\til W_k$
 is as follows:
\begin{lem}\label{lem-w-c}
 The $(i,j)$-entry $(\til W_k)_{i,j}$ is:
\begin{equation}
(\til W_k)_{i,j}=\begin{cases}
\xxi(\ovl i,n-1,k+1)(c^{[2:n]})&\text{ if }j=1,\\
\xxi(\ovl i,n-1,k-j+1)(c^{[2:n]})&\text{ if }j>1,
\end{cases}
\end{equation}
where for $(\ci(k,l))_{1\leq k,l\leq n}$ we set 
$c^{(l)}:=(\ci(1,l),\cd,\ci(n,l))$ and 
$c^{[a:b]}:=(c^{(a)},c^{(a+1)},\cd,c^{(b)})$ for $1\leq a<b\leq n$. 
Then, $c=c^{[1:n]}$ as above.
Note that $\xxi(\ovl i,p,k)(c)$ depends only on $c^{[1:p]}$.
\end{lem}
{\sl Proof.}
We shall show 
\begin{equation}
\xxi(\ovl i,n,j)(c)-\frac{\ci(i,1)}{\ci(i-1,1)}\xxi(\ovl{i+1},n,j)(c)
=\frac{\ci(i,1)}{\ci(i-1,1)} \xxi(\ovl i,n-1,j)(c^{[2,n]}).
\label{xi-xi}
\end{equation}
For $M=M_1\sqcup\cd\sqcup M_{n-j+1}\in\cM_j^{(n)}$ such that
$M_1=\{2,3,\cd,a\}$ is non-empty, let $M':=M\setminus\{2\}$ be an element
in $\cM_j^{(n-1)}[1]$. Then 
we have 
\[
 C^M_{i-1,i_3,\cd,i_a}\cdot D^M=\frac{\ci(i,1)}{\ci(i-1,2)}
C^{M\setminus\{2\}}_{i_3,\cd,i_a}\cdot D^{M\setminus\{2\}},
\]
where $M\setminus\{2\}$ is considered as an element in
$\cM^{(n-1)}_j[1]$ and by \eqref{shift-cd}
the left-hand side of (\ref{xi-xi}) is written as
\[
\frac{\ci(i,1)}{\ci(i-1,1)\ci(i-1,2)} 
\sum_{\tiny\begin{array}{l}i-1\leq i_3\leq\cd\leq i_a\leq n\\
M'(=M\setminus\{2\})=M'_1
\sqcup\cd\sqcup M'_{n-j+1}\in{\mathcal M}^{(n-1)}_j[1]
\end{array}}
C^{M'}_{i_3,\cd,i_a}\cdot D^{M'},
\]
which shows (\ref{xi-xi}).\qed
%Note 15 pp25\& Note 17 pp24-25

Applying the above elementary transformations to the matrix $W_k$
repeatedly, we have the following:
\begin{cor}
For $k,j$ such that $1\leq k<j\leq n$, we get
\begin{equation}
\xxi(\ovl i,n,j)(c^{[k:n]})-\frac{\ci(i,1)}{\ci(i-1,1)}
\xxi(\ovl{i+1},n,j)(c^{[k:n]})
=\frac{\ci(i,1)}{\ci(i-1,1)} \xxi(\ovl i,n-1,j)(c^{[k+1,n]}).
\label{xi-xi-k}
\end{equation}
\end{cor}
Then we 
find the following equalities of determinants:
\begin{eqnarray*}
\det W_k&=&
\begin{vmatrix}
\ci(1,1)\xxi(\ovl 1,n-1,k+1)(c^{[2:n]})
&\ci(1,1)\xxi(\ovl 1,n-1,k-1)(c^{[2:n]})
&\cd
&\ci(1,1)\xxi(\ovl 1,n-1,2)(c^{[2:n]})&0\\
\frac{\ci(2,1)}{\ci(1,1)}\xxi(\ovl 2,n-1,k+1)(c^{[2:n]})
&\frac{\ci(2,1)}{\ci(1,1)}\xxi(\ovl 2,n-1,k-1)(c^{[2:n]})&\cd
&\frac{\ci(2,1)}{\ci(1,1)}\xxi(\ovl 2,n-1,2)(c^{[2:n]})&0\\
\vdots&\vdots&\cd&\vdots&\vdots\\
\frac{\ci(k-1,1)}{\ci(k-2,1)}\xxi(\ovl{k-1},n-1,k+1)(c^{[2:n]})
&\frac{\ci(k-1,1)}{\ci(k-2,1)}\xxi(\ovl{k-1},n-1,k-1)(c^{[2:n]})&\cd
&\frac{\ci(k-1,1)}{\ci(k-2,1)}\xxi(\ovl{k-1},n-1,2)(c^{[2:n]})&0\\
\xxi(\ovl k,n,k+1)(c^{[2:n]})&\xxi(\ovl k,n,k-1)(c^{[2:n]})&\cd
&\xxi(\ovl k,n,2)(c^{[2:n]})&
{\ci(k-1,1)}^{-1}
\end{vmatrix}\\
&=&
\begin{vmatrix}
\xxi(\ovl 1,n-1,k+1)(c^{[2:n]})
&\xxi(\ovl 1,n-1,k-1)(c^{[2:n]})
&\cd
&\xxi(\ovl 1,n-1,2)(c^{[2:n]})\\
\xxi(\ovl 2,n-1,k+1)(c^{[2:n]})
&\xxi(\ovl 2,n-1,k-1)(c^{[2:n]})&\cd
&\xxi(\ovl 2,n-1,2)(c^{[2:n]})\\
\vdots&\vdots&\cd&\vdots\\
\xxi(\ovl{k-1},n-1,k+1)(c^{[2:n]})
&\xxi(\ovl{k-1},n-1,k-1)(c^{[2:n]})&\cd
&\xxi(\ovl{k-1},n-1,2)(c^{[2:n]})
\end{vmatrix}\\
&=&\cd\cd\cd\\
&=&
\begin{vmatrix}
\xxi(\ovl 1,n-k+2,k+1)(c^{[k-1:n]})
&\xxi(\ovl 1,n-k+2,k-1)(c^{[k-1:n]})\\
\xxi(\ovl 2,n-k+2,k+1)(c^{[k-1:n]})
&\xxi(\ovl 2,n-k+2,k-1)(c^{[k-1:n]})
\end{vmatrix}
=\xxi(\ovl 1,n-k+1,k+1)(c^{[k:n]})
\end{eqnarray*}
Thus, it follows from  (\ref{c-xi-1})  that for $k=1,2,\cd,n-1$
\begin{equation}
\Del_{w_0\Lm_k,s_k\Lm_k}(\Theta^-_\bfi(c))
=\xxi(\ovl 1,n-k+1,k+1)(c^{[k:n]})=
\ci(1,k)+\frac{\ci(2,k)}{\ci(1,k+1)}+\cd+
\frac{\ci(n,k)}{\ci(n-1,k+1)}+
\frac{\ci(n-1,k+1)}{\ci(n,k+1)}+\frac{\ci(n-2,k+2)}{\ci(n-1,k+2)}+\cd+
\frac{\ci(k,n)}{\ci(k+1,n)}
\end{equation}
The case $k=n$ is easily obtained by the formula in \cite[Sect.4]{BZ2}:
\begin{equation}
\Del_{w_0\Lm_n,s_n\Lm_n}(\Theta^-_\bfii0(c))=\ci(1,n).
\end{equation}
Now, the proof of Theorem \ref{thm-c-1} has been accomplished.
\qed

%%%%%%%%%%%%%%%%%%%%%%%%%%%%%%
\subsection {Proof of Theorem \ref{thm-c-2}}%7.3 

For $k=1,2,\cd,n-1$,  set 
\[
U_k:=\begin{pmatrix}
\xxi(\ovl 1,n,k)&\xxi(\ovl 1,n,k-1)&\cd&\xxi(\ovl 1,n,2)&
\xxi(\ovl 1,n,1)\\
\xxi(\ovl 2,n,k)&\xxi(\ovl 2,n,k-1)&\cd&\xxi(\ovl 2,n,2)&
\xxi(\ovl 2,n,1)\\
\vdots&\vdots&\cd&\vdots&\vdots\\
\xxi(\ovl{k-1},n,k)&\xxi(\ovl{k-1},n,k-1)&\cd&\xxi(\ovl{k-1},n,2)&
\xxi(\ovl{k-1},n,1)\\
\xxi(\ovl{k+1},n,k)&\xxi(\ovl{k+1},n,k-1)&\cd&\xxi(\ovl{k+1},n,2)&
\xxi(\ovl{k+1},n,1)
\end{pmatrix}.
\]
Since 
$\Del_{w_0s_k\Lm_k,\Lm_k}(\Theta_\bfii0(c))=\lan u_{\Lm_k}\,,\,
X^{(n)}\cd X^{(1)}\ovl s_kv_{\Lm_k}\ran$, 
the function $\Del_{w_0s_k\Lm_k,\Lm_k}(\Theta_\bfii0(c))$ is given as
the coefficient of the vector $v_1\ot v_2\ot \cd\ot v_k$ in 
$X^{(n)}\cd X^{(1)}\ovl s_kv_{\Lm_k}$.
Thus, by the argument in the proof of Theorem \ref{thm-c-1},
we obtain:
\begin{equation}
\Del_{w_0s_k\Lm_k,\Lm_k}(\Theta^-_\bfii0(c))=
\det U_k\q(k=1,2,\cd,n-1).
\label{vmat-c2}
\end{equation}
Using the formula (\ref{xi-xi}), we have for $j=1,2,\cd,k$
\begin{equation}
\xxi(\ovl{k-1},n,j)(c)-\frac{\ci(k,1)}{\ci(k-2,1)}\xxi(\ovl{k+1},n,j)(c)
=\frac{\ci(k-1,1)}{\ci(k-2,i)}\xxi(\ovl{k-1},n-1,j)(c^{[2:n]})
+\frac{\ci(k,1)}{\ci(k-2,1)}\xxi(\ovl k,n-1,j)(c^{[2:n]}).
\label{xi-xi2}
\end{equation}
Here, applying this formula to the above determinant, we have 
%See Note 15 pp70
\begin{eqnarray*}
&&\Del_{w_0s_k\Lm_k,\Lm_k}(\Theta^-_\bfii0(c))\\
&&=\frac{\ci(k-1,1)}{\ci(k,1)}
\begin{vmatrix}
\xxi(\ovl 1,n-1,k)(c^{[2:n]})&\xxi(\ovl 1,n-1,k-1)(c^{[2:n]})
&\cd&\xxi(\ovl 1,n-1,2)(c^{[2:n]})\\
\xxi(\ovl 2,n-1,k)(c^{[2:n]})&\xxi(\ovl 2,n-1,k-1)(c^{[2:n]})
&\cd&\xxi(\ovl 2,n-1,2)(c^{[2:n]})\\
\vdots&\vdots&\cd&\vdots\\
\xxi(\ovl{k-1},n-1,k)(c^{[2:n]})&\xxi(\ovl{k-1},n-1,k-1)(c^{[2:n]})
&\cd&\xxi(\ovl{k-1},n-1,2)(c^{[2:n]})
\end{vmatrix}\\
&&\qq \,\,+
\begin{vmatrix}
\xxi(\ovl 1,n-1,k)(c^{[2:n]})&\xxi(\ovl 1,n-1,k-1)(c^{[2:n]})
&\cd&\xxi(\ovl 1,n-1,2)(c^{[2:n]})\\
\xxi(\ovl 2,n-1,k)(c^{[2:n]})&\xxi(\ovl 2,n-1,k-1)(c^{[2:n]})
&\cd&\xxi(\ovl 2,n-1,2)(c^{[2:n]})\\
\vdots&\vdots&\cd&\vdots\\
\xxi(\ovl{k-1},n-1,k)(c^{[2:n]})&\xxi(\ovl{k-1},n-1,k-1)(c^{[2:n]})
&\cd&\xxi(\ovl{k-1},n-1,2)(c^{[2:n]})\\
\xxi(\ovl{k+1},n-1,k)(c^{[2:n]})&\xxi(\ovl{k+1},n-1,k-1)(c^{[2:n]})
&\cd&\xxi(\ovl{k+1},n,2)(c^{[2:n]})
\end{vmatrix}.
\end{eqnarray*}
In the above formula, let $Z_{k-1}$  be the matrix 
in the first determinant. Thus, we have 
\begin{equation}
\det U_k=\frac{\ci(k-1,1)}{\ci(k,1)}\det Z_{k-1}+\det U_{k-1}.
\label{recur-c-2}
\end{equation}
Repeating these steps, we can derive 
\[
\det U_k=\det U_1+\sum_{j=1}^{k-1}\frac{\ci(k-j,j)}{\ci(k-j+1,j)}\det Z_j.
\]
Carrying out the elementary transformations above to the matrix 
$Z_{k-1}$ we easily know that $\det Z_{j}=1$ for $j=1,\cd,k-1$ and 
$\det U_1=\xxi(\ovl2,n-k+1,k)(c^{[k:n]})={\ci(1,k)}^{-1}$, which show
\begin{equation}
\Del_{w_0s_k\Lm_k,\Lm_k}(\Theta^-_\bfii0(c))
=\frac{1}{\ci(1,k)}+\sum_{j=1}^{k-1}\frac{\ci(k-j,j)}{\ci(k-j+1,j)},
\end{equation}
and then (\ref{del-c-2}).

Next, to show (\ref{del-c-3}), we shall see
$\Del_{w_0\Lm_n,s_n\Lm_n}(\Theta^-_{\bfii0^{-1}}(c))$ 
since for $g\in B^-_{w_0}$ we have 
\begin{equation}
\Del_{w_0s_i\Lm_i,\Lm_i}(g)=\Del_{w_0\Lm_i,s_i\Lm_i}(\eta(g)),
\label{eta-si-c}
\end{equation}
where $\bfii0^{-1}=(n\, n-1 \cd 21)^n$ and 
$\eta(\Theta^-_{\bfii0^{-1}}(c))=\Theta^-_{\bfii0}(\ovl c)$
($-:\ci(i,j)\mapsto{\ci(n-i+1,j)}$).
%See Note 16 p90-91.
These mean
\[
 \Del_{w_0s_i\Lm_i,\Lm_i}(\Theta^-_{\bfii0}(c))
=\Del_{w_0\Lm_i,s_i\Lm_i}(\Theta^-_{\bfii0^{-1}}(c')),
\]
where $c'=(c_N,c_{N-1},\cd,c_1)$ for $c=(c_1,\cd,c_N)$.
As we have seen $\omega({\pmb y}_i(c))=\al_i^\vee(c^{-1})x_i(c)=
{\pmb x}_i(c)$, we have 
$\omega(\Theta^-_{\bfii0^{-1}}(c'))=\ovl X^{(n)}\ovl X^{(n-1)}\cd\ovl
X^{(1)}$
where $\ovl X^{(p)}={\pmb x}_1(\ci(1,p)){\pmb x}_{2}(\ci(2,p))\cd 
{\pmb x}_n(\ci(n,p))$.
We define $\xsi(\ovl i,p,j)$ and $\xsi(i,p,j)$ by 
\[
\ovl X^{(p)}\ovl X^{({p-1})}\cd \ovl X^{(1)}v_{\ovl i}
=\sum_{j=1}^n{\xsi(\ovl i,p,j)}v_j
+\sum_{j=1}^n{\xsi(\ovl i,p,\ovl j)}v_{\ovl j}\in V(\Lm_1).
\]
To describe these coefficients explicitly, we 
 define the similar objects to the segments as above.
For $1\leq p,i\leq n$, set $L:=p-n+i-1$ and $S:=n-i+2$, 
\[
\ovl\cM_i^{(p)}:=\{M=\{m_1,\cd,m_L\}|1\leq m_1\leq \cd\leq m_L\leq p\}.
\]
For an element $M$, define the segments $M_1,\cd,M_S$ by the same way as
before, thus, $M=M_1\sqcup\cd\sqcup M_S$. 
For an element $M=M_1\sqcup\cd\sqcup M_S\in \ovl\cM_i^{(p)}$ writing
$M_S=\{q,q+1,\cd,p-1,p\}$, set
\begin{eqnarray*}
&&D^M:=\prod_{j=1}^S\prod_{m\in
 M_j}\frac{\ci(i+j-1,m)}{\ci(i+j-2,m)},
\qq
%See Note 15, pp84
F^{M_S}_{j_q,j_{q+1},\cd,j_p}:=\frac{\ci(j_q-1,q-1)\ci(j_{q+1}-1,q)\cd
\ci(j_{p-1}-1,p-2)\ci(j_p-1,p-1)}{\ci(j_q,q)\ci(j_{q+1},q+1)\cd
\ci(j_{p-1},p-1)\ci(j_p,p)}.
\end{eqnarray*}
\begin{lem}We have
\begin{equation}
\xsi(\ovl i,p,k)=\ci(k-1,p)
\sum_{\tiny\begin{array}{l}k\leq j_p\leq j_{p-1}\cd\leq j_q\leq n\\
M=M_1\sqcup\cd\sqcup M_S\in\ovl\cM^{(p)}_i
\end{array}}
D^M\cdot F^{M_S}_{j_q,j_{q+1},\cd,j_p},
\end{equation}
where $M_S=\{q,q+1,\cd,p-1,p\}$ and $S=n-i+2$.
Note that $\xsi(\ovl 1,n,k)=\ci(k-1,n)$.
\end{lem}
By the similar argument as the above cases, we have 
\begin{equation}
\Del_{w_0\Lm_n,s_n\Lm_n}(\Theta^-_{\bfii0^{-1}}(c))
=\begin{vmatrix}
\xsi(\ovl 1,n,\ovl n)&\xsi(\ovl 1,n,n-1)&\cd&\xsi(\ovl 1,n,2)&
\xsi(\ovl 1,n,1)\\
\xsi(\ovl 2,n,\ovl n)&\xsi(\ovl 2,n,n-1)&\cd&\xsi(\ovl 2,n,2)&
\xsi(\ovl 2,n,1)\\
\vdots&\vdots&\cd&\vdots&\vdots\\
\xsi(\ovl n,n,\ovl n)&\xsi(\ovl n,n,n-1)&\cd&\xsi(\ovl n,n,2)&
\xsi(\ovl n,n,1)
\end{vmatrix}
\label{det-c-3}
\end{equation}
%See note 15, pp87
To calculate this determinant we need some preparations.
Let us define the functions:
\begin{eqnarray}
&&\xE(\ovl i,k,j)(c):=
\sum_{\tiny\begin{array}{l}k\leq j_k\leq j_{k-1}\cd\leq j_q\leq n\\
M=M_1\sqcup\cd\sqcup M_S\in\ovl\cM^{(p)}_i,\\
M_S=\{q,q+1,\cd,k\}\ne\emptyset,\\
j_k<u_{k+1}<u_{k+2}\cd<u_j\leq n
\end{array}}
D^M\cdot F^{M_S}_{j_q,j_{q+1},\cd,j_k}
\cdot F^{\{k+1,k+2,\cd,j\}}_{u_{k+1},u_{k+2},\cd,u_j}
\q(j\geq k),
\label{def-e}\\
%\frac{\ci(u_{k+1}-1,k)\ci(u_{k+2}-1,k+1)\cd\ci(u_j-1,j-1)}{
%\ci(u_{k+1},k+1)\ci(u_{k+2},k+2)\cd\ci(u_j,j)},
&&\xF(\ovl i,k)(c):=
\sum_{\tiny\begin{array}{l}
M=\{m_1,\cd,m_{k-n+i-1}\}\\
=M_1\sqcup\cd\sqcup M_S\in\ovl\cM^{(p)}_i,\q
M_S=\emptyset,\\
1\leq l_1\leq l_2\leq\cd\leq l_{n-k+1}\leq m_{k-n+i-1}\end{array}}
D^M\cdot\frac{\ci(l_1,l_1)\ci(l_2+1,l_2)\ci(l_3+2,l_3)\cd
\ci(l_{n-k+1}+n-k,l_{n-k+1})}{\ci(l_1-1,l_1)\ci(l_2,l_2)\ci(l_3+1,l_3)\cd
\ci(l_{n-k+1}+n-k-1,l_{n-k+1})}.
\label{def-f}
\end{eqnarray}
%See note 15, pp92-94-->teisei note 16, pp11-13
\begin{lem}
We have the following formula for $2\leq k\leq n$ and $2\leq p\leq n$:
\begin{eqnarray}
&&\xsi(\ovl i,p,k)-\frac{\ci(k-1,p)}{\ci(k,p)}\xsi(\ovl i,p,k+1)=
\frac{\ci(k-1,p)}{\ci(k,p)}\xsi(\ovl i,p-1,k),\q \xsi(\ovl 1,n,k)=\ci(k-1,n),
\label{4-1}\\
&&\xsi(\ovl i,n,n-1)-\frac{\ci(n-2,n)}{\xsi(\ovl 1,n,\ovl n)}\cdot
\xsi(\ovl i,n,\ovl n)=
\frac{\ci(n-2,n)}{\ci(n-1,n)}\cdot\xsi(\ovl i,n-1,n-1)+\ci(n-2,n)
\cdot\xE(\ovl i,n,n)+\frac{\ci(n-2,n)}{\xsi(\ovl 1,n,n)}\cdot
\xF(\ovl i,n),
\label{4-2}\\
%See note 15 pp94
&&\xF(\ovl{n-k+2},k)\xsi(\ovl i,k-1,k-1)-\ci(k-2,k-1)
\xF(\ovl i,k)=\ci(k-2,k-1)\cdot\xF(\ovl{n-k+1},k)\cdot\xE(\ovl i,k-1,k-1)
+\ci(k-2,k-1)\xF(\ovl i,k-1),
\label{4-3}\\
%See note 16 pp6
&&\xE(\ovl{n-k+2},k,j)\cdot\xsi(\ovl i,k-1,k-1)-\ci(k-2,k-1)\cdot
\xE(\ovl i,k,j)=\ci(k-2,k-1)\cdot\xE(\ovl i,k-1,j).
\label{4-4}
%See note 16 pp14
\end{eqnarray}
\end{lem}
Direct inspections show this lemma.\qed

Let us denote the right-hand side of \eqref{4-2} by $X_i$ $(i=2,\cd,n)$.
For the right-hand side of \eqref{det-c-3} subtracting the 
$\frac{\ci(k-1,n)}{\ci(k,n)}\times (n-k)$-th 
column from $n-k+1$th column for $k=1,2,\cd,n-1$ 
and applying the formula \eqref{4-1} and \eqref{4-2}, we obtain
\begin{eqnarray}
&&\hspace{0.5cm}\text{R.H.S of (\ref{det-c-3})}\label{sigma-3}\\
&&=\frac{1}{\ci(n-2,n)}\begin{vmatrix}
\xsi(\ovl 1,n,\ovl n)&0&\cd&0&0\\
\xsi(\ovl 2,n,\ovl n)&X_2&\cd&\xsi(\ovl 2,n-1,2)&
\xsi(\ovl 2,n-1,1)\\
\vdots&\vdots&\cd&\vdots&\vdots\\
\xsi(\ovl n,n,\ovl n)&X_n&\cd&\xsi(\ovl n,n-1,2)&
\xsi(\ovl n,n-1,1)
\end{vmatrix}
=\frac{\xsi(\ovl 1,n,\ovl n)}{\ci(n-1,n)}
\begin{vmatrix}
\xsi(\ovl 2,n-1,n-1)&\cd&\xsi(\ovl 2,n-1,2)&\xsi(\ovl 2,n-1,1)\\
\vdots&\cd&\vdots&\vdots\\
\xsi(\ovl n,n-1,{n-1})&\cd&\xsi(\ovl n,n-1,2)&
\xsi(\ovl n,n-1,1)
\end{vmatrix}\nn\\
&&+
\xsi(\ovl 1,n,\ovl n)
\begin{vmatrix}
\xE(\ovl 2,n,n)&\xsi(\ovl 2,n-1,n-2)
&\cd&\xsi(\ovl 2,n-1,2)&\xsi(\ovl 2,n-1,1)\\
\vdots&\vdots&\cd&\vdots&\vdots\\
\xE(\ovl n,n,n)&\xsi(\ovl n,n-1,n-2)&\cd&\xsi(\ovl n,n-1,2)&\xsi(\ovl n,n-1,1)
\end{vmatrix}
+
\begin{vmatrix}
\xF(\ovl 2,n)&\xsi(\ovl 2,n-1,n-2)
&\cd&\xsi(\ovl 2,n-1,2)&\xsi(\ovl 2,n-1,1)\\
\vdots&\vdots&\cd&\vdots&\vdots\\
\xF(\ovl n,n)&\xsi(\ovl n,n-1,n-2)&\cd&\xsi(\ovl n,n-1,2)&\xsi(\ovl
 n,n-1,1)
\nn
\end{vmatrix}
\end{eqnarray}
To complete this calculations, we see the following lemma:
\begin{lem}\label{lem:eta}
Set 
\begin{eqnarray}
\eta_k[j]&&:=\begin{vmatrix}
%See note 16 pp14
\xE(\ovl{n-k+2},k,j)&\xsi(\ovl{n-k+2},k-1,k-2)
&\cd&\xsi(\ovl{n-k+2},k-1,2)&\xsi(\ovl{n-k+2},k-1,1)\\
\vdots&\vdots&\cd&\vdots&\vdots\\
\xE(\ovl{n},k,j)&\xsi(\ovl n,k-1,k-2)
&\cd&\xsi(\ovl n,k-1,2)&\xsi(\ovl n,k-1,1)
\end{vmatrix},
\label{eta-k-def}\\
\phi_k&&:=
\begin{vmatrix}
\xF(\ovl{n-k+2},k)&\xsi(\ovl{n-k+2},k-1,k-2)
%See note 16 pp8
&\cd&\xsi(\ovl{n-k+2},k-1,2)&\xsi(\ovl{n-k+2},k-1,1)\\
\vdots&\vdots&\cd&\vdots&\vdots\\
\xF(\ovl{n},k)&\xsi(\ovl n,k-1,k-2)
&\cd&\xsi(\ovl n,k-1,2)&\xsi(\ovl n,k-1,1)
\end{vmatrix}.
\label{phi-k-def}
\end{eqnarray}
Then we have 
\begin{eqnarray}
&&\eta_k[j]=\eta_{k-1}[j]+\frac{\xE(\ovl{n-k+2},k,j)(c)}{\ci(k-2,k-1)},
\label{eta-k}\\
%see note 16 pp14-15
&&\phi_k=\phi_{k-1}+(\frac{1}{\ci(k-2,k-1)}
+\eta_{k-1}[k-1])\xF(\ovl{n-k+2},k)(c).
\label{phi-k}
%see note 16 pp8
\end{eqnarray}
\end{lem}
{\sl Proof of Lemma\ref{lem:eta}} 
By the formula (\ref{4-1}) and (\ref{4-4}), we obtain 
\begin{equation}
\xE(\ovl{n-k+2},k,j)(\xsi(\ovl i,k-1,k-2)-\ci(k-3,k-1)\cdot
\xE(\ovl i,k,j))=\frac{\ci(k-3,k-1)}{\ci(k-2,k-1)}
\xE(\ovl{n-k+2},k,j)\cdot\xsi(\ovl i,k-2,k-2)+
\ci(k-3,k-1)\cdot \xE(\ovl i,k-1,j).
\label{4-5}
%Note 16,pp14
\end{equation}
Considering the column transformations as above to the determinant
in (\ref{eta-k-def}) and 
applying the formula (\ref{4-1}) and (\ref{4-5}) repeatedly, 
we obtain (\ref{eta-k}).
Similarly, applying (\ref{4-1}) and (\ref{4-3}) to the determinant
in (\ref{phi-k-def}) repeatedly, we obtain (\ref{phi-k}).\qed

From \eqref{sigma-3}, we have 
\begin{equation}
\Del_{w_0\Lm_n,s_n\Lm_n}(\Theta^-_{\bfii0^{-1}}(c))
=\frac{\xsi(\ovl 1,n,\ovl n)}{\ci(n-1,n)}+
\xsi(\ovl 1,n,\ovl n)\cdot\eta_n[n]+\phi_n.
\label{3sum}
\end{equation}
By (\ref{eta-k}) and (\ref{phi-k}), we get that 
\begin{eqnarray}
&&\eta_k[j]=\sum_{l=2}^k{\ci(l-2,l-1)}^{-1}\xE(\ovl{n-l+2},l,j)(c),\\
&&\phi_k=\phi_2+\sum_{j=3}^k\xF(\ovl{n-j+2},j)(c)
({\ci(j-2,j-1)}^{-1}+\eta_{j-1}[j-1]).
%See note 16 pp 15
\end{eqnarray}
Substituting these to (\ref{3sum}), we have 
\begin{equation}
\Del_{w_0\Lm_n,s_n\Lm_n}(\Theta^-_{\bfii0^{-1}}(c))
=\phi_2+\sum_{j=3}^{n+1}\xF(\ovl{n-j+2},j)(c)
\cdot({\ci(j-2,j-1)}^{-1}+
\sum_{l=2}^{j-1}{\ci(l-2,l-1)}^{-1}\xE(\ovl{n-l+2},l,j-1)(c)),
\label{EF}
\end{equation}
where we understand $\xF(\ovl 1,n+1)(c)=1$.
Using (\ref{eta-k}) and the explicit form of $\xE(\ovl i,k,j)$
in (\ref{def-e}), we get the following theorem:
\begin{thm}
We have
\begin{equation}
\Del_{w_0\Lm_n,s_n\Lm_n}(\Theta^-_{\bfii0^{-1}}(c))
=\sum_{(*)}
\di(l_1,l_1)\cdot\di(l_2,l_2-1)\cd \di(l_k,{l_k-k+1})\cdot
\dbi(u_1,2)\dbi(u_2,3)\cd \dbi(u_m,m+1),
\label{thm-c-inv}
\end{equation}
where $\di(i,j)=\frac{\ci(i,j)}{\ci(i-1,j)}$, 
$\dbi(i,j)=\frac{\ci(i-1,j-1)}{\ci(i,j)}$ and 
$(*)$ is the conditions:
$k+m=n$, $0\leq m<n$,\\
$1\leq l_1< l_2<\cd<l_k\leq n,$ and 
$1\leq u_1< u_2<\cd<u_m\leq n.$
\end{thm}
{\sl Proof.}
For $M\in \cM_i^{(p)}$
if $i=n-k+2$ and $M_S\ne\emptyset$, then $L=1$ and then 
$M=M_S=\{k\}$. Thus, it follows from \eqref{def-e} and \eqref{def-f}  
that
\begin{eqnarray}
&&\xE(\ovl{n-k+2},k,j)(c)=\sum_{\begin{array}{l}
\syl k\leq u_l<u_{k+1}\\ \syl <\cd<u_j\leq n\end{array}}
\frac{\ci(u_k-1,k-1)\ci(u_{k+1}-1,k)\cd \ci(u_j-1,j-1)}
{\ci(u_k,k)\ci(u_{k+1},k+1)\cd \ci(u_j,j)}
=\sum_{\begin{array}{l}\syl
k\leq u_l<u_{k+1}\\ \syl<\cd<u_j\leq n\end{array}}
\ovl D^{(k)}_{u_k}\ovl D^{(k+1)}_{u_{k+1}}\cd \ovl D^{(j)}_{u_j},
\label{n-k+2-e}\\
%Note16,pp15
&&\xF(\ovl{n-k+2},k)(c)=\sum_{1\leq l_1<l_2<\cd<l_{n-k+2}\leq n}
D^{(l_1)}_{l_1}D^{(l_2-1)}_{l_2}\cd D^{(l_{n-k+2}-(n-k+1))}_{l_{n-k+2}}.
\label{n-k+2-f}\\
&&\frac{1}{\ci(l-2,l-1)}=\ovl D^{(2)}_1\ovl D^{(3)}_2\cd\ovl
 D^{(l-1)}_{l-2}.
\label{1/c}
%Note 16, pp19
\end{eqnarray}
Then, applying these to \eqref{EF} 
we get \eqref{thm-c-inv}
\qed

Since we have
\[
\Del_{w_0s_n\Lm_n,\Lm_n}(\Theta^-_{\bfii0}(c))
= \Del_{w_0\Lm_n,s_n\Lm_n}(\Theta^-_{\bfii0^{-1}}(\ovl c)),
\]
where for $c=(\ci(i,j))$, we define $\ovl c=({\ci(i,n-j+1)}^{-1})$,
we get the following:
\begin{equation}
\Del_{w_0s_n\Lm_n,\Lm_n}(\Theta^-_\bfii0(c))
=\sum_{(*)}
{\Ci(u_1,1)}{\Ci(u_2,2)}\cd {\Ci(u_m,m)}
{\bCi(l_1,l_1-1)}\cdot{\bCi(l_2,l_2-2)}\cd 
{\bCi(l_{k},l_{k}-k)}
\label{Del-w0sn}
\end{equation}
%note 16 pp 37 red pen's
where 
$\Ci(i,j)=\frac{\ci(i,n-j)}{\ci(i-1,n-j+1)}$,
$\bCi(i,j)=\frac{\ci(i-1,n-j)}{\ci(i,n-j)}$ and 
%$\di(i,j)=\frac{\ci(i,j)}{\ci(i-1,j)}$, 
%$\dbi(i,j)=\frac{\ci(i-1,j-1)}{\ci(i,j)}$ and 
$(*)$ is the conditions:
$k+m=n$, $0\leq m<n$,\\
$1\leq l_1< l_2<\cd<l_k\leq n,$ and 
$1\leq u_1< u_2<\cd<u_m\leq n.$

This shows (\ref{del-c-3}) and then we accomplish the proof
of Theorem \ref{thm-c-2}.
\qed

%%%%%%%%%%%%%%%%%%%%%%%%%%%%%%%%%%%%%%%%%%%%%%%%%%%%%%%
\subsection{Correspondence to the monomial realizations}
Let us present the affirmative answer to
the conjecture for type $C_n$ in this subsection.
First, we see the monomial realization of $B(\Lm_1)$ associated with 
the cyclic oder $\cd(12\cd n)(12\cd n)\cd$, which means that the sign
$p_0=(p_{i,j})$ is given by $p_{i,j}=1$ if $i<j$ and $p_{i,j}=0$ if 
$i>j$. 
The crystal $B(\Lm_1)$ is described as follows:
We abuse the notation $B(\Lm_1):=\{v_i,v_{\ovl i}|1\leq i\leq n\}$ if 
there is no confusion.
Then the actions of $\eit$ and $\fit$ are defined as 
$\fit=f_i$ and $\eit=e_i$ in (\ref{c-f1}) and (\ref{c-f2}).
To describe the monomial realizations, we write down the monomials $A_{i,m}$
associated with $p_0$:
\begin{equation}
A_{i,m}=\begin{cases}
\ci(1,m){\ci(2,m)}^{-1}\ci(1,m+1),&\text{ for }i=1,\\
\ci(i,m){\ci(i+1,m)}^{-1}{\ci(i-1,m+1)}^{-1}\ci(i,m+1),
&\text{ for }1<i\leq n-1,\\
\ci(n,m){\ci(n-1,m+1)}^{-2}\ci(n,m+1)&\text{ for }i=n.
\end{cases}
\label{aim-c-1}
\end{equation}
Here the monomial realization of $B(\Lm_1)$ is described explicitly:
\begin{equation}
B(\ci(1,k))(\cong B(\Lm_1))=\left\{\frac{\ci(j,k)}{\ci(j-1,k+1)}=m_1^{(k)}(v_j),
\frac{\ci(j-1,k+n-j+1)}{\ci(j,k+n-j+1)}=m_{1}^{(k)}(v_{\ovl j})\,\,
|\,\,1\leq j\leq n\,\,\right\},
\end{equation}
where $m^{(k)}_i:B(\Lm_i)\hookrightarrow \cY(p)$ ($u_{\Lm_i}\mapsto 
\ci(i,k)$)
is the embedding of crystal as in Sect.5
and we understand $\ci(0,k)=1$.
Now, Theorem \ref{thm-c-1} claims the following:
\begin{thm}
We obtain  $\Del_{w_0\Lm_n,s_n\Lm_n}(\Theta^-_\bfii0(c))
=m_n^{(n)}(u_{\Lm_n})=\ci(n,n)$
 and 
\begin{eqnarray}
&&\Del_{w_0\Lm_k,s_k\Lm_k}(\Theta^-_\bfii0(c))
=\sum_{j=1}^nm^{(k)}_1(v_j)+\sum_{j=k+1}^nm^{(k)}_1(v_{\ovl j}),\qq
(k=1,2\cd,n-1)
\label{thm-mono-c1}\\
&&\Del_{w_0s_k\Lm_k,\Lm_k}(\Theta^-_\bfii0(c))
=\sum_{j=1}^km^{(k-n)}_1(v_{\ovl j}), \qq(k=1,2,\cd,n-1).
\label{thm-mono-c2}
\end{eqnarray}
\end{thm}

To mention the result for 
$\Del_{w_0s_n\Lm_n,\Lm_n}(\Theta^-_\bfii0(c))$, 
we define the 
following set of monomials:
%{\color{red}????????? highest lowest ?????????????}
\begin{equation}
\bbB:=\left\{
\Ci(u_1,1)\Ci(u_2,2)\cd \Ci(u_m,m)
\bCi(l_1,l_1-1)\cdot\bCi(l_2,l_2-2)\cd 
\bCi(l_{k},l_{k}-k)
\begin{array}{|l}
0\leq m\leq n,\,\,k+m=n, \\%m=n -> highest weight vec \ci(n,0)
1\leq l_1< l_2<\cd<l_k\leq n,\\ 
1\leq u_1< u_2<\cd<u_m\leq n,\\
l_a=u_b\Rightarrow l_a\geq a+b
\end{array}
\right\}
\end{equation}
%note 16 pp 37 red pen's
where $\Ci(i,j)=\frac{\ci(i,n-j)}{\ci(i-1,n-j+1)}$,
$\bCi(i,j)=\frac{\ci(i-1,n-j)}{\ci(i,n-j)}$ 
as above.
Note that if $k=0$ (resp. $m=0$) for an element in $\bbB$, then 
it is just $\ci(n,0)$ (resp. $\frac{1}{\ci(n,n)}$).
We find the following fact:
\begin{pro}
Let $m^{(0)}_n:B(\Lm_n)\to \cY(p)$ be the embedding of crystals 
such that $m^{(0)}_n(u_{\Lm_n})=\ci(n,0)$.
Then, we obtain ${\rm Im}(m^{(0)}_n)=\bbB$ 
and then $B(\ci(n,0))=\bbB$.
\end{pro}
{\sl Proof.}
The map $m_n^{(0)}$ is described explicitly:
\begin{equation}
m^{(0)}_n([u_1,\cd,u_m,\ovl{l_k},\ovl{l_{k-1}},\cd,\ovl{l_1}])
=\Ci(u_1,1)\Ci(u_2,2)\cd \Ci(u_m,m)
\bCi(l_1,l_1-1)\cdot\bCi(l_2,l_2-2)\cd 
\bCi(l_{k},l_{k}-k),
\label{m0}
\end{equation}
where $[u_1,\cd,u_m,\ovl{l_k},\ovl{l_{k-1}},\cd,\ovl{l_1}]$ is an
element in $B(\Lm_n)$ (see \eqref{Ytab-c}), which satisfies 
\[
u_1< u_2\cd< u_m\leq  n< 
\ovl{l_k}<\ovl{l_{k-1}}<\cd<\ovl{l_1}.
\]
Note that $m^{(0)}_n([1,2,\cd,n])=\ci(n,0)$ and 
$m^{(0)}_n([\ovl n,\ovl{n-1},\cd,
\ovl 2,\ovl 1])=\frac{1}{\ci(n,n)}$, which are the highest weight vector 
and the lowest weight vector in $B(\ci(n,0))$ respectively.
So, it suffices to see the compatibility 
$\eit\circ m^{(0)}_n=m^{(0)}_n\circ\eit$ and 
$\fit\circ m^{(0)}_n=m^{(0)}_n\circ\fit$.
To see these, we find the parts including ${\ci(i,j)}^{\pm}$.
Indeed, $\Ci(i,j), \Ci(i+1,j), \bCi(i,j)$ and $\bCi(i+1,j)$ include 
$\ci(i,n-j),{\ci(i,n-j+1)}^{-1},{\ci(i,n-j)}^{-1}$ and $\ci(i,n-j)$
respectively.
There can be $16=2^4$ cases according that 
each of the 4parts exists or not, which correspond to the vectors
in $B(\Lm_n)$ including $i,i+1,\ovl{i+1}$ and $\ovl i$.
For example, if we consider 
$v=[\cd,i,\cd,]$ , that is, if 
$v$ includes only $i$ not $i+1,\ovl{i+1},\ovl i$, then we have 
$m^{(0)}_n(v)$ includes the part $\Ci(i,k)$ where the entry $i$ in $v$ is the
$k$-th entry. In this case, $m^{(0)}_n(v)$ has $\ci(i,n-k)$ and then 
$\vp_i(m^{(0)}_n(v))=1$ and $\vep_i(m^{(0)}_n(v))=0$. 
Then, applying $\fit$ on $m^{(0)}_n(v)$, we have 
$\fit(m^{(0)}_n(v))=m^{(0)}_n(v)\cdot A_{i,n-k}^{-1}$. 
Here we find 
\begin{equation}
\Ci(i,k)\cdot A_{i,n-k}^{-1}=
\begin{cases}\displaystyle\ci(1,n-k)\cdot
\frac{\ci(2,n-k)}{\ci(1,n-k)\ci(1,n-k+1)}
=\frac{\ci(2,n-k)}{\ci(1,n-k+1)}=\Ci(2,k),&\text{ for }i=1,\\
\displaystyle
\frac{\ci(i,n-k)}{\ci(i-1,n-k+1)}\cdot
\frac{\ci(i+1,n-k)\ci(i-1,n-k+1)}{\ci(i,n-k)\ci(i,n-k+1)}
=\frac{\ci(i+1,n-k)}{\ci(i,n-k+1)}=\Ci(i+1,k)
&\text{ for }1<i< n,\\
\displaystyle
\frac{\ci(n,n-k)}{\ci(n-1,n-k+1)}\cdot
\frac{{\ci(n-1,n-k+1)}^{2}}{\ci(n,n-k)\ci(n,n-k+1)}
=\frac{\ci(n-1,n-k+1)}{\ci(n,n-k+1)}=\bCi(n,k-1)&\text{ for }i=n.
\end{cases}
\end{equation}
These correspond to 
\[
\fit(v)=\begin{cases}
[2,\cd]&\text{ for }i=1,\\
[\cd,i+1,\cd]&\text{ for }1<i<n,\\
[\cd,\ovl n,\cd]&\text{ for }i=n.
\end{cases}
\]
Then, in this case we obtain $\fit(m^{(0)}_n(v))=m^{(0)}_n(\fit(v))$.

Next, we consider 
$v=[\cd,\ovl i,\cd,]$ , that is, if 
$v$ includes only $\ovl i$ not $i,i+1,\ovl{i+1}$, then we have 
$m^{(0)}_n(v)$ includes the part $\bCi(i,k)$ where the entry $\ovl i$ 
in $v$ is the
$i-k$-th entry from the bottom. 
In this case, $m^{(0)}_n(v)$ has ${\ci(i,n-k)}^{-1}$ and then 
$\vep_i(m^{(0)}_n(v))=1$. Then, applying $\eit$ on $m^{(0)}_n(v)$, we have 
$\eit(m^{(0)}_n(v))=m^{(0)}_n(v)\cdot A_{i,n-k-1}$. 
Here we find 
\begin{equation}
\bCi(i,k)\cdot A_{i,n-k-1}=
\begin{cases}\displaystyle{\ci(1,n-k)}^{-1}\cdot
\frac{\ci(1,n-k-1)\ci(1,n-k)}{\ci(2,n-k-1)}
=\frac{\ci(1,n-k-1)}{\ci(2,n-k-1)}=\bCi(2,k+1),&\text{ for }i=1,\\
\displaystyle
\frac{\ci(i-1,n-k)}{\ci(i,n-k)}\cdot
\frac{\ci(i,n-k-1)\ci(i,n-k)}{\ci(i+1,n-k-1)\ci(i-1,n-k)}
=\frac{\ci(i,n-k-1)}{\ci(i+1,n-k-1)}=\bCi(i+1,k+1)
&\text{ for }1<i< n,\\
\displaystyle
\frac{\ci(n-1,n-k)}{\ci(n,n-k)}\cdot
\frac{\ci(n,n-k-1)\ci(n,n-k)}{{\ci(n-1,n-k)}^{2}}
=\frac{\ci(n,n-k-1)}{\ci(n-1,n-k)}=\Ci(n,k+1)&\text{ for }i=n.
\end{cases}
\end{equation}
These correspond to 
\[
\eit(v)=\begin{cases}
[\cd,\ovl 2]&\text{ for }i=1,\\
[\cd,\ovl{i+1},\cd]&\text{ for }1<i<n,\\
[\cd,n,\cd]&\text{ for }i=n.
\end{cases}
\]
Then, in this case we obtain $\eit(m^{(0)}_n(v))=m^{(0)}_n(\eit(v))$.
As for other cases, we can discuss similarly and obtain 
$\eit(m^{(0)}_n(v))=m^{(0)}_n(\eit(v))$ and 
$\fit(m^{(0)}_n(v))=m^{(0)}_n(\fit(v))$.
\qed 

Finally, we shall show 
\begin{pro}\label{til-b}
Let $\wtil\bbB$ be the set of monomials:
\begin{equation}
\wtil\bbB:=
\left\{
\Ci(u_1,1)\Ci(u_2,2)\cd \Ci(u_m,m)
\bCi(l_1,l_1-1)\cdot\bCi(l_2,l_2-2)\cd 
\bCi(l_{k},l_{k}-k)
\begin{array}{|l}
0\leq m\leq n,\,\,k+m=n, \\%m=n -> highest weight vec \ci(n,0)
1\leq l_1< l_2<\cd<l_k\leq n,\\ 
1\leq u_1< u_2<\cd<u_m\leq n.
\end{array}
\right\}
\end{equation}
Then, we get $\bbB=\wtil\bbB$. 
\end{pro}
Here, note that the difference between $\bbB$ and $\wtil\bbB$ is
the condition ``$l_a=u_b\Rightarrow l_a\geq a+b$''.
Thus, the inclusion $\bbB\subset\wtil\bbB$ is trivial.

To show Proposition \ref{til-b}, we see the following:
Let $\cB$ be the set
\begin{equation}
\cB:=\left\{[j_1,\cd,j_n]\,\,\begin{array}{|l}
1\leq  j_1<\cd< j_k\leq \ovl 1.
\end{array}
\right\}.
\end{equation}
Then, it follows from (\ref{Ytab-c}) that $B(\Lm_n)\subset \cB$. 
We extend the map $m^{(0)}_n$ to the map on $\cB$, say also
$m^{(0)}_n$.  So, for 
$v=[u_1,\cd,u_m,\ovl{l_k},\ovl{l_{k-1}},\cd,\ovl{l_1}]\in\cB$,
we have the monomial 
\begin{equation}
m^{(0)}_n([u_1,\cd,u_m,\ovl{l_k},\ovl{l_{k-1}},\cd,\ovl{l_1}])
=\Ci(u_1,1)\Ci(u_2,2)\cd \Ci(u_m,m)
\bCi(l_1,l_1-1)\cdot\bCi(l_2,l_2-2)\cd 
\bCi(l_{k},l_{k}-k),
\label{m0-out}
\end{equation}
which appears in the summation (\ref{Del-w0sn}) and 
may not necessarily belong to $\bbB$. But, indeed, we can show that 
it belongs to $\bbB$, which means that the conjecture is positive for this case.

For $v=[u_1,u_2,\cd,\ovl l_2,\ovl l_1]\in \cB$, if 
$u_a=i=l_b$ and $a+b\leq i$, then we say the pair $(u_a,\ovl l_b)$ in
$v$ is in $i$-{\it configuration} or simply, $v$ is in  
$i$-{\it configuration}. Thus, 
if $v\in \cB$ is in $i$-configuration for any $i$,
 then $v\in B(\Lm_n)$. 
For this $v$, let $(j:k)=j,j+1,\cd,k-1,k$
(resp. $\ovl{(s:t)}=\ovl s,\ovl{s-1},\cd,\ovl t$) be a consecutive subsequence 
of $u_1,\cd,u_m$ (resp. $\ovl l_k,\ovl{l_{k-1}},\cd,\ovl{l_1}$),
which is also called a segment of $v$.

\begin{lem}\label{lem:jmk}
For $v=[u_1,u_2,\cd,\ovl l_2,\ovl l_1]\in \cB$, suppose that 
there exist $a,b$ such that 
$u_a=j=l_b$ and $v$ is not in
 $j$-configuration. 
\begin{enumerate}
\item
For $v=(\cd(j:m)\cd\ovl{(k:j)}\cd)$ we have $n>m+k-j$, where $(j:m)$ and 
$\ovl{(k,j)}$ are segments of $v$.
\item
We assume that $a$ is  the smallest among the elements satisfying
$u_a=l_b=j$ and $a+b>j$ for some $b$.
Then we have $a+b=j+1$.
\end{enumerate}
\end{lem}
{\sl Proof.} 
(i) Since $v$ is not in $j$-configuration, we have $a+b>j$.
If $u_c=m$ and $l_d=k$, $c+d\leq n$. 
We know that $c-a=m-j$ and $d-b=k-j$, which mean that 
\[
 n\geq c+d=(a+m-j)+(b+k-j)>m+k-j.
\]
(ii)
If there is no $c,d$  such that $u_c=l_d$, $c<a$ and
$d<b$.
Then we have $(a-1)+(b-1)\leq j-1$ and then $a+b\leq j+1$. Since $v$ is
not in $j$-configuration, we obtain $a+b=j+1$.
If there are $c,d$ such that $c<a$, $d<b$, $u_c=i=l_d$, $i< j$ and 
$v$ is in $i$-configuration. We may assume that $c,d$ are the nearest
to $a,b$, that is, there is no pair $(k,\ovl k)$ in $v$
such that  $i<k<j$. 
Indeed, if there is such a pair $(k,\ovl k)$ with $i<k<j$
and $v$ is in $k$-configuration, then we may 
replace $i$ with $k$. Unless $v$ is in $k$-configuration, 
it contradicts the minimality of $a,b$.
The last assumption means $(a-c)+(b-d)\leq j-i+1$. 
Since $v$ is in $i$-configuration, we have $c+d\leq i$. 
Then we get $a+b\leq j+1$. Thus, by the assumption $v$ is not in 
$j$-configuration, we have $a+b>j$ and then we obtain $a+b=j+1$.
\qed
%note 16, pp32-33
\begin{df}\label{def-B*}
\begin{enumerate}
\item 
If a pair $(u_a,\ovl l_b)=(i,\ovl i)$ in $v$ satisfies 
$a+b=i+1$, we call such a pair is in $(+1)$-{\it configuration}.
\item 
We define the set as follows:
\begin{equation}
\cB^*:=
\left\{
v=[u_1,\cd,u_m,\ovl l_k,\cd,\ovl l_1]
\begin{array}{|l}
v\text{ satisfies the conditions (A1)-(A6) below.} 
\end{array}
\right\}
\end{equation}
\begin{enumerate}
\item[{(A1)}] $0\leq m\leq n,\,\,k+m=n$.
\item[{(A2)}] There exists $x\in\{0,1,\cd,m\}$ such that 
$1\leq u_1< \cd <u_x$, $u_{x+1}<\cd<u_m\leq n$ and 
there exists $y\in\{0,1,\cd,k\}$ such that 
$1\leq l_1< \cd <l_y$, $l_{y+1}<\cd<l_k\leq n$,  and 
$x,y$ satisfy the following (A3)-(A6).
\item[{(A3)}] If there exist $p\in\{1,\cd,x\}$ and $q\in\{1,\cd,y\}$
	     such that $u_p=l_q$, then $p+q\leq u_p$.
\item[{(A4)}] If $u_x\geq u_{x+1}$, then there exists $z\in\{1,\cd,y\}$
such that $l_z=u_x$ and $x+z\leq u_x$, 
and there exists $w\in\{1,2,\cd,z\}$ such that 
$l_w=u_{x+1}$, $\{l_z,l_{z+1},\cd,l_w\}$ is consecutive 
and $(x+1)+w=u_{x+1}+1$, namely, the pair $(u_{x+1},\ovl l_w)$ is in 
$(+1)$-configuration.
\item[{(A5)}] If $l_y\geq l_{y+1}$, then there exists $z'\in\{1,\cd,x\}$
such that $u_{z'}=l_y$ and $y+z'\leq l_y$,
and there exists $w'\in\{1,2,\cd,z'\}$ such that 
$u_{w'}=l_{y+1}$, $\{u_{z'},u_{z'+1},\cd,u_{w'}\}$ is consecutive 
and $w'+(y+1)=l_{y+1}+1$, namely, the pair $(u_{w'},\ovl l_{y+1})$ is in 
$(+1)$-configuration.
\item[{(A6)}] If $u_x\geq u_{x+1}$ and $l_y\geq l_{y+1}$, then 
$u_x=l_y$.
\end{enumerate}
\end{enumerate}
\end{df}
%See note 17 pp40
By the definition, it is evident $\cB\subset \cB^*$. 
For an element $v=[u_1,\cd,u_m,\ovl l_k,\cd,\ovl l_1]
\in\cB^*$, we define its level $l(v)$ as follows:
\begin{equation}
l(v)=\begin{cases}\min(a+b-1| u_a=l_b
\text{ and }  a+b=u_a+1)&\text{if there exist }a,b \text{ such that }
u_a=l_b, a+b=u_a+1,\\
n&\hbox{otherwise}.
\end{cases}
\end{equation}
The following lemma is obtained from Lemma \ref{lem:jmk} and 
the definition of the level.
\begin{lem}\label{lem:level}
\begin{enumerate}
\item
Assume that for $v\in\cB^*$ there exist $x$ (resp. $y$) such that 
$u_{x}\geq u_{x+1}$ (resp. $l_y\geq l_{y+1}$) and there is no $y$
(resp. $x$) such that $l_y\geq l_{y+1}$ (resp. $u_{x}\geq u_{x+1}$). 
Then we have 
$l(v)=u_{x+1}$ (resp. $l(v)=l_{y+1}$).
\item
Assume that for $v\in\cB^*$ there exist $x, y$ such that 
$u_{x}\geq u_{x+1}$ and $l_y\geq l_{y+1}$. Then we have 
$l(v)=\min(u_{x+1},l_{y+1})$.
\item
For an element $v$ in $\cB^*$, 
$v\in B(\Lm_n)$ if and only if $l(v)=n$.
\end{enumerate}
\end{lem}
{\sl Proof.} The statements (i) and (ii) are evident from the
definition of $\cB^*$. 
Let us show (iii). Write $v=[u_1,\cd,u_k,\ovl l_m,\cd,\ovl l_1]$.
If $v\in B(\Lm_n)$, then any $(i,\ovl i)$-pair in $v$ 
is in $i$-configuration, which means that $l(v)=n$.
Conversely, assume that $v\not\in B(\Lm_n)$. 
Consider the case that for $v$ there exists $x$ as in the condition (A4) above.
In this case, there exists $z$ such that $u_{x+1}=l_z$ and
$(x+1)+z=u_{x+1}+1$,
 which implies $l(v)\leq u_{x+1}<n$. The case that there exists $y$
 satisfying the condition (A5) is treated by the similar way.
Now, assume that there are no $x,y$ satisfying 
the conditions (A4) and (A5) respectively.
In this case, $v$ is an element in $\cB$. It follows from  Lemma
\ref{lem:jmk} that there are $a,b$ such that $u_a=l_b$ and $a+b=u_a+1$.
This shows that $l(v)\leq u_a<n$ 
since $(n,\ovl n)$-pair is always in $n$-configuration.
\qed

\begin{df}\label{tau-def}
Define the transformation $\tau_i$ $(i=1,2,\cd,n-1)$ on $\cB^*$ as follows:
%See note17,pp26-27
\begin{enumerate}
\item
For $v=[u_1,\cd,u_m,\ovl{l_{m'}},\ovl{l_{m'-1}},\cd,\ovl{l_1}]$, 
find $i$ which is the smallest such that $u_a=i=l_b$ for some $a,b$
and $a+b=i+1$.
If it does not exist, then
$\tau_i$ is nothing but the identity. Indeed, In this case 
by Lemma \ref{lem:level} $v\in B(\Lm_n)$. 
Consider the case that such $a$ exists. Then let
$(i:j)$ and $\ovl{(k:i)}$ be the segments in $v$ including $i,\ovl i$
and $j,k$ are the largest ones.
\item
Suppose there are no $x,y$ satisfying the condition (A4) and (A5)
     respectively, that is, $v\in\cB$. 
Let $(j_1:j_2)$ (resp. $\ovl{(k_2:k_1)}$) be the right-next 
(resp. left-next)segments to $(i:j)$ (resp. $\ovl{(k:i)}$), that is, 
$v=(\cd(i:j)(j_1:j_2)\cd\ovl{(k_2:k_1)}\,\,\ovl{(k:i)}\cd)$.
Let 
\[
 \tau_i:v\mapsto(\cd(k+1:k+j-i+1)(j_1:j_2)\cd\ovl{(k_2:k_1)}\,\,
\ovl{(k+j-i+1:j+1)}\cd).
\]
If $j\ne i$, then $\tau_j$ is the identity.
%Here, if $\min(j_1,k_1)>k+j-i+1$, then stop the process and
%     $\tau=\tau_i$.
%Note that though $\tau_i(v)$ is not necessarily in $\cB$, 
%it is easy to see that $\tau_i(v)$ is in $p$-configuration
%with respect to any 
%pair $(p.\ovl p)$ in $(k+1:k+j-i+1)\cd\ovl{(k+j-i+1:j+1)}$.
%Note that by Lemma \ref{lem:jmk} $n\geq k+j-i+1$, which means that 
%$\tau_i$ is well-defined.
\item
Suppose that there is $x$ (resp. $y$)satisfying (A4) (resp. (A5))
and no $y$ (resp. $x$)satisfying (A5) (resp. (A4)). 
In this case, it follows from the above argument that 
set $i=u_{x+1}$ (resp. $i=l_{y+1}$) and then 
$\tau_i$ is defined as same as the previous one
and $\tau_j$ is identity for $j\ne i$.
\item 
Suppose that there are $x$ and $y$ satisfying (A4) and (A5)
respectively and $w$ and $w'$ are as in (A4) and (A5) respectively,
namely, $u_{x+1}=l_z$ and $x+w=u_a$, and $u_{w'}=l_{y+1}$ and $w'+y=l_{y+1}$.
If $u_{x+1}\ne l_{y+1}$, then set $i=\min(u_{x+1},l_{y+1})$ and $\tau_i$ 
is defined as the previous one. 
If $u_{x+1}=l_{y+1}(=i)$, then $\tau_i$ is the composition of the previous
$\tau_{u_{x+1}}$ and $\tau_{l_{y+1}}$, that is, 
$\tau_i:=\tau_{l_{y+1}}\circ \tau_{u_{x+1}}$.
If $j\ne i$, set $\tau_j={\rm id}$.
\end{enumerate}
\end{df}
We obtain the following lemma:
\begin{lem}\label{tau}
\begin{enumerate}
\item
The transformation $\tau_i$ defined above is well-defined, that is, 
$\tau_i(\cB^*)\subset \cB^*$ for any $i$.
\item 
If $\tau_i\ne\id$, then we have $l(\tau_i(v))>l(v)$ for any $v\in
      \cB^*$.
\end{enumerate}
\end{lem}

{\sl Proof.}
First let us see the case (ii) in Definition \ref{tau-def}.
In this case, we have 
\[
  \tau_i(v)=(\cd(k+1:k+j-i+1)(j_1:j_2)\cd\ovl{(k_2:k_1)}\,\,
\ovl{(k+j-i+1:j+1)}\cd).
\]
In this formula, if $k+j-i+1<\min(j_1,k_1)$,
$\tau_i(v)\in\cB\subset\cB^*$. 
We can see that the level of $\tau_i(v)$ is equal to $j_1$ and it 
is greater than $i$, which implies $l(\tau_i(v))>l(v)$.
Consider the case $j_1\leq k+j-i+1$, which corresponds to the
condition (A4) in Definition \ref{def-B*}. By the definition of the segment
we know that $j_1>j+1$. Thus, 
we have $\ovl{j_1}\in\ovl{(k+j-i+1:j+1)}$, 
which means there is the pair $(j_1,\ovl{j_1})$'s in $\tau_i(v)$.
Let us see that this pair is in $(+1)$-configuration.
Set $\tau_i(v):=[u'_1,\cd,u'_k,\ovl l'_m,\cd,\ovl l'_1]$. 
And let $u'_{a}=k+1$, $u'_b=j_1$, $l'_c=\ovl{j+1}$ and $l'_d=\ovl{j_1}$.
Thus, we have $b-a=(k+j-i+1)-(k+1)+1=j-i+1$ and $d-c=j_1-(j+1)$.
Here note that $u_a=i$ and $l_c=\ovl i$ in $v$, which means $a+c=i+1$.
Thus, we get
\[
 b+d=a+(j-i+1)+c+(j_1-j-1)=a+c+j_1-i=j_1+1.
\]
We can easily see that any pair $(i,\ovl i)$ 
in the parts $(u'_1,\cd,{k+j-i+1})$ and
$(\ovl{k+j-i+1},\cd,\ovl l'_c)$ is in $i$-configuration.
Then, $\tau_i(v)$ is in $\cB^*$. 
As for the level of $\tau_i(v)$, we can show as the previous case.
The case $k_1\leq k+j-i+1$ is shown similarly, which corresponds to the 
condition (A5) in Definition \ref{def-B*}.

Next, let us see (iii) in Definition \ref{tau-def}.
Suppose that there is $x$ such that $u_x\geq u_{x+1}$ and 
there is no $y$ such that $l_y\geq l_{y+1}$. 
Let $w$ be the number such that $u_{x+1}=l_w$ and $x+w=u_{x+1}$.
Set $i:=u_x\geq u_{x+1}=:j=l_w$ and let $(j:A)$ (resp. 
$\ovl{(B:j)}$) be the segment including $j=u_{x+1}$ 
(resp. $\ovl j=\ovl l_w$) and $A=u_z$ (resp. $B=l_{z'}$), that is, 
\[
 v=[\cd i,(j:A),A', \cd, \ovl B',\ovl{(B:j)}\cd ]
\]
Note that there is $\ovl i$ between $\ovl B$ and $\ovl j$.
Then, we obtain 
\begin{equation}
\tau_j(v)=[\cd i,(B+1:A+B-j+1),A',\cd,\ovl B',\ovl{(A+B-j+1:A+1)}\cd]
\end{equation}
Here  $i<B+1$ and 
there are the following two cases:\\
(1) $A+B-j+1<A',B'$. (2) $A+B-j+1\geq\min(A',B')$.\\
In the first case, we do not need to consider the conditions (A4)--(A6)
and then $\tau_j(v)\in\cB\subset\cB^*$.
Let us see the case (2). 
%Set $\tau_j(v)=[u'_1,\cd,u'_m,\ovl l'_k,\cd,\ovl l'_1]$. 
By the similar argument to the above cases (ii) and (iii), we know that 
$\tau_j(v)\in\cB^*$ and $l(\tau_j(v))\geq \min(A',B')>j\geq l(v)$.
We can also check the case that there is $y$ such that $l_y\geq l_{y+1}$
and there is no $x$ such that $u_x\geq u_{x+1}$.

Finally, let us see (iv) in Definition \ref{tau-def}.
%Note 17, pp42
Let $x,y$ be numbers satisfying $u_x\geq u_{x+1}$ and $l_y\geq l_{y+1}$
and 
set $i:=u_x=l_y$, $j:=u_{x+1}$, $k:=l_{y+1}$. Note that 
$i\geq j,k$. By the condition (A6), there exist
$r,s$ such that $u_r=k$, $l_s=j$, $(x+1)+s=j+1$ and $(y+1)+r=k+1$.
Note that $l(v)=\min(j,k)$.
Now, $v$ is in the form:
\[
 v=[\cd (k:i),(j:A),B,\cd \ovl C,\ovl{(D: k)},\ovl{(i:j)}\cd ]
\]
where $(j:A)$ (resp. $\ovl{(D:k)}$) is the segment 
including $j$ (resp. $\ovl k$), $A<B+1$ and $C+1>D$.
We consider the following three cases:
(1) $j<k$. (2) $j>k$. (3) $j=k$.\\
The cases (1) and (2) are treated in the similar manner to
the above cases. Thus, we see the case (3).
In this case, $\tau_{j}$ is defined as the last one in 
Definition \ref{tau-def} (iv).
So, we have 
\begin{eqnarray}
&&\qq \tau_j(v)=\tau_{l_{y+1}}\circ \tau_{u_{x+1}}(v) \label{tau-long}
\\
&&\q
=\tau_{l_{y+1}}([\cd,(j:i)(i+1:A+i-j+1),(B:E),F,\cd,\ovl C,\ovl{(D:j)},
\ovl{(A+i-j+1:j+1)},\cd])\nn \\
&&\q=\begin{cases}
[\cd (D+1:G),(B:E),F,\cd,\ovl C,\ovl{(G:A+1)},\cd]
&\text{ if }A+i-j\ne B,\\
[\cd (D+1:H),F,\cd,\ovl C,\ovl{(H:E+1)},\cd]
&\text{ if }A+i-j=B,
\end{cases}\nn
\end{eqnarray}
where $G=A+i-2j+D+1$ and 
$H=E+D-j+1$. 
Here, we can easily see that $\tau_j(v)$ is in $\cB^*$ by
the similar way to the previous cases.
By the formula \eqref{tau-long} we also know that 
\[
 l(v)\leq j(=u_{x+1})<\min(B,C,F)\leq l(\tau_j(v)).
\]
Now, we completed proving the lemma \ref{tau}.\qed

\begin{ex}
For $n=10$ and $v=[12356\ovl6\, \ovl5\,\ovl3\,\ovl2\, \ovl1]\in\cB$,  
we have
\[
[12356\ovl6\, \ovl5\,\ovl3\,\ovl2\, \ovl1]\mapright{\tau_1}
[45656\ovl6\, \ovl5\,\ovl6\,\ovl5\, \ovl4]\mapright{\tau_5}
[45678\ovl6\, \ovl5\,\ovl8\,\ovl7\, \ovl4]
\mapright{\tau_5}[47878\ovl8\, \ovl7\,\ovl8\,\ovl7\, \ovl4]
\mapright{\tau_7}[4789\,10\,\ovl{10}\, \ovl9\,\ovl8\,\ovl7\, \ovl4]
\in B(\Lm_n).
\]
\end{ex}
\begin{pro}\label{pro-tau}
For any element $v\in \cB^*$, there exists
a sequence of indices $i_1,\cd,i_k$ such  that 
$\tau_{i_1}\circ \cd\circ \tau_{i_k}(v)\in B(\Lm_n)$.
\end{pro}
{\sl Proof. }
It follows from Lemma \ref{tau}(ii) 
that for any $v\in \cB^*\setminus B(\Lm_n)$, 
there exists $j$ such that 
$l(v)<l(\tau_j(v))$. This fact means that there exists 
$i_1,\cd,i_k$ such that
\[
 l(\tau_{i_1}\circ \cd\circ \tau_{i_k}(v))=n,
\]
which is equivalent to $\tau_{i_1}\circ \cd\circ \tau_{i_k}(v)
\in B(\Lm_n)$ by Lemma \ref{lem:level}.
\qed

Next lemma is the key for the relations to the monomial realization.
\begin{lem}\label{tau-inv}
For any $v\in\cB^*$ and any $i\in I$, 
we have $m^{(0)}_n(\tau_i(v))=m^{(0)}_n(v)$.
\end{lem}
{\sl Proof.}
If $\tau_i=\id$, there is nothing to show. So, we consider 
an element in $\cB^*\setminus B(\Lm_n)$.
For $v=[u_1,\cd,u_m,\ovl l_m,\cd,\ovl l_1]\in \cB^*\setminus B(\Lm_n)$, 
set  $u_a=i=l_b$, $a+b=i+1$ and let $(i:j)$ (resp. $(\ovl{(k:i)})$) 
be a segment including $i$ (resp. $\ovl i$), that is, 
$v=[\cd,(i:j)\cd,\ovl{(k:i)}\cd]$.
Let 
\[
P:= (\Ci(i,a)\Ci(i+1,a+1)\cd\Ci(j,j-i+a))\cdot
( \bCi(k,{b-i},)\bCi(k-1,b-i)
\cd\bCi(i,i-b))=\frac{\ci(j,n-j+i-a)}{\ci(k,n-a+1)}
\]
be the part of the monomial $m^{(0)}_n(v)$ related to 
the segments $(i:j)$ and $\ovl{(k:i)}$.
Note that for the last equality we use 
$a+b=i+1$.
For the element 
\[
\tau_i(v)=
 [\cd(k+1:k+j-i+1)\cd\ovl{({k+j-i+1}:{j+1})}\cd],
\] 
we also get the part of $m^{(0)}_n(\tau_i(v))$:
\[
 Q=(\Ci(k+1,a)\Ci(k+2,a+1)\cd\Ci(k+j-i+1,j-i+a))\cdot
( \bCi(k+j-i+1,{j-b+1},)\bCi(k+j-i,j-b+1)
\cd\bCi(j+1,j-b+1))=\frac{\ci(j,n-j+i-a))}{\ci(k,n-a+1)}, 
\]
which shows $P=Q$ and then $m^{(0)}_n(v)=m^{(0)}_n(\tau_i(v))$.
\qed

%???????? Uniqueness of rectification ?????????????????
Here, let us see that 
\begin{lem}\label{rect}
For any $v\in\cB^*$ 
there exists a unique element $v'\in B(\Lm_n)$ 
which is obtained by applying $\tau_j$'s to $v$.
\end{lem}
{\sl Proof.}
For $v\in \cB^*\setminus B(\Lm_n)$, suppose that 
there exist $i_1,\cd,i_r$ and 
$j_1,\cd,j_s$  such that 
\[
v_1:= \tau_{i_1}\cd \tau_{i_r}(v)\ne 
\tau_{j_1}\cd\tau_{j_s}(v)=:v_2,
\]
and $v_1,v_2\in B(\Lm_n)$. 
By Lemma \ref{tau-inv} we know that 
\begin{equation}
m_n^{(0)}(v_1)=m_n^{(0)}(v)=m_n^{(0)}(v_2).
\label{vv1v2}
\end{equation}
In the meanwhile, the restricted map ${m_n^{(0)}}_{|B(\Lm_n)}$ 
is bijective. Thus, it follows from  \eqref{vv1v2} that $v_1=v_2$, 
which contradicts the assumption $v_1\ne v_2$.
\qed

Indeed, in the example above, the second step and the third step 
can be exchanged. But, the result turns out to be the same.

For $v\in\cB^*$, we showed that there exists a unique element 
$v'\in B(\Lm_n)$ obtained by applying $\tau_i$'s.
Let us denote $v'$ by ${\rm Rect}(v)$ and call it the 
{\it rectification} of $v$.
Indeed, we know that ${\rm Rect}(\cB^*)={\rm Rect}(\cB)=B(\Lm_n)$.

{\sl Proof }of Proposition  \ref{til-b}. 
We have $m^{(0)}_n(\cB)=\wtil\bbB$ and we know that 
Rect$(\cB)=B(\Lm_n)$ and $m^{(0)}_n({\rm Rect}(v))=m^{(0)}_n(v)$ for any
$v\in\cB$.
Therefore, we obtain 
\[
\q\qq\qq\qq\qq\qq\q\wtil\bbB=m^{(0)}_n(\cB)=m^{(0)}_n({\rm Rect}(\cB))
=m^{(0)}_n(B(\Lm_n))=\bbB.\qq\qq\qq\qq\qq\q\qed
\]

Hence, owing to Theorem \ref{thm-c-2} (ii)  and 
Proposition \ref{til-b} one gets the following result:
\begin{thm}
For $b\in B(\Lm_n)$ let $n(b)$ be the multiplicity of $b$ defined as 
$n(b):=\sharp\{v\in\cB|{\rm Rect}(v)=b\}$. Then, we have
\begin{equation}
\Del_{w_0s_n\Lm_n,\Lm_n}(\Theta^-_\bfii0(c))
=\sum_{b\in B(\Lm_n)\setminus\{u_{\Lm_n}\}}n(b)\cdot{m^{(0)}_n(b)},
\end{equation}
where $u_{\Lm_n}=[1,2,\cd,n-1,n]$ is the highest weight vector.
\end{thm}
Thus, we know that Conjecture \ref{conj1} is affirmative for 
type $C_n$ and the sequence $\bfii0=(12\cd n)^n$.

%%%%%%%%%%%%%%%%%%%%%%%%%%%%%%%%%
\renewcommand{\thesection}{\arabic{section}}
\section{Explicit form of 
$f_B(t\Theta_{\bfii0}^-(c))$ for $B_n$}
\setcounter{equation}{0}
\renewcommand{\theequation}{\thesection.\arabic{equation}}

\subsection{Main theorems}
In case of type $B_n$ fix the sequence 
$\bfii0=(12\cd n)^n$. 
\begin{thm}\label{thm-Bn1}
For $k=1,\cd,n$ and 
$c=(\ci(j,i))_{1\leq,i,j\leq n}=
(\ci(1,1),\ci(2,1),\cd,\ci(n-1,n),\ci(n,n))\in(\bbC^\times)^{n^2}$, we have
\begin{eqnarray*}
&&\Del_{w_0\Lm_k,s_k\Lm_k}(\Theta^-_\bfii0(c))\\
&&\qq=
\ci(1,k)+\frac{\ci(2,k)}{\ci(1,k+1)}+\cd+
\frac{\ci(n-1,k)}{\ci(n-2,k+1)}+\frac{{\ci(n,k)}^2}{\ci(n-1,k+1)}
+2\frac{\ci(n,k)}{\ci(n,k+1)}+\frac{\ci(n-1,k+1)}{{\ci(n,k+1)}^2}
+\frac{\ci(n-2,k+2)}{\ci(n-1,k+2)}+\cd+
\frac{\ci(k,n)}{\ci(k+1,n)},
\end{eqnarray*}
where note that 
$\Del_{w_0\Lm_n,s_n\Lm_n}(\Theta^-_\bfii0(c))=\ci(n,n).$
%Note 15,pp29
\end{thm}
\begin{thm}\label{thm-Bn2}
Let $k$ be an index running over $\{1,2,\cd,n-1\}$ and $c$ be as in the
previous theorem.
Then we have
\begin{eqnarray}
\Del_{w_0s_k\Lm_k,\Lm_k}(\Theta^-_\bfii0(c))
=\frac{1}{\ci(1,k)}+\sum_{j=1}^{k-1}\frac{\ci(k-j,j)}{\ci(k-j+1,j)}.
\end{eqnarray}
\end{thm}
The case $k=n$ will be presented in \ref{Bn-del-n}.
We shall prove the above theorems in \ref{proof-b-1} and 
\ref{subsec-thm-Bn2}.
%%%%%%%%%%%%%%%%%%%%%%%%%%%%%%%%%%%%%%%%%
\subsection{Proof of Theorem \ref{thm-Bn1}}\label{proof-b-1}

Considering similarly to type $C_n$ as in \ref{proof-c-1}, 
we can write 
\begin{eqnarray*}
{\pmb x}_i(c):=\al_i^\vee(c^{-1})x_i(c)&=&
\begin{cases}c^{-h_i}(1+c\cdot e_i)&i\ne n,\\
c^{-h_n}(1+c\cdot e_n+\frac{c^2}{2}e_n^2)&i=n,
\end{cases}\\
{\pmb y}_i(c):=y_i(c)\al_i^\vee(c^{-1})&=&
\begin{cases}(1+c\cdot f_i)c^{-h_i}&i\ne  n,\\
(1+c\cdot f_n+\frac{c^2}{2}f_n^2)c^{-h_n}&i=n,
\end{cases}
\end{eqnarray*}
since $f_i^2=e_i^2=0$ ($i\ne n$) and $e_n^3=f_n^3=0$
on the vector representation $V(\Lm_1)$.
We also have $\omega({\pmb y}_i(a))=\al_i^\vee(a^{-1})x_i(a)=
{\pmb x}_i(a)$ and define $\xxi(\ovl i,p,j)=\xxi(\ovl i,p,j)(c^{(p)})$ 
and $\xxi(i,p,j)=\xxi(i,p,j)(c^{(p)})$ for
$p,j\in I$ by 
\begin{eqnarray*}
&&X^{(p)}X^{({p-1})}\cd X^{(1)}v_{i}
=\sum_{j=0}^n{\xxi(i,p,j)}v_j
+\sum_{j=1}^n{\xxi(i,p,\ovl j)}v_{\ovl j}\in V(\Lm_1)
\q(i=0,1,\cd,n),\\
&&X^{(p)}X^{({p-1})}\cd X^{(1)}v_{\ovl i}
=\sum_{j=0}^n{\xxi(\ovl i,p,j)}v_j
+\sum_{j=1}^n{\xxi(\ovl i,p,\ovl j)}v_{\ovl j}\in V(\Lm_1)
\q(i=1,2,\cd,n),
\end{eqnarray*}
where $c^{(p)}=(\ci(1,1),\ci(2,1),\cd,\ci(n-1,p),\ci(n,p))$
and $X^{(p)}={\pmb x}_n(\ci(n,p)){\pmb x}_{n-1}(\ci(n-1,p))\cd 
{\pmb x}_1(\ci(1,p))$.

By (\ref{minor-bilin}) and $\omega(\Theta_{{\bf i}_0}(c))=X^{(n)}\cd
X^{(1)}$, 
same as \eqref{del-x} we have 
\begin{equation*}
 \Del_{w_0\Lm_i,s_i\Lm_i}(\Theta_\bfi(c))=\lan \ovl s_i\cdot
 u_{\Lm_i}\,,\,
X^{(n)}\cd X^{(1)}v_{\Lm_i}\ran.
\end{equation*}

To describe $\xxi(\ovl i,p,j)$ explicitly
let us use the {\it segments} as in Sect.\ref{cn}.

For $m\in M=M_1\sqcup\cd\sqcup M_S\in{\mathcal M}^{(p)}_k$, define 
$n(m):=n-j+1$ if $m\in M_j$.
For $M=M_1\sqcup\cd\sqcup M_S\in{\mathcal M}^{(p)}_k$, write
$M_1=\{2,3,\cd,a\}$. 
For 
$i-1\leq i_2\leq i_3\leq\cd\leq i_a\leq n$ and for 
$1\leq b<c\leq a$ and $i\leq j_2\leq \cd\leq j_b\leq n$, 
define the monomials in $(\ci(j,i))_{1\leq i,j\leq n}$ by 
\begin{eqnarray}
&&C_{i_2,i_3,\cd,i_a}^M:=\frac{
\left(\ci(i_2+1-2\ep_{i_2},1+\ep_{i_2})\right)^{1+\ep_{i_2+1}}\cd
\left(\ci(i_a+1-2\ep_{i_a},a+\ep_{i_a}-1)\right)^{1+\ep_{i_a+1}}}
{{\ci(i_2,2)}^{1+\ep_{i_2}}\cd {\ci(i_a,a)}^{1+\ep_{i_a}}},\qq
D^M:=\prod_{m\in M\setminus M_1}\frac{\ci(n(m)-1,m)}{\ci(n(m),m)},\\
&&
\wtil C^M_{j_2,j_3,\cd,j_b;b,c}:=\frac{
{\ci(j_2,1)}^{1+\ep_{j_2}}\cd {\ci(j_b,b-1)}^{1+\ep_{j_b}}
\ci(n,b)\ci(n-1,c+1)\cd \ci(n-1,a)}
{\ci(j_2-1,2)\cd \ci(j_b-1,b)
\ci(n,c){\ci(n,c+1)}^2\cd {\ci(n,a)}^2}.
\label{c-bar}
\end{eqnarray}
%Note15 pp44, Note17 pp1
where $\ep_i=\del_{i,n}$ and $C_{i_2,i_3,\cd,i_a}^M=1$ 
(resp. $D^M=1$) if $M_1=\emptyset$ (resp. $M\setminus M_1=\emptyset$).
Note that in (\ref{c-bar}) if $b=1$ (resp. $c=a$), 
then ${\ci(j_2,1)}^{1+\ep_{j_2}}\cd {\ci(j_b,b-1)}^{1+\ep_{j_b}}
=\ci(j_2-1,2)\cd \ci(j_b-1,b)=1$ (resp. 
$\ci(n-1,c+1)\cd \ci(n-1,a)=
{\ci(n,c+1)}^2\cd {\ci(n,a)}^2=1$).
\begin{pro}
In the setting above, we have
\begin{eqnarray}
&&\xxi(\ovl i,p,\ovl k)={\ci(i-1,1)}^{-1}
\sum_{i\leq i_2\leq \cd\leq i_p\leq k}
(\ci(i_2-1,2)\cd \ci(i_p-1,p))^{-1}
({\ci(i_2,1)}^{1+\ep_{i_2}}{\ci(i_3,2)}^{1+\ep_{i_3}}\cd
{\ci(i_p,p-1)}^{1+\ep_{i_p}}{\ci(k,p)}^{1+\ep_k}).
%See note 15 pp39 or note 16 pp99
\label{b-xi-bar}\\
&& \xxi(\ovl i,p,0)={\ci(i-1,1)}^{-1}\sum_{q=1}^p
\sum_{i\leq i_2\leq \cd\leq i_q\leq n}
(\ci(i_2-1,2)\cd \ci(i_q-1,q))^{-1}
({\ci(i_2,1)}^{1+\ep_{i_2}}{\ci(i_3,2)}^{1+\ep_{i_3}}\cd
{\ci(i_q,q-1)}^{1+\ep_{i_q}}{\ci(n,q)}),
\label{b-xi-0}
%note 15 pp39 or 16 pp99
\\
&& \xxi(\ovl i,p,k)
={\ci(i-1,1)}^{-1}\left(
\sum_{(A)}
C_{i_2,i_3,\cd,i_a}^M\cdot D^M
+2\sum_{(B)}
\wtil C_{j_2,\cd,j_b;b,c}^M\cdot D^M
\right)
\label{b-xi}
%See note 17 pp1
\end{eqnarray}
where $\ci(0,1)=1$ and the conditions  (A)  and (B) are as follows:
\begin{enumerate}
\item[(A)]
$i-1\leq i_2\leq\cd\leq i_a\leq n$, 
$M=M_1\sqcup\cd\sqcup M_S\in{\mathcal M}^{(p)}_k,$
$M_1=\{2,\cd,a\}$.
\item[(B)]
$i\leq j_2\leq\cd\leq j_b\leq n, $
$M=M_1\sqcup\cd\sqcup M_S\in{\mathcal M}^{(p)}_k,$
$M_1=\{2,\cd,a\}\ne\emptyset,\,\,1\leq b<c\leq a.$
\end{enumerate}
\end{pro}

{\sl Proof.}
Set $\cX:={\pmb x}_n(c_n)\cd{\pmb x}_1(c_1)$. Calculating directly
we have the formula:
\begin{eqnarray}
&&\cX v_i=\begin{cases}c_1^{-1}v_1&\text{ if }i=1,\\
c_{i-1}c_i^{-1}v_i+v_{i-1}&\text{ if }i=2,\cd,n,
\end{cases}\\
&&\cX v_0=v_0+2c_n^{-1}v_n,\\
&&\cX v_{\ovl i}=
c_{i-1}^{-1}(c_iv_{\ovl i}+c_{i+1}v_{\ovl{i+1}}+\cd +c_n^2v_{\ovl n}+
c_nv_0+v_n),
\end{eqnarray}
where we understand $c_0=1$. Using these, we get for $i=1,2,\cd,n$, 
\begin{eqnarray}
&&\xxi(\ovl i,p,\ovl k)=\sum_{j=i}^k\xxi(\ovl i,p-1,\ovl j)
\frac{{\ci(k,p)}^{1+\ep_k}}{\ci(j-1,p)},
\q(k=1,\cd,n),
\label{recur-b1}\\
&&
\xxi(\ovl i,p,k)=\xxi(\ovl i,p-1,k+1)+\xxi(\ovl i,p-1,k)
\frac{\ci(k-1,p)}{\ci(k,p)},
\q(k=1\cd,n-1),
\label{recur-b2}\\
&&\xxi(\ovl i,p,0)=\sum_{j=i}^n\xxi(\ovl i,p-1,\ovl j)
\frac{\ci(n,p)}{\ci(j-1,p)}+\xxi(\ovl i,p-1,0),
%See note 15 pp39
\label{recur-b3}\\
&&\xxi(\ovl i,p,n)=\sum_{j=i}^n\xxi(\ovl i,p-1,\ovl j)
{\ci(j-1,p)}^{-1}+2 \xxi(\ovl i,p-1,0)
{\ci(n,p)}^{-1}+\xxi(\ovl i,p-1,n)\frac{\ci(n-1,p)}{{\ci(n,p)}^2}.
\label{recur-b4}
%note 15 pp40 
\end{eqnarray}
Indeed, the formulae (\ref{b-xi-bar}) and 
(\ref{b-xi-0})
are easily shown by the induction
on $p$ using the formulae (\ref{recur-b1}) and (\ref{recur-b3}).

To obtain (\ref{b-xi}) we see the segments of elements in $\cM_{k}^{(p)}$, 
$\cM_{k+1}^{(p-1)}$ and $\cM_{k}^{(p-1)}$ as the case $C_n$ and 
apply the recursions (\ref{recur-b2}) and (\ref{recur-b4}) to the
induction hypothesis. Arguing similarly to the previous case, 
we obtain the desired results.\qed

Thus, for example, we have 
\[
\Del_{w_0\Lm_1,s_1\Lm_1}=\xxi(\ovl 1,n,2)
=\sum_{j=1}^{n-1}\frac{\ci(j,1)}{\ci(j-1,2)}+
\frac{{\ci(n,1)}^2}{\ci(n-1,2)}+2\frac{\ci(n,1)}{\ci(n,2)}+
\frac{\ci(n-1,2)}{{\ci(n,2)}^2}+
\sum_{j=3}^n\frac{\ci(n-j+1,j)}{\ci(n-j+2,j)}.
\]

The following is the same as Lemma \ref{lem-det-c}.
\begin{lem}\label{lem-det-b}
For $k=1,\cd, n-1$ we define the matrix $W_k$ by 
\begin{equation}
W_k:=
\begin{pmatrix}
\xxi(\ovl 1,n,k+1)&\xxi(\ovl 1,n,k-1)&\cd&\xxi(\ovl 1,n,2)&
\xxi(\ovl 1,n,1)\\
\xxi(\ovl 2,n,k+1)&\xxi(\ovl 2,n,k-1)&\cd&\xxi(\ovl 2,n,2)&
\xxi(\ovl 2,n,1)\\
\vdots&\vdots&\cd&\vdots&\vdots\\
\xxi(\ovl k,n,k+1)&\xxi(\ovl k,n,k-1)&\cd&\xxi(\ovl k,n,2)&
\xxi(\ovl k,n,1)
\end{pmatrix}.
\label{det-b}
\end{equation}
Then we have $\Del_{w_0\Lm_k,s_k\Lm_k}(\Theta^-_\bfii0(c))=\det W_k$.
\end{lem}

The last column of the matrix $W_k$ is just 
\[
{^t}(\xxi(\ovl 1,n,1),\cd,\xxi(\ovl k,n,1))=
 {}^t(1,{\ci(1,1)}^{-1},{\ci(2,1)}^{-1},\cd,{\ci(k-1,1)}^{-1}).
\]
Then, applying the elementary transformations on $W_k$ by 
($i$-th row) $-\frac{\ci(i,1)}{\ci(i-1,1)}\times$($i+1$-th row)
for $i=1,\cd,k-1$, in the transformed matrix $\til W_k$
its $(i,j)$-entry is as follows:
\begin{lem} For $c=(\ci(k,l))_{1\leq,k,l\leq n}$ 
we set $c^{(l)}:=(\ci(1,l),\cd,\ci(n,l))$ and 
$c^{[a:b]}:=(c^{(a)},c^{(a+1)},\cd,c^{(b)})$ $(a\leq b)$.
For $i=1,\cd, k-1$  the $(i,j)$-entry $(\til W_k)_{i,j}$ is:
\begin{equation}
(\til W_k)_{i,j}=\begin{cases}
\xxi(\ovl i,n-1,k+1)(c^{[2:n]})&\text{ if }j=1,\\
\xxi(\ovl i,n-1,k-j+1)(c^{[2:n]})&\text{ if }1<j<k,\\
0&\text{ if }j=k,
\end{cases}
\end{equation}
\end{lem}
Note that by definition we have $c=c^{[1:n]}$.\\
{\sl Proof.}
The proof is similar to the one for Lemma \ref{lem-w-c}. Indeed, 
for type $B_n$ we also get the same formula as \eqref{xi-xi}:
\begin{equation*}
\xxi(\ovl i,n,j)(c)-\frac{\ci(i,1)}{\ci(i-1,1)}\xxi(\ovl{i+1},n,j)(c)
=\frac{\ci(i,1)}{\ci(i-1,1)} \xxi(\ovl i,n-1,j)(c^{[2,n]}).
\end{equation*}
Then, this shows the lemma.
\qed

Applying the above elementary transformations to the matrix $W_k$, we
find 
\begin{eqnarray*}
\det W_k
=\xxi(\ovl 1,n-k+1,k+1)(c^{[k:n]}).
\end{eqnarray*}
Thus, it follows from  (\ref{b-xi})  that 
\begin{eqnarray}
&&\Del_{w_0\Lm_k,s_k\Lm_k}(\Theta^-_\bfi(c))
=\xxi(\ovl 1,n-k+1,k+1)(c^{[k:n]})\\
&&\qq=\ci(1,k)+\frac{\ci(2,k)}{\ci(1,k+1)}+\cd+
\frac{\ci(n-1,k)}{\ci(n-2,k+1)}+\frac{{\ci(n,k+1)}^2}{\ci(n-1,k+1)}+
2\frac{\ci(n,k)}{\ci(n,k+1)}+
\frac{\ci(n-1,k+1)}{{\ci(n,k+1)}^2}
+\frac{\ci(n-2,k+2)}{\ci(n-1,k+2)}+\cd+
\frac{\ci(k,n)}{\ci(k+1,n)}.\nn
\end{eqnarray}
The case $k=n$ is easily obtained by the formula in \cite[(4.18)]{BZ2}:
\begin{equation}
\Del_{w_0\Lm_n,s_k\Lm_n}(\Theta^-_\bfii0(c))=\ci(n,n).
\end{equation}
Now, the proof of Theorem \ref{thm-Bn1} has been accomplished.
\qed
%%%%%%%%%%%%%%%%%%%%%%%%%%%%%%%%%%%%%%%%%
\subsection{Proof of Theorem \ref{thm-Bn2}}\label{subsec-thm-Bn2}%8.3

The proof of Theorem \ref{thm-Bn2} is the same as that of 
Theorem \ref{thm-c-2} (ii). Indeed, we obtain the formula 
(\ref{vmat-c2}) and (\ref{xi-xi2}) for type $B_n$. 
Defining the matrix $U_k$ and $Z_k$, same as type $C_n$
we have $\det Z_j=1$ ($j=1,\cd,k-1$)and $\det U_1={\ci(1,k)}^{-1}$.
Thus, we obtain the desired result for type $B_n$.\qed

%%%%%%%%%%%%%%%%%%%%%%%%%%%%%%%%%%%%%%%%%%%
\subsection{Correspondence to the monomial realizations}%8.4

Except for $\Del_{w_0s_n\Lm_n,\Lm_n}(\Theta^-_\bfii0(c))$ 
(see \ref{Bn-del-n}), we shall see
the positive answer to the conjecture for type $B_n$.

First, we see the monomial realization of $B(\Lm_1)$ associated with 
the cyclic oder $\cd(12\cd n)(12\cd n)\cd$, which means that the sign
$p_0=(p_{i,j})$ is given by $p_{i,j}=1$ if $i<j$ and $p_{i,j}=0$ if 
$i>j$ as in the previous section.
The crystal $B(\Lm_1)$ is described as follows:
We abuse the notation 
$B(\Lm_1):=\{v_i,v_{\ovl i}|1\leq i\leq n\}\sqcup\{v_0\}$ if 
there is no confusion.
Then the actions of $\eit$ and $\fit$ $(1\leq i<n)$ are defined as 
$\fit=f_i$ and $\eit=e_i$ in (\ref{b-f1}).
The actions of $\til f_n$ and $\til e_n$ are given as:
\begin{equation}
\til f_nv_n=v_0,\q \til f_nv_0=v_{\ovl n},\q 
\til e_nv_0=v_n,\q \til e_nv_{\ovl n}=v_0.
\label{b-til-f2}
\end{equation}

To see the monomial realization $B(\ci(1,k))$, 
we describe the monomials $A_{i,m}$ explicitly:
% See Note 17 pp4-5
\begin{equation}
A_{i,m}=\begin{cases}
\ci(1,m){\ci(2,m)}^{-1}\ci(1,m+1)&\text{ for }i=1,\\
\ci(i,m){\ci(i+1,m)}^{-1}{\ci(i-1,m+1)}^{-1}\ci(i,m+1)
&\text{ for }1<i< n-1,\\
\ci(n-1,m){\ci(n,m)}^{-2}{\ci(n-2,m+1)}^{-1}
\ci(n-1,m+1)&\text{ for }i=n-1,\\
\ci(n,m){\ci(n-1,m+1)}^{-1}\ci(n,m+1)&\text{ for }i=n.
\end{cases}
\label{aim-b-1}
\end{equation}
Here the monomial realization $B(\ci(1,k))$ for 
 $B(\Lm_1)$ associated with 
$p_0$ is described explicitly:
\begin{equation}
B(\ci(1,k))=\left\{\frac{{\ci(j,k)}^{\ep_j}}{\ci(j-1,k+1)}=m^{(k)}_1(v_j),
\frac{\ci(n,k)}{\ci(n,k+1)}=m^{(k)}_1(v_0),
\frac{\ci(j-1,k+n-j+1)}{{\ci(j,k+n-j+1)}^{\ep_j}}=m^{(k)}_1(v_{\ovl j})\,\,
|\,\,1\leq j\leq n\,\,\right\},
\end{equation}
where 
$m^{(k)}_i:B(\Lm_i)\hookrightarrow \cY(p)$ ($u_{\Lm_i}\mapsto 
\ci(i,k)$)
is the embedding of crystal as in Sect.5
and we understand $\ci(0,k)=1$.
Now, Theorem \ref{thm-Bn1} and Theorem \ref{thm-Bn2} mean the following:
\begin{thm}
We obtain $\Del_{w_0\Lm_n,s_n\Lm_n}(\Theta^-_\bfii0(c))
=m_n^{(n)}(u_{\Lm_n})(=\ci(n,n))$
 and 
\begin{eqnarray}
&&\Del_{w_0\Lm_k,s_k\Lm_k}(\Theta^-_\bfii0(c))
=\sum_{j=1}^nm^{(k)}_1(v_j)+2m^{(k)}_1(v_0)
+\sum_{j=k+1}^nm^{(k)}_1(v_{\ovl j}),\qq
(k=1,2\cd,n-1)
\label{thm-mono-b1}\\
&&\Del_{w_0s_k\Lm_k,\Lm_k}(\Theta^-_\bfii0(c))
=\sum_{j=1}^km^{(k-n)}_1(v_{\ovl j}), \qq(k=1,2,\cd,n-1).
\label{thm-mono-b2}
\end{eqnarray}
\end{thm}
Note that the second result is derived from the fact that
$B(\ci(1,k))=B({\ci(1,n+k)}^{-1})$, which is the connected component
including ${\ci(1,n+k)}^{-1}$ as the lowest monomial.

%%%%%%%%%%%%%%%%%%%%%%%%%%%%%%%%%%%%%%%%%%
\subsection{Triangles and $\Del_{w_0s_n\Lm_n,\Lm_n}(\Theta^-_\bfii0(c))$}
\label{Bn-del-n}%8.5
To state the result for $\Del_{w_0s_n\Lm_n,\Lm_n}(\Theta^-_\bfii0(c))$,
we need certain preparations.
The set of {\it triangles} $\tri_n$ is defined as follows:
\begin{equation}
\tri_n:=\{{(\ji(k,l))}_{1\leq k\leq l\leq n}|
1\leq \ji(k,l+1)\leq \ji(k,l)<\ji(k+1,l+1)\leq n+1
\q(1\leq k\leq l<n)\}.
\end{equation}
We visualize a triangle $(\ji(k,l))$ in $\tri_n$ as follows:
\[(\ji(k,l))=
 \begin{array}{c}
\ji(1,1)\\
\ji(2,2)\ji(1,2)\\
\ji(3,3)\ji(2,3)\ji(1,3)\\
\cd\cd\cd\cd\\
\ji(n,n)\ji(n-1,n)\cd\ji(2,n)\ji(1,n)
\end{array}
\]
By the definition of $\tri_n$, we easily obtain 
\begin{lem}
For any $k\in\{1,2,\cd,n\}$ there exists a unique $j$
($1\leq j\leq k+1$) such that the $k$th row of 
a triangle $(\ji(k,l))$ in $\tri_n$ is in the following form:
\begin{equation}
k\text{-th row}\q
(\ji(k,k),\ji(k-1,k),\cd,\ji(2,k),\ji(1,k))
=(k+1, k,k-1,\cd,j+1,j-1,j-2,\cd,2,1),
\end{equation}
that is, we have $\ji(m,k)=m$ for $m<j$ and 
$\ji(m,k)=m+1$ for $m\geq j$.
\end{lem}
For a triangle $\del=(\ji(k,l))$, we list $j$'s as in the lemma:
$s(\del):=(s_1,s_2,\cd,s_n)$, which we call the {\it label} of 
a triangle $\del$.
Here we have
\begin{lem}\label{s4}
For $\del\in\tri_n$ let $s(\del)=(s_1,\cd,s_n)$ be its label. 
Then, we have
\begin{enumerate}
\item
The label $s(\del)$  satisfies 
$1\leq s_k\leq k+1$,
and $s_{k+1}=s_{k}$ or $s_{k}+1$ for $k=1,\cd,n$.
\item
Each  $k$-th row of a triangle $\del$ is in 
one of the following {\rm I, II, III, IV:}
\begin{enumerate}[\rm I.]
\item $s_{k+1}=s_k+1$ and $s_k=s_{k-1}$.
\item $s_{k+1}=s_k$ and $s_k=s_{k-1}$.
\item $s_{k+1}=s_k+1$ and $s_k=s_{k-1}+1$.
\item $s_{k+1}=s_k$ and $s_k=s_{k-1}+1$.
\end{enumerate}
Here we suppose that $s_0=1$ and $s_{n+1}=s_{n-1}+1$, which means that 
the 1st row must be in {\rm I,II or IV} 
and the $n$-th row must be in {\rm I or IV}.
\end{enumerate}
\end{lem}
Now, we associate a Laurant 
monomial $m(\del)$ in variables $(\ci(i,j))_{i\in I,j\in\bbZ}$ 
with a triangle $\del=(\ji(k,l))$ by the
following recipe.
\begin{enumerate}
\item Let $s=(s_1\cd,s_n)$  be the label of $\del$.
\item If  $i$-th row is in the form I, then associate $\ci(i,s_i)$.
\item If  $i$-th row is in the form IV, then associate
      ${\ci(i,s_i)}^{-1}$.
\item If  $i$-th row is in the form II or III, then associate 1.
\item  Take the product of all monomials as above for $1\leq i \leq n$, then
we obtain the monomial $m(\del)$ associated with $\del$. This defines a map 
$m:\tri_n\to \cY$, where $\cY$ is the set of Laurant monomials 
in $(\ci(i,j))_{i\in I,j\in\bbZ}$.
\end{enumerate}
Let us denote the special triangle $\del=(\ji(k,l))$ such that 
$\ji(k,l)=k+1$(resp. $\ji(k,l)=k$) for any $k,l$ by $\del_h$ 
(reps. $\del_l$). 
Then, we can present the result for type $B_n$.
\begin{thm}\label{thm-Lm-n} For the 
type $B_n$, we obtain the explicit form:
\begin{equation}
\Del_{w_0s_n\Lm_n,\Lm_n}(\Theta^-_\bfii0(c))=
\sum_{\del\in\tri_n\setminus\{\del_l\}}\ovl{m(\del)},
\end{equation}
where the monomial $\ovl{m(\del)}$ is obtained by applying 
$\ovl{\,\,}:\ci(i,j)\to{\ci(i,n-j+1)}^{-1}$. 
\end{thm}
The proof of this theorem will be given in \ref{proof-thm-Lm-n}.
\begin{ex} 
The set of triangles $\tri_4$ is as follows:
\begin{eqnarray*}
&&\begin{array}{c}
1\\21\\321\\4321
\end{array}\q
\begin{array}{c}
1\\21\\321\\5321
\end{array}\q
\begin{array}{c}
1\\21\\421\\5321
\end{array}\q
\begin{array}{c}
1\\31\\421\\5321
\end{array}\q
\begin{array}{c}
2\\31\\421\\5321
\end{array}\q
\begin{array}{c}
2\\31\\421\\5421
\end{array}\q
\begin{array}{c}
1\\31\\421\\5421
\end{array}\q
\begin{array}{c}
1\\21\\421\\5421
\end{array}\\
&&\begin{array}{c}
1\\31\\431\\5421
\end{array}\q
\begin{array}{c}
2\\31\\431\\5421
\end{array}\q
\begin{array}{c}
2\\32\\431\\5421
\end{array}\q
\begin{array}{c}
1\\31\\431\\5431
\end{array}\q
\begin{array}{c}
2\\31\\431\\5431
\end{array}\q
\begin{array}{c}
2\\32\\431\\5431
\end{array}\q
\begin{array}{c}
2\\32\\432\\5431
\end{array}\q
\begin{array}{c}
2\\32\\432\\5432
\end{array}
\end{eqnarray*}
and their labels $s(\del)$ are
\begin{eqnarray*}
&&(2,3,4,5),\q (2,3,4,4),\q(2,3,3,4),\q(2,2,3,4),\q(1,2,3,4),\q
(1,2,3,3),\q(2,2,3,3),\q(2,3,3,3),\\
&&(2,2,2,3),\q (1,2,2,3),\q(1,1,2,3),\q(2,2,2,2),\q(1,2,2,2),\q
(1,1,2,2),\q(1,1,1,2),\q(1,1,1,1).
\end{eqnarray*}
Then, we have the corresponding monomials ${m(\del)}$:
\begin{eqnarray*}
&& \frac{1}{\ci(4,5)},\q
\frac{\ci(4,4)}{\ci(3,4)},\q \frac{\ci(3,3)}{\ci(2,3)\ci(4,4)},\q
\frac{\ci(2,2)}{\ci(1,2)\ci(4,4)},\q\frac{\ci(1,1)}{\ci(4,4)},\q
\frac{\ci(1,1)\ci(4,3)}{\ci(3,3)},\q\frac{\ci(2,2)\ci(4,3)}{\ci(1,2)\ci(3,3)},\q
\frac{\ci(4,3)}{\ci(2,3)},\\
&&\frac{\ci(3,2)}{\ci(1,2)\ci(4,3)},\q
\frac{\ci(1,1)\ci(3,2)}{\ci(2,2)\ci(4,3)},\q\frac{\ci(2,1)}{\ci(4,3)},\q
\frac{\ci(4,2)}{\ci(1,2)},\q \frac{\ci(1,1)\ci(4,2)}{\ci(2,2)},\q
\frac{\ci(2,1)\ci(4,2)}{\ci(3,2)},\q \frac{\ci(3,1)}{\ci(4,2)},\q
{\ci(4,1)},
\end{eqnarray*}
and we have the monomials $\ovl{m(\del)}$:
\begin{eqnarray*}
&& \ci(4,0),\q
\frac{\ci(3,1)}{\ci(4,1)},\q \frac{\ci(2,2)\ci(4,1)}{\ci(3,2)},\q
\frac{\ci(1,3)\ci(4,1)}{\ci(2,3)},\q\frac{\ci(4,1)}{\ci(1,4)},\q
\frac{\ci(3,2)}{\ci(1,4)\ci(4,2)},\q\frac{\ci(1,3)\ci(3,2)}{\ci(2,3)\ci(4,2)},\q
\frac{\ci(2,2)}{\ci(4,2)},\\
&&\frac{\ci(1,3)\ci(4,2)}{\ci(3,3)},\q
\frac{\ci(2,3)\ci(4,2)}{\ci(1,4)\ci(3,3)},\q\frac{\ci(4,2)}{\ci(2,4)},\q
\frac{\ci(1,3)}{\ci(4,3)},\q \frac{\ci(2,3)}{\ci(1,4)\ci(4,3)},\q
\frac{\ci(3,3)}{\ci(2,4)\ci(4,3)},\q \frac{\ci(4,3)}{\ci(3,4)},\q
\frac{1}{\ci(4,4)}.
\end{eqnarray*}
Then, the total sum of all monomials $\ovl{m(\del)}$ except $\ci(4,0)$
is $\Del_{w_0s_4\Lm_4,\Lm_4}(\Theta^-_\bfii0(c))$ for $B_4$.
\end{ex}

%%%%%%%%%%%%%%%%%%%%%%%%%%%%%%%%%%%%%%
\subsection{Crystal structure on $\tri_n$}%8.6

Let $\del$ be a triangle. 
Then the actions of $\fit$ and $\eit$ are defined as follows:
Let $J_k=(\ji(k,k),\ji(k-1,k),\cd,\ji(2,k),\ji(1,k))$ be the $k$-th
row of $\del$. Thus, we denote $\del=(J_1,\cd,J_n)$.
It follows from  Lemma \ref{s4} that
there exists a unique $j$ such that 
$J_k=(k+1,k,\cd,j+1,j-1,\cd,2,1)$ and $J_k$ is in one of I,II,III,IV.
Set $J'_k=(k+1,\cd,j+2,j,j-1,\cd,2,1)$ and 
$J''_k=(k+1,\cd,j+1,j,j-2,\cd,2,1)$.
Then, we have
\begin{eqnarray}
&&\fit\del=\begin{cases}
(\cd,J_{i-1},J'_i,J_{i+1},\cd)&\text{ if $J_i$ is in I},\\
\qq 0&\text{otherwise,}
\end{cases}\label{Cr-del-f}\\
&&\eit\del=\begin{cases}
(\cd,J_{i-1},J''_i,J_{i+1},\cd)&\text{ if $J_i$ is in IV},\\
\qq 0&\text{otherwise.}
\end{cases}\label{Cr-del-e}
\end{eqnarray}
The weight of $\del=(J_1,\cd,J_n)$ is defined as follows:
Let $s=(s_k)_{k=1,\cd,n}$ be the label of $\del$, that is, 
$J_k=(k+1,k,\cd, j+1,j-1,\cd,2,1)$ for $j=s_k$:
\begin{equation}
\wt(\del)=\Lm_n-\sum_{k=1}^{n}(s_k-1)\al_k.
\label{Cr-del-wt}
\end{equation}
We can easily check that $\tri_n$ is equipped with 
the crystal structure by 
\eqref{Cr-del-f}, \eqref{Cr-del-e} and 
\eqref{Cr-del-wt}, and obtain:
\begin{pro}\label{del-spin-b}
As a crystal, $\tri_n$ is isomorphic to $B(\Lm_n)$.
The highest (resp. lowest) weight crystal is 
$\del_h\hbox{(resp. $\del_l$)}\in\tri_n$.
\end{pro}
{\sl Proof.}
As was given in \ref{Bn}, we know the explicit form of the crystal 
$B(\Lm_n)\cong B^{(n)}_{sp}$. So, let us describe the 
one-to-one correspondence between $\tri_n$ and $B^{(n)}_{sp}$:
Let $s=(s_1,\cd,s_n)$ be the label of a triangle $\del\in \tri_n$.
Now, let us associate $\ep=(\ep_1,\cd,\ep_n)\in B^{(n)}_{sp}$ with $s$ by
\begin{equation}
\ep_1=\begin{cases}+&\text{ if }s_1=1,\\
-&\text{ if }s_1=2,
\end{cases}\qq
\ep_k=\begin{cases}+&\text{ if }s_{k}=s_{k-1},\\
-&\text{ if }s_{k}=s_{k-1}+1,
\end{cases}(k>1),
\end{equation}
which defines the map $\psi:\tri_n\to B^{(n)}_{sp}$. 
Then, e.g., the vector $\del_h$ (resp. $\del_l$) corresponds to the 
highest (resp. lowest) weight vector $(+,+,\cd,+)$
(resp. $(-,-,\cd,-)$).
It is clear to find that the map $\psi$ is bijective.
Now, if $k$-th row of $\del$ is in type I, then
$(s_{k-1},s_k,s_{k+1})=(j,j,j+1)$
 for some $j$ and in the corresponding $\ep=\psi(\del)$, we have
 $\ep_k=+,\ep_{k+1}=-$. For $\til f_k(\del)$ let $s'$ be its label
and $\ep':=\psi(\til f_k\del)$.
It is easy to see that  $(s'_{k-1},s'_k,s'_{k+1})=(j,j+1,j+1)$ and $\ep'_k=-$,
$\ep'_{k+1}=+$, which shows 
that the map $\psi$ is compatible with the
action of $\til f_k$, that is, 
$\ep'=\psi(\til f_k\del)=\til f_k\ep=\til f_k\psi(\del)$.
For the types II,III and IV we can see $\til f_k(\del)=0$ and 
$\til f_k \psi(\del)=0$, thus we find the compatibility of $\psi$.
By arguing similarly, we can also see the compatibility of $\psi$ with
the action of $\til e_k$'s. 
Thus, we find that the map $\psi$ is an isomorphism of crystals.\qed

Next, let us show that the map $m:\tri_n\to\cY$ gives an 
isomorphism between $\tri_n$ and the monomial
realization of $B(\Lm_n)$.
Consider the crystal structure on $\cY$
by taking $\ovl p=(\ovl p_{i,j})$ such that 
\begin{equation}
\ovl p_{i,j}=\begin{cases}1&\text{ if }i>j,\\
0&\text{ if }i<j.
\end{cases}
\label{ovl-pij}
\end{equation}
Indeed, this corresponds to the cyclic sequence of $I$ such as:
$\cd(n n-1 \cd 21)(n n-1 \cd 21)\cd$.
So, by the prescription in Sect.\ref{sect-mono} 
we obtain the monomial realization $\cY(\ovl p)$. 
Here, by the definition of the map $m$ we can get that 
$m(\del_h)=\ci(n,1)$, which is one of the highest monomials in
$\cY(\ovl p)$ with the highest weight $\Lm_n$. 
\begin{pro}\label{mono-del-b}
Let $B(\ci(n,1))$ be the connected subcrystal of $\cY(\ovl p)$ whose
highest monomial is $\ci(n,1)$ and which is isomorphic to
 $B(\Lm_n)$
of type $B_n$. Then we have $m(\tri_n)= B(\ci(n,1))$.
\end{pro}
{\sl Proof.}
In the setting $\cY(\ovl p)$, we have 
\begin{equation}\label{aim}
A_{i,m}=\begin{cases}
\ci(1,m){\ci(2,m+1)}^{-1}\ci(1,m+1),&\text{ for }i=1,\\
\ci(i,m){\ci(i-1,m)}^{-1}{\ci(i+1,m+1)}^{-1}\ci(i,m+1),
&\text{ for }1<i<n-1,\\
\ci(n-1,m){\ci(n-2,m)}^{-1}{\ci(n,m+1)}^{-2}\ci(n-1,m+1),&\text{ for }i=n-1,\\
\ci(n,m){\ci(n-1,m)}^{-1}\ci(n,m+1)&\text{ for }i=n.
\end{cases}
\end{equation}
For $\del\in\tri_n$ and its label $s=(s_1,\cd,s_n)$, 
suppose the $k$-th row of $\del$ is in I $(k\ne1,n)$, that is, 
there is some $j$ such that $(s_{k-1},s_k,s_{k+1})=(j,j,j+1)$,
where as we mentioned above $s_0=1$ and $s_{n+1}=s_{n-1}+1$.
In this case, the $k-1$(resp. $k+1$)-th row is in II or IV 
(resp. III or IV).
In the case $k\ne,1,n-1,n$, 
by the action of $\til f_k$ we have $(s_{k-1},s_k,s_{k+1})=
(j,j,j+1)\to(j,j+1,j+1)$. Then we get the $k$-th row of $\til f_k(\del)$
is in IV and the $k-1$ (resp. $k+1$)-th row is in I or III 
(resp. I or II), i.e., the following four cases:
\[
\til f_k: (s_{k-2},s_{k-1},s_k,s_{k+1},s_{k+2})=
\begin{cases}(j,j,j,j+1,j+1)&\to(j,j,j+1,j+1,j+1),\\
(j-1,j,j,j+1,j+1)&\to(j-1,j,j+1,j+1,j+1),\\
(j,j,j,j+1,j+2)&\to(j,j,j+1,j+1,j+2),\\
(j-1,j,j,j+1,j+2)&\to(j-1,j,j+1,j+1,j+2).
\end{cases}
\]
In these cases, the types of the $k-1,k,k+1$-th rows and 
the parts of monomials related to the 
variables $\ci(k-1,l),\ci(k,l),\ci(k+1,l)$
are turn out to be
as follows:
\[
\til f_k:\begin{array}{cccccc}
\rm (II,I,IV)&\to&\rm(I,IV,II)\qq&
\frac{\ci(k,j)}{\ci(k+1,j+1)}&\to&\frac{\ci(k-1,j)}{\ci(k,j+1)}\\
\rm (IV,I,IV)&\to&\rm(III,IV,II)\qq&
\frac{\ci(k,j)}{\ci(k-1,j)\ci(k+1,j+1)}&\to&\frac{1}{\ci(k,j+1)}\\
\rm (II,I,III)&\to&\rm(I,IV,I)\qq&
\ci(k,j)&\to&\frac{\ci(k-1,j)\ci(k+1,j+1)}{\ci(k,j+1)}\\
\rm (IV,I,III)&\to&\rm(III,IV,I)\qq&
\frac{\ci(k,j)}{\ci(k-1,j)}&\to&\frac{\ci(k+1,j+1)}{\ci(k,j+1)}.
\end{array}
\]
In all these cases, the action by $\til f_k$ on monomial $m(\del)$
is described as multiplying the
monomial
$A_{k,j}^{-1}={\ci(k,j)}^{-1}\ci(k-1,j)\ci(k+1,j+1){\ci(k,j+1)}^{-1}$, 
which implies $m(\til f_k\del)=\til f_k m(\del)$.

In the case $k=1$, if the 1st row is in I, then we get 
$(s_0,s_1,s_2)=(1,1,2)$. This is changed by the action of $\til
f_1$ to $(1,2,2)$, which means that as for the labels of 
the 0th,1st,2nd,3rd rows and the parts of the monomials
related to variables $\ci(1,1),\ci(1,2),\ci(2,2)$. Hence, we have 
\[\til f_1:\q
\begin{array}{ccccccc}
\rm(I,IV)&\to&\rm(IV,II)&\qq
 &\frac{\ci(1,1)}{\ci(2,2)}&\to&\frac{1}{\ci(1,2)},\\
\rm(I,III)&\to&\rm(IV,I)&\qq &\ci(1,1)&\to&\frac{\ci(2,2)}{\ci(1,2)}.
\end{array}
\]
In both cases, the action by $\til f_1$ is given
by multiplying the monomials 
${A_{1,1}}^{-1}={\ci(1,1)}^{-1}\ci(2,2){\ci(2,1)}^{-1}$.

The cases $k=n$ are also done by the similar way:
\[
 \til f_{n}:(s_{n-2},s_{n-1},s_n)=
\begin{array}{ccccccc}
 &(j-1,j,j)&\to&(j-1,j,j+1)&\q \ci(n-1,j)\ci(n,j)&\to&\frac{1}{\ci(n,j+1)}\\
 &(j,j,j)&\to&(j,j,j+1)&\q \ci(n,j)&\to&\frac{\ci(n-1,j)}{\ci(n,j+1)}.
\end{array}
\]
This means that the action by $\til f_{n}$ is given by multiplying the
monomial  ${A_{n,j}}^{-1}={\ci(n,j)}^{-1}\ci(n-1,j){\ci(n,j+1)}^{-1}$.
Let us see the last case $k=n-1$.  If the $i$-th row of $\del$ is in I,
then we have 
$(s_{n-3},s_{n-2},s_{n-1},s_n)=(j,j,j,j+1)$ or $(j-1,j,j,j+1)$ whose
types are (II,I,IV) or (IV,I,IV) respectively if the $(n-1)$-th row is
in I.
Then we have
\[
 \til f_{n-1}:\q
\begin{array}{ccccccc}
&\rm(II,I,IV)&\to&\rm(I,IV,I)&\qq\frac{\ci(n-1,j)}{\ci(n,j+1)}&\to
&\frac{\ci(n-2,j)\ci(n,j+1)}{\ci(n-1,j+1)},\\
&\rm(IV,I,IV)&\to&\rm(III,IV,I)&\qq\frac{\ci(n-1,j)}{\ci(n-2,j)\ci(n,j+1)}&
\to&\frac{\ci(n,j+1)}{\ci(n-1,j+1)}.
\end{array}
\]
This implies that the action by $\til f_{n-1}$ is described by
multiplying the monomial \\
${A_{n-1,j}}^{-1}=
{\ci(n-1,j)}^{-1}\ci(n-2,j){\ci(n,j+1)}^2{\ci(n-1,j+1)}^{-1}.
$
As for the actions of $\eit$, we can show that $m(\eit \del)=\eit
m(\del)$ by the similar manner to the cases of $\fit$. \qed

Then, by these results we have 
\begin{thm}
$B(\ci(n,1))=m(\tri_n))\cong B(\Lm_n)$ for type $B_n$.
\end{thm}

%%%%%%%%%%%%%%%%%%%%%%%%%%%%%%%%%%%
\subsection{Proof of Theorem \ref{thm-Lm-n}}\label{proof-thm-Lm-n}

By the explicit descriptions in (\ref{B-fi}) and (\ref{B-ei}),
we know that $e_i^2=f_i^2=0$ on $V^{(n)}_{sp}$.
Thus, we can write 
${\pmb x}_i(c):=\al_i^\vee(c^{-1})x_i(c)=c^{-h_i}(1+c\cdot e_i)$ 
and ${\pmb y}_i(c):=y_i(c)\al_i^\vee(c^{-1})=(1+c\cdot f_i)c^{-h_i}$ 
on 
$V^{(n)}_{sp}$
and then for $\ep=(\ep_1,\cd,\ep_n)\in B^{(n)}_{sp}$
\begin{equation}
\pmbx_i(c)\ep=\begin{cases}
c\ep+\ep'&\text{ if }
(\ep_i,\ep_{i+1})=(-,+),\\
c^{-1}\ep&\text{ if }(\ep_i,\ep_{i+1})=(+,-),\\
\ep&\text{ otherwise},
\end{cases}\,\,(i<n),\qq
\pmbx_n(c)\ep=\begin{cases}c\ep+\ep''&\text{ if }\ep_n=-,\\
c^{-1}\ep&\text{ if }\ep_n=+,
\end{cases}
\label{xi-Bn}
\end{equation}
where $\ep'=(\cd,\ep_{i-1},+,-,\ep_{i+2},\cd)$ for $i<n$ and 
$\ep''=(\cd,\ep_{n-1},+)$.

For $c^{(i)}=(\ci(1,i),\cd,\ci(n,i))\in (\bbC^\times)^n$ set
$X^{(i)}(c^{(i)}):=\pmbx_1(\ci(1,i))\cd \pmbx_n(\ci(n,i))$ and 
$X(c):=X^{(n)}(c^{(n)})\cd X^{(1)}(c^{(1)})$ where 
$c:=(c^{(n)},\cd,c^{(1)})$.

To calculate 
$\Del_{w_0s_n\Lm_n,\Lm_n}(\Theta^-_\bfii0(c))$ explicitly, first
let us see 
$\Del_{w_0\Lm_n,s_n\Lm_n}(\Theta^-_{\bfii0^{-1}}(c))$
($\bfii0^{-1}=(n\,\, n-1\cd 21)^n$) 
since
for $g\in U\ovl w_0 U$ we have 
\begin{equation}\label{eta-si}
\Del_{w_0s_i\Lm_i,\Lm_i}(g)=\Del_{w_0\Lm_i,s_i\Lm_i}(\eta(g)),
\end{equation}
%See Note 16 p90-91.
and $\eta(\Theta^-_{\bfii0}(c))=\Theta^-_{{\bfii0}^{-1}}(\ovl c)$
($-:\ci(j,i)\mapsto\ci(n-j+1,i)$).
By the above formula and $\omega(\Theta^-_{\bfii0^{-1}}(c))=X(c)$, 
we have
\[
\Del_{w_0\Lm_n,s_n\Lm_n}(\Theta^-_{\bfii0^{-1}}(c))
=\lan s_n\ep_h,X(c)\ep_l\ran
\] 
where $s_n\ep_h=s_n(+,\cd,+)=(+,\cd,+,-)$ and $\ep_l=(-,\cd,-)$.
Thus, we would like to see the coefficient of the vector 
$(+,\cd,+,-)$ in $X(c)\ep_l$.

For the sequence $\bfii0=(12\cd n)^n$, 
define $J$ to be the set of all subsequences of $\bfii0$.
Since each $\pmbx_i(c)$ can be written in the form $c^{h_i}(1+ce_i)$ on 
$V_{sp}$,
the operator $X(c)$ is expanded in the form:
\begin{equation}\label{expand}
X(c)=\sum_{\xi\in J}\al_\xi e_\xi,
\end{equation}
%Note 16 pp48
where $e_\xi:=e_{i_1}\cd e_{i_m}$ for 
$\xi=(i_1,\cd,i_m)\in J$ and $\al_\xi$ 
is a product of $\al_i^\vee(c)$'s
and some scalars. 
Now, we shall find for which subsequence $(i_1,\cd,i_k)\in J$ we get
\begin{equation}\label{non0}
 e_{i_1}\cd e_{i_k}\ep_l=(+,\cd,+,-).
\end{equation}
Here note that 
\begin{equation}\label{wt-dif}
\wt(+,\cd,+,-)-\wt(\ep_l)=\sum_{i=1}^ni\al_i-\al_n,
\end{equation}
and then for $(i_1,\cd,i_k)$ in (\ref{non0}), we know that
if $i\ne n$ the number of $i$'s is $i$ and the number of $n$ is 
$n-1$. 

Let us associate an element in $U({\mathfrak n}_+)$ with 
a triangle $\del=(\ji(k,l))\in \tri_n$ by the following way:
\begin{enumerate}
\item
For $i=1,\cd,n$ define $r_i(\del)$ be the set of indices as
$r_i(\del):=\{(k,l)|\ji(k,l)=i\}$ and 
$n_i(\del):=\sharp r_i(\del)$.
\item
Set $L=n_i(\del)$. 
For $r_i(\del)=\{(k_{a_1},l_{a_1}),\cd,(k_{a_{L}},l_{a_L})\}$
define  $l_{a_1},\cd,l_{a_L}$ with $l_{a_1}<l_{a_2}<\cd<l_{a_L}$.
\item For $\del$ let $r_i(\del)$ be as in (ii).  We set 
\[
 E_i(\del):=e_{l_{a_1}}e_{l_{a_2}}\cd e_{l_{a_L}}
\]
and $E(\del):=E_n(\del)\cdot E_{n-1}(\del)\cd E_1(\del)$.
\end{enumerate}
\begin{ex}\label{EX-B}
For $n=4$ and the triangle 
\[
\del=\begin{array}{c}
2\\32\\431\\5421
\end{array}\qq
\text{we obtain }\q
\begin{array}{ll}
r_1(\del)=\{(1,3),(1,4)\},&E_1(\del)=e_3e_4,\\
r_2(\del)=\{(1,1),(1,2),(2,4)\},&E_2(\del)=e_1e_2e_4\\
r_3(\del)=\{(2,2),(2,3)\},&E_3(\del)=e_2e_3\\
r_4(\del)=\{(3,3),(3,4)\},&E_4(\del)=e_3e_4,
\end{array}
\]
and then $E(\del)=(e_3e_4)(e_2e_3)(e_1e_2e_4)(e_3e_4)$.
\end{ex}

\begin{lem}\label{tri-ep}
Let $(j_1,\cd,j_m)\in J$ be a subsequence of $\bfii0$ 
and $E=e_{j_1}\cd e_{j_m}$ be the associated
monomial of $e_i$'s. Then we have 
\begin{equation}
E\ep_l=(+,\cd,+,-)\text{ if and only if there exists } 
\del\in \tri_n\setminus\{\del_l\}\text{ such that }
E=E(\del).
\end{equation}
\end{lem}
Note that $\del=(\ji(k,l))\ne\del_l$ if and only if $\ji(n,n)=n+1$.

\nd
{\sl Proof.}
First, for $\del=(\ji(k,l))=\del_l$, that is, $\ji(k,l)=k$, we have
\[
 E(\del_l)=e_n(e_{n-1}e_n)(e_{n-2}e_{n-1}e_n)\cd(e_1e_2\cd e_{n-1}e_{n}),
\]
and then 
$E(\del_l)\ep_l=(+,+,\cd,+,+)\ne(+,+,\cd,+,-)$. 

Next, assuming that there exists $\del\in\tri_n\setminus\{\del_l\}$ 
such that 
$E=E(\del)$, show that $E\ep_l=s_n\ep_h=(+,\cd,+,-)$.
Indeed, if $E\del\ne0$, we have that 
$E\ep_l=s_n\ep_h$ by the explicit actions of $e_i$'s 
and the weight counting.
So it suffices to show $E\ep_l\ne0$.
Now, for $\del=(\ji(i,k))$ write 
\[
E= E(\del)=e_{l_m}\cd e_{l_2}e_{l_1}\q(m=\frac{n(n+1)}{2}-1).
\]
Let us show $e_{l_k}\cd e_{l_1}\ep_l\ne0$ 
($1\leq k\leq m$) by the induction on $k$.
By the definition of $E(\del)$ we know that $e_{l_1}=e_n$ and then 
this implies $e_{l_1}\ep_l=e_n\ep_l\ne0$.
Assume $e_{l_{k-1}}\cd e_{l_1}\ep_l\ne0$ and  $e_{l_k}$ 
appears in $E_j(\del)$. Set $i=l_k$.
If $i\leq n-1$, then the $i-1,i,i+1$th row of the triangle $\del$ 
around this $j$ is in one of the following forms:
\begin{eqnarray*}
&&(a) 
\begin{array}{cc}
\syl i-1\text{-th row}&\syl *\q j-1\q\q\\
\syl i\text{-th row}&\q \textcircled{$\scriptstyle j$}\q \syl j-1\q \\
\syl i+1\text{-th row}&\syl \q\q j\q\q j-1
\end{array}\qq
(b)
\begin{array}{cc}
\syl i-1\text{-th row}&\syl *\q j-1\q\q\\
\syl i\text{-th row}&\q \textcircled{$\scriptstyle j$}\q \syl j-1\q \\
\syl i+1\text{-th row}&\syl \q\q j\q\q j-2
\end{array}\qq
(c)
\begin{array}{cc}
\syl i-1\text{-th row}&\syl *\q j-1\q\q\\
\syl i\text{-th row}&\q \textcircled{$\scriptstyle j$}\q \syl j-2\q \\
\syl i+1\text{-th row}&\syl \q\q j\q\q j-2
\end{array}\\
&&(d) 
\begin{array}{cc}
\syl i-1\text{-th row}&\syl *\q j-2\q\q\\
\syl i\text{-th row}&\q \textcircled{$\scriptstyle j$}\q \syl j-2\q \\
\syl i+1\text{-th row}&\syl \q\q j\q\q j-2
\end{array}\qq
(e)
\begin{array}{cc}
\syl i-1\text{-th row}&\syl *\q j-1\q\q\\
\syl i\text{-th row}&\q \textcircled{$\scriptstyle j$}\q \syl j-2\q \\
\syl i+1\text{-th row}&\syl \q\q j-1\q **
\end{array}\qq
(f)
\begin{array}{cc}
\syl i-1\text{-th row}&\syl *\q j-2\q\q\\
\syl i\text{-th row}&\q \textcircled{$\scriptstyle j$}\q \syl j-2\q \\
\syl i+1\text{-th row}&\syl \q\q j-1\q **
\end{array},
\end{eqnarray*}
where $*$ means $j$ or $j+1$, $**$ means $j-2$ or $j-3$ and 
$\textcircled{$\syl j$}$ is the entry which gives $e_{l_k}=e_i$ 
in $E(\del)$.
As for (a), we see that $E_j=\cd e_ie_{i+1}\cd$ and 
$E_{j-1}=\cd e_{i-1}e_ie_{i+1}\cd$, and then
 we have $e_{l_{k-1}}=e_{i+1}$.
Suppose that $j-1$ in the $i+1$-th row 
gives $e_{i+1}=e_{l_p}$ $(1\leq p\leq m)$ in $E_{j-1}(\del)$. 
Since $e_{l_{k-1}}\cd e_{l_1}\ep_l\ne0$, we know that 
the vector before applying $e_{i+1}$ in $E_{j-1}$, say $\ep'=(\ep'_i)$, 
is not zero, that is, 
\[
 \ep'=e_{l_{p-1}}\cd e_{l_1}\ep_l\ne0.
\]
Since the vector $e_{i-1}e_ie_{i+1}\ep'\ne0$, 
the vector $\ep'$ satisfies $\ep'_{i-1}=\ep'_i=\ep'_{i+1}=-,\ep'_{i+2}=+$. 
If $i\leq n-2$, there exists $j$ in $i+2$-th row of $\del$ 
and this means that there is $e_{i+2}$ in $E_j(\del)$ since 
there are $j-1$ and $j$ in the $i+1$th row. Then, 
we know that the vector before being applied $e_{l_{k-1}}=e_{i+1}$, say 
$\ep''=(\ep''_i)$, has the signs $\ep''_i=\ep''_{i+1}=-$ 
and $\ep''_{i+2}=+$, which implies $e_ie_{i+1}\ep''\ne0$:
\begin{eqnarray*}
&& (\ep'_{i-1},\ep'_i,\ep'_{i+1},\ep'_{i+2})=(-,-,-,+)\mapright{e_{i+1}}
(-,-,+,-)\mapright{e_i}(-,+,-,-)\mapright{e_{i-1}}
(+,-,-,-)\\
&&\cd \mapright{e_{i+2}}(*,-,-,+)\mapright{e_{i+1}}
(*,-,+,-)\mapright{e_i}(*,+,-,-).
\end{eqnarray*}
%See Note 17,pp30
Thus, we know that $e_{l_k}\cd e_{l_1}\ep_l\ne 0$ in the case (a). 
For the case (b), since we can deduce that there is only $j-1$ or $j$ 
in the $i+2$th row, the following two cases can occur:
\[
 E_{j}(\del)E_{j-1}(\del)
=(\cd e_ie_{i+1}\cd)(\cd e_{i-1}e_ie_{i+2}\cd), \text{or }
 E_{j}(\del)E_{j-1}(\del)
=(\cd e_ie_{i+1}e_{i+2}\cd)(\cd e_{i-1}e_i\cd),
\]
where in the R.H.S. of the above equations 
``$\cd$'' means that there is no
$e_{i-1},e_{i},e_{i+1}$ or $e_{i+2}$.
In both cases, before applying $e_i$ in $E_j(\del)$, we applied 
$e_{i-1}$ in $E_{j-1}(\del)$ and then $e_{i+1}$ in $E_j(\del)$, 
and we know that 
the resulting vector $\ep$
does not vanish by the induction hypothesis and is in
the form $\ep=(\ep_k)$ such that $\ep_i=-$ and $\ep_{i+1}=+$, which
means  $e_i\ep\ne0$.

For the cases (c)--(f), arguing similarly we obtain 
\begin{enumerate}
\item[(c)]
$E_j(\del)\cdot E_{j-1}(\del)=(\cd e_ie_{i+1}\cd)(\cd e_{i-1}e_{i+2}\cd)$ or 
$(\cd e_ie_{i+1}e_{i+2}\cd)(\cd e_{i-1}\cd)$.
\item[(d)]
$E_j(\del)\cdot E_{j-1}(\del)=(\cd e_ie_{i+1}\cd)(\cd e_{i+2}\cd)$ or 
$(\cd e_ie_{i+1}e_{i+2}\cd)(\cd)$.
\item[(e)]
$E_j(\del)\cdot E_{j-1}(\del)=(\cd e_i\cd)(\cd
	  e_{i-1}e_{i+1}e_{i+2}\cd)$
or $(\cd e_i\cd)(\cd e_{i-1}e_{i+1}\cd)$.
\item[(f)]
$E_j(\del)\cdot E_{j-1}(\del)=(\cd e_i\cd)(\cd e_{i+1}e_{i+2}\cd)$ or 
$(\cd e_i\cd)(\cd e_{i+1}\cd)$.
\end{enumerate}
Thus, in all cases we find that
$e_{l_k}\cd e_{l_1}\ep_l\ne 0$. 
The case $i=n-1$ is done by the similar
way.  So, finally we assume that $i=l_k=n$. 
Then, the $n-1$-th row and the $n$-th row of $\del$ around $j$ are as
follows:
\[
 (i)
\begin{array}{cc}
\syl n-1\text{-th row}&\syl *\q j-1\q\\
\syl n\text{-th row}&\q \textcircled{$\scriptstyle j$}\q \syl j-1
\end{array}\qq
(ii)
\begin{array}{cc}
\syl n-1\text{-th row}&\syl *\q j-1\q\\
\syl n\text{-th row}&\q \textcircled{$\scriptstyle j$}\q \syl j-2
\end{array}\qq
(iii)
\begin{array}{cc}
\syl n-1\text{-th row}&\syl *\q j-2\q\\
\syl n\text{-th row}&\q \textcircled{$\scriptstyle j$}\q \syl j-2
\end{array}
\]
where $*$ means $j$ or $j+1$ and 
$\textcircled{$\syl j$}$ is the entry which gives $e_{l_k}=e_n$. 
As for (i), we know that $E_j(\del)=\cd e_n\cd$ and 
$E_{j-1}(\del)=\cd e_{n-1}e_n\cd$.
Suppose that $j-1$ in the $n$-th row gives $e_n=e_{l_p}$ 
in $E_{j-1}(\del)$.
Denote the vector before being applied 
$e_n=e_{l_p}$ by $\ep'=(\ep'_i)$, that is, 
$\ep'=e_{l_{p-1}}\cd e_{l_1}\ep_l$, 
which satisfies $\ep'_{n-1}=\ep'_{n}=-$ and then 
\[
E_{j-1}(\del):
 (\ep'_{n-1},\ep'_n)=(-,-)\mapright{e_n}(-,+)\mapright{e_{n-1}}(+,-).
\]
This means that the vector before being applied $e_n=e_{l_k}$, 
say $\ep''=(\ep''_i)$, has the form $\ep''_n=-$ and then 
we have $e_n\ep''=e_{l_k}\cd e_{l_1}\ep_l\ne 0$.
Considering similarly, we have 
\begin{enumerate}
\item[(ii)] $E_j(\del)\cdot E_{j-1}(\del)\cdot E_{j-2}(\del)
=(\cd e_n\cd)(\cd e_{n-1}\cd)(\cd e_n\cd)$.
\item[(iii)] $E_j(\del)\cdot E_{j-1}(\del)\cdot E_{j-2}(\del)
=(\cd e_n\cd)(\cd \cd)(\cd e_{n-1}e_n\cd)$.
\end{enumerate}
These imply that in the both cases the vector applied $e_n$ in 
$E_j(\del)$ never vanish. Thus, we also find that 
$e_n\ep''=e_{l_k}\cd e_{l_1}\ep_l\ne 0$ and then 
$E(\del)\ep_l\ne0$. Now, we find that $E(\del)\ep_l\ne0$ for any 
$\del\in \Del_n\setminus\{\del_h\}$.

Next, assume that $E\ep_l=(-,-,\cd,-,+)$.
Under this assumption, by (\ref{wt-dif}) we have that each 
$e_i$ ($i\ne n$) should appear  $i$-times in $E$ and 
$e_n$ should appear $n-1$-times in $E$.
In the sequence $\bfii0$ there are $n$ cycles $(12\cd n)$ just as
\[
\begin{array}{ccccc}
\bfii0= \underbrace{12\cd n}& \underbrace{12\cd n}&
\cd  &\underbrace{12\cd n}&\underbrace{12\cd n}\\
n\text{th cycle}&n-1\text{th cycle}&&2\text{nd cycle}&1
\text{st cycle}
\end{array}
\]
Since there are $n$ positions of the index $n$ in $\bfii0$,
we should choose $n-1$ from them, which is written as 
$(n,n-1,\cd,j+1,j-1,\cd,2,1)$ for some $j$.
This defines 
$(\ji(n-1,n),\ji(n-2,n),\cd,\ji(2,n),\ji(1,n))$, which is 
the bottom row($=n$-th row) of an element in $\tri_n$
except the first entry $\ji(n,n)$. Here, we fix $\ji(n,n)=n+1$.
Next, let us see $e_{n-1}$.
We can apply $e_{n-1}$ after applying $e_n$ and can apply 
$e_n$ after applying $e_{n-1}$. 
Thus, we get the positions of $n-1$'s in $\bfii0$, 
which must be in between two $n$'s chosen in the previous step.
Then, we obtain the list
$(\ji(n-1,n-1),\ji(n-2,n-1),\cd,\ji(2,n-1),\ji(1,n-1))$, which 
should satisfy the condition $\ji(k,n)\leq \ji(k,n-1)<\ji(k+1,n)$
$(1\leq k\leq n-1)$ and becomes the $n-1$th 
row of the corresponding $\del$.
This condition just coincides with the ones for $\tri_n$.
Repeating these steps, we obtain a triangle $\del$ 
and find that there exists $\del\in\tri_n$ such that $E=E(\del)$.
\qed

To show Theorem \ref{thm-Lm-n},
let us find the coefficient of the vector 
generated from $\ep_l$ by 
applying $a_\xi e_{i_1}\cd e_{i_k}$ 
in (\ref{expand}).
Set ${\bf E}_i(c)=c\al_i^\vee(c^{-1})e_i$. Then we can write
$\pmbx_i(c)=\al_i^{\vee}(c^{-1})+{\bf E}_i(c)$.
Here note that for $\ep=(\ep_i)\in B_{sp}^{(n)}$:
\[
{\bf E}_i(c)\ep=c\al_i^\vee(c^{-1})e_i(\ep_1,\cd,\ep_n)=\begin{cases}
(\cd,\buildrel{i}\over+,\buildrel{i+1}\over-,\cd)
&\text{ if }\ep_i=-,\,\,\ep_{i+1}=+,\,\,
i\ne n,\\
(\cd\cd\cd,\buildrel{n}\over+)
&\text{ if }\ep_n=-,\,\,i=n,\\
\qq0&\text{otherwise.}
\end{cases}
\]
We consider the four cases I-IV: If the $i$-th row of $\del$ is in I,
there exists $j$ such that the $i-1$, $i$, $i+1$th rows of $\del$
are in the form:
\[
\begin{array}{lcccc}
i-1\text{-th row}&\cd j+3\q j+2\q j+1\q j-1\q j-2\\
i\text{-th row}&\cd \q\q j+3\q j+2\q j+1\q j-1\q j-2\\
i+1\text{-th row}&\cd\q\q\q\q j+3\q j+2\q j\q\q j-1\q j-2
\end{array}
\]
This means that in the expansion (\ref{expand}) we know that 
this $\del$ gives the following monomial in which $\al_i^{\vee}$
appears once:
\begin{equation}
\cd \al_i^\vee{({\ci(i,j)}^{-1})}{\bf E}_{i+1}(\ci(i+1,j))\cd
{\bf E}_{i-1}(\ci(i-1,j-1)){\bf E}_{i}(\ci(i,j-1))
{\bf E}_{i+1}(\ci(i+1,j-1))\cd 
\end{equation}
Thus, by the action of $\al_i^\vee{({\ci(i,j)}^{-1})}$
we obtain the coefficient $\ci(i,j)$.
In  the other cases II,III,IV, discussing similarly
we obtain the coefficient 1 for II and III, and we have 
the coefficient ${\ci(i,j)}^{-1}$ for IV, which coincides with the recipe
after Lemma \ref{s4} and then we know that 
the desired coefficient is the same as $m(\del)$.

As has been shown above, the set of monomials 
$m(\tri_n)$ has the crystal structure 
isomorphic to $B(\ci(n,1)) (\cong B(\Lm_n))$.
The following lemma is direct from Theorem \ref{thm-dual}.
\begin{lem}
The set of monomials $\ovl{m(\tri_n)}$ is a crystal isomorphic to
$B(\ci(n,0))\cong B(\Lm_n)$.
\end{lem}
{\sl Proof.}
This is the case $a=n+1$ in Theorem \ref{thm-dual}.
Indeed, the explicit form of $A_{i,l}$ is given in (\ref{aim}), 
which induces the formula
$\ovl{A_{i,l}^{-1}}=A_{i,n-l}$.\qed

Thus, by the formula (\ref{eta-si}), we obtain 
\[
 \Del_{w_0s_n\Lm_n,\Lm_n}(\Theta^-_\bfii0(c))=
\Del_{w_0\Lm_n,s_n\Lm_n}(\Theta^-_{\bfii0^{-1}}(\ovl c)),
\]
where for $c=(c^{(1)},\cd,c^{(n)})$ ($c^{(i)}=(\ci(l,i))_{1\leq l\leq n}$)
we define $\ovl c=(\ovl{c^{(n)}},\cd,\ovl{c^{(1)}})$ 
$(\ovl{c^{(i)}}=({\ci(n-l+1,i)}^{-1})_{1\leq l\leq n})$.
Therefore, we have completed the proof of 
Theorem \ref{thm-Lm-n}.
\qed 

Hence, we obtain the affirmative answer to our conjecture for type
$B_n$.

%%%%%%%%%%%%%%%%%%%%%%%%%%%%%%%%%%%%%%%%%
\renewcommand{\thesection}{\arabic{section}}
\section{Explicit form of 
$f_B(t\Theta_{\bfii0}^-(c))$ for $D_n$}
\setcounter{equation}{0}
\renewcommand{\theequation}{\thesection.\arabic{equation}}

\subsection{Main Theorems}
In case of type $D_n$, fix the cyclic reduced longest word
$\bfii0=(1 2 \cd n-1\,n)^{n-1}$.
\begin{thm}\label{thm-Dn1}
For $k\in\{1,2,\cd,n\}$ and $c=(\ci(i,j))=(\ci(1,1),\ci(2,1),\cd,
\ci(n-1,n-1),\ci(n,n-1))\in(\bbC^\times)^{n(n-1)}$, we have
\begin{eqnarray*}
&&\Del_{w_0\Lm_k,s_k\Lm_k}(\Theta^-_\bfii0(c))\\
&&\hspace{-20pt}=
\ci(1,k)+\frac{\ci(2,k)}{\ci(1,k+1)}+\cd+
\frac{\ci(n-2,k)}{\ci(n-3,k+1)}
+\frac{\ci(n-1,k)\ci(n,k)}{\ci(n-2,k+1)}
+\frac{\ci(n-2,k+1)}{\ci(n-1,k+1)\ci(n,k+1)}+
\frac{\ci(n,k)}{\ci(n-1,k+1)}+\frac{\ci(n-1,k)}{\ci(n,k+1)}
+\frac{\ci(n-3,k+2)}{\ci(n-2,k+2)}+\cd+
\frac{\ci(k,n-1)}{\ci(k+1,n-1)}\\
&&\qq\qq (k=1,2,\cd,n-2),\\
&&\Del_{w_0\Lm_{n-1},s_{n-1}\Lm_{n-1}}(\Theta^-_\bfii0(c))=\ci(n-1,n-1),\qq
\Del_{w_0\Lm_{n},s_n\Lm_{n}}(\Theta^-_\bfii0(c))=\ci(n,n-1).
\end{eqnarray*}
\end{thm}

\begin{thm}\label{thm-Dn2}
Let $k$ be in $\{1,2,\cd,n-2\}$.
Then we have
\begin{eqnarray}
\Del_{w_0s_k\Lm_k,\Lm_k}(\Theta^-_\bfii0(c))
=\frac{1}{\ci(1,k)}+\sum_{j=1}^{k-1}\frac{\ci(k-j,j)}{\ci(k-j+1,j)}.
\end{eqnarray}
\end{thm}
The proof of Theorem \ref{thm-Dn2} is the almost same as the ones
for Theorem \ref{thm-c-2} and Theorem \ref{thm-Bn2}. 

The cases $k=n-1,n$ will be presented in \ref{Dn-del-n}.
We shall prove the above theorems in the next section.
%%%%%%%%%%%%%%%%%%%%%%%%%%%%%%%%%%%%%%%%%
\subsection{Proof of Theorem \ref{thm-Dn1}}

Considering similarly to type $C_n$ as in \ref{proof-c-1}, 
we can write 
\begin{eqnarray*}
{\pmb x}_i(c):=\al_i^\vee(c^{-1})x_i(c)&=&
c^{-h_i}(1+c\cdot e_i),\\
{\pmb y}_i(c):=y_i(c)\al_i^\vee(c^{-1})&=&
(1+c\cdot f_i)c^{-h_i},
\end{eqnarray*}
since $f_i^2=e_i^2=0$ 
on the vector representation $V(\Lm_1)$.
We also have $\omega({\pmb y}_i(c))=\al_i^\vee(c^{-1})x_i(c)=
{\pmb x}_i(c)$ and define the coefficients 
$\xxi(\ovl i,p,j)=\xxi(\ovl i,p,j)(c^{[1:p]})$ and 
$\xxi(i,p,j)=\xxi(i,p,j)(c^{[1:p]})$ for $j\in I$ and
$p\in\{1,2,\cd,n-1\}$
by 
\begin{eqnarray*}
&&X^{(p)}X^{({p-1})}\cd X^{(1)}v_{i}
=\sum_{j=1}^n{\xxi(i,p,j)}v_j
+\sum_{j=1}^n{\xxi(i,p,\ovl j)}v_{\ovl j}\in V(\Lm_1)
\q(i=1,\cd,n),\\
&&X^{(p)}X^{({p-1})}\cd X^{(1)}v_{\ovl i}
=\sum_{j=1}^n{\xxi(\ovl i,p,j)}v_j
+\sum_{j=1}^n{\xxi(\ovl i,p,\ovl j)}v_{\ovl j}\in V(\Lm_1)
\q(i=1,2,\cd,n),
\end{eqnarray*}
where $c^{[1:p]}=(\ci(1,1),\ci(2,1),\cd,\ci(n-1,p),\ci(n,p))$
and $X^{(p)}={\pmb x}_n(\ci(n,p)){\pmb x}_{n-1}(\ci(n-1,p))\cd 
{\pmb x}_1(\ci(1,p))$.

It follows from  (\ref{minor-bilin}) and 
$\omega(\Theta_{{\bf i}_0}(c))=X^{(n-1)}\cd X^{(1)}$ that 
we also have 
\begin{equation}
 \Del_{w_0\Lm_i,s_i\Lm_i}(\Theta_\bfi(c))=\lan \ovl s_i\cdot
 u_{\Lm_i}\,,\,
X^{(n-1)}\cd X^{(1)}v_{\Lm_i}\ran, 
\label{del-x-d}
\end{equation}
which is almost same as (\ref{del-x}).
To describe $\xxi(\ovl i,p,j)$ explicitly
let us define the {\it segments} of type $D_n$, which are similar to 
the ones for $B_n$ and $C_n$.
For $k\in I$  and $p\in\{1,2,\cd,n-1\}$
set $L:=p-n+k+1$ and $S:=n-k$, which are slightly different 
from the other types.
For  
$M\in {\mathcal M}^{(p)}_k[r]
:=\{M=\{m_2,m_3,\cd,m_L\}|2+r\leq m_2<\cd m_L\leq p+r\}$
$(r\in\bbZ)$.
As in the previous sections, we denote $\cM^{(p)}_k[0]$  by
$\cM^{(p)}_k$ .

For $m\in M_j\subset M=M_1\sqcup\cd\sqcup M_S\in{\mathcal M}^{(p)}_k$, 
define $n(m):=n-j+1$.
For $M=M_1\sqcup\cd\sqcup M_S\in{\mathcal M}^{(p)}_k$, write
$M_1=\{2,3,\cd,a\}$. 
For 
$1\leq b\leq c\leq a\leq p$ and $i\leq i_2\leq \cd\leq i_b\leq n-1$
define the monomials in $(\ci(j,i))_{1\leq i<n,1\leq j\leq n}$, 
\begin{eqnarray}
&&C_{i_2,i_3,\cd,i_b}^M:=
\frac{(\ci(i_2,1)\cd\ci(i_b,b-1))({\ci(n,1)}^{\ep_{i_2+1}}\cd 
{\ci(n,b-1)}^{\ep_{i_b+1}})}
{\ci(i_2-1,2)\cd \ci(i_b-1,b)},\qq
D^M:=\prod_{m\in M\setminus M_1}\frac{\ci(n(m)-1,m)}{\ci(n(m),m)},
\label{CD}\\
&&
E^M_{b,c}:=\frac{
\ci(n,b)\ci(n,b+1)\cd\ci(n,c-1)}{\ci(n-1,b+1)\ci(n-1,b+2)\cd\ci(n-1,c)}
+\frac{
\ci(n-1,b)\ci(n-1,b+1)\cd\ci(n-1,c-1)}{\ci(n,b+1)\ci(n,b+2)\cd\ci(n,c)},
\label{EE}\\
&&F^M_{c,a}:=\frac{
\ci(n-2,c+1)\ci(n-2,c+2)\cd \ci(n-2,a)}
{(\ci(n-1,c+1)\ci(n,c+1))(\ci(n-1,c+2)\ci(n,c+2))\cd 
(\ci(n-1,a)\ci(n,a))},
\label{FF}
\end{eqnarray}
where $\ep_i=\del_{i,n}$ and $C_{i_2,i_3,\cd,i_a}^M=1$ 
(resp. $D^M=1$, $E^M_{b,c}=1$, $F^M_{c,a}=1$) if 
$M_1=\emptyset$ or $b=1$ (resp. $M\setminus M_1=\emptyset$, 
$M_1=\emptyset$ or $b=c$, $M_1=\emptyset$ or $c=a$).
%See note 15 pp56, note 17 pp7
\begin{pro}
In the setting above, we have
\begin{eqnarray}
&&\xxi(\ovl i,p,\ovl k)=
\sum_{i\leq i_2\leq \cd\leq i_p\leq k}
\frac{(\ci(i_2,1)\ci(i_3,2)\cd
\ci(i_p,p-1)\ci(k,p){\ci(n,p)}^{\ep'_{k}})
{\ci(n,1)}^{\ep'_{i_2}}{\ci(n,2)}^{\ep'_{i_3}}\cd
{\ci(n,p-1)}^{\ep'_{i_p}}}{\ci(i-1,1)\ci(i_2-1,2)\cd \ci(i_p-1,p)}
%See note 15 pp55 or note 17 pp7-8
\label{d-xi-bar}\\
&&\qq\qq\qq\qq\qq\qq\q(k=1,2,\cd,n), \nn \\
&&\xxi(\ovl i,p,n)=
\sum_{\begin{array}{l}\scriptstyle i\leq i_2\leq \cd\leq i_p\leq n+1\\
\scriptstyle i_2,\cd,i_p\ne n\end{array}}
\frac{(\ci(i_2-2\ep''_{i_2},1)\ci(i_3-2\ep''_{i_3},2)\cd
\ci(i_p-2\ep''_{i_p},p-1)\ci(n-1,p))
{\ci(n,1)}^{\ep'_{i_2}}{\ci(n,2)}^{\ep'_{i_3}}\cd
{\ci(n,p-1)}^{\ep'_{i_p}}}{\ci(i-1,1)\ci(i_2-1,2)\cd \ci(i_p-1,p)}
%See note 15 pp55 or note 17 pp7-8
\label{d-xi-n}\\
&& \xxi(\ovl i,p,k)
={\ci(i-1,1)}^{-1}
\sum_{(A)}
C_{i_2,\cd,i_b}^M\cdot D^M\cdot E^M_{b,c}\cdot F^M_{c,a}\qq
(k=1,2,\cd,n-1),
\label{d-xi}
%See note 17 pp7
\end{eqnarray}
where $\ep'_i=\del_{i,n-1}$, $\ep''_i=\del_{i,n+1}$ and 
the condition (A) is as follows:
\begin{enumerate}
\item[(A)]
$M=M_1\sqcup\cd\sqcup M_S\in{\mathcal M}^{(p)}_k,$
$i\leq i_2\leq\cd\leq i_b\leq n-1, $
$M_1=\{2,\cd,a\},\,\,1\leq b<c\leq a.$
\end{enumerate}
\end{pro}

{\sl Proof.} 
Set $\cX:={\pmb x}_n(c_n)\cd{\pmb x}_1(c_1)$. By calculating directly
we have the formula:
\begin{eqnarray}
&&\cX v_i=\begin{cases}c_1^{-1}v_1&\text{ if }i=1,\\
c_{i-1}c_i^{-1}v_i+v_{i-1}&\text{ if }i=2,\cd,n-2,\\
c_{n-2}c_{n-1}^{-1}c_n^{-1}v_{n-1}+v_{n-2}&\text{ if }i=n-1,\\
c_{n-1}c_n^{-1}v_n+c_{n}^{-1}v_{n-1}&\text{ if }i=n,\\
\end{cases}\\
&&\cX v_{\ovl i}=\begin{cases}
c_{i-1}^{-1}(c_iv_{\ovl i}+\cd +c_{n-2}v_{\ovl{n-2}}+
c_{n-1}c_nv_{\ovl{n-1}}+c_nv_{\ovl n}+c_{n-1}v_n+v_{n-1})&\text{ if }i<n,
\\
c_nc_{n-1}^{-1}v_{\ovl n}+c_{n-1}^{-1}v_{n-1}&\text{ if }i=n.
\end{cases}
\end{eqnarray}
%note 15 pp52
where we understand $c_0=1$. Using these, we get 
\begin{eqnarray}
&&\xxi(\ovl i,p,\ovl k)=\sum_{j=i}^k\xxi(\ovl i,p-1,\ovl j)
\frac{\ci(k,p){\ci(n,p)}^{\ep'_{k}}}{\ci(j-1,p)},
\q(k=1,\cd,n),
\label{recur-d1}\\
&&
\xxi(\ovl i,p,k)=\xxi(\ovl i,p-1,k+1)+\xxi(\ovl i,p-1,k)
\frac{\ci(k-1,p)}{\ci(k,p)},
\q(k=1\cd,n-2),
\label{recur-d2}\\
&&\xxi(\ovl i,p,n-1)=\sum_{j=i}^n\xxi(\ovl i,p-1,\ovl j)
{\ci(j-1,p)}^{-1}+\xxi(\ovl i,p-1,n){\ci(n,p)}^{-1}
+\xxi(\ovl i,p-1,n-1)\frac{\ci(n-2,p)}{\ci(n-1,p)\ci(n,p)},
%See note 15 pp56 note 17 pp8-9
\label{recur-d3}\\
&&\xxi(\ovl i,p,n)=\left(\sum_{j=i}^{n-1}\xxi(\ovl i,p-1,\ovl j)
{\ci(j-1,p)}^{-1}\ci(n-1,p)\right)+
\xxi(\ovl i,p-1,n){\ci(n-1,p)}{{\ci(n,p)}^{-1}}.
\label{recur-d4}
%note 15 pp56 or note 17 pp7
\end{eqnarray}
Indeed, the formulae (\ref{d-xi-bar}) and 
(\ref{d-xi-n})
are easily proved by the induction
on $p$ using the formulae (\ref{recur-d1}) and (\ref{recur-d4})
as the other types.

Considering  
the segments in $\cM_{k}^{(p)}$, 
$\cM_{k+1}^{(p-1)}$ and $\cM_{k}^{(p-1)}$ as the cases for $C_n$ and 
applying the recursions (\ref{recur-d2}) and (\ref{recur-d3}) to the
induction hypothesis, 
we obtain (\ref{d-xi}).\qed

Thus, for example, we have 
\[
\Del_{w_0\Lm_1,s_1\Lm_1}(\Theta^-_{\bfii0}(c))=\xxi(\ovl 1,n-1,2)
=\sum_{j=1}^{n-2}\frac{\ci(j,1)}{\ci(j-1,2)}+
\frac{{\ci(n-1,1)}\ci(n,1)}{\ci(n-2,2)}+\frac{\ci(n,1)}{\ci(n-1,2)}+
\frac{\ci(n-1,1)}{{\ci(n,2)}}+\frac{\ci(n-2,2)}{\ci(n-1,2)\ci(n,2)}+
\sum_{j=3}^{n-1}\frac{\ci(n-j,j)}{\ci(n-j+1,j)}.
\]

The following is similar to Lemma \ref{lem-det-c} 
and Lemma \ref{lem-det-b}.
\begin{lem}
For $k=1,\cd, n-2$ we define the matrix $W_k$ by 
\begin{equation}
W_k:=
\begin{pmatrix}
\xxi(\ovl 1,n-1,k+1)&\xxi(\ovl 1,n-1,k-1)&\cd&\xxi(\ovl 1,n-1,2)&
\xxi(\ovl 1,n-1,1)\\
\xxi(\ovl 2,n-1,k+1)&\xxi(\ovl 2,n-1,k-1)&\cd&\xxi(\ovl 2,n-1,2)&
\xxi(\ovl 2,n-1,1)\\
\vdots&\vdots&\cd&\vdots&\vdots\\
\xxi(\ovl k,n-1,k+1)&\xxi(\ovl k,n-1,k-1)&\cd&\xxi(\ovl k,n-1,2)&
\xxi(\ovl k,n-1,1)
\end{pmatrix}.
\label{det-d}
\end{equation}
Then we have $\Del_{w_0\Lm_k,s_k\Lm_k}(\Theta^-_\bfii0(c))=\det W_k$.
\end{lem}

Similar to the previous cases,
the last column of the matrix $W_k$ is given as 
\[
{}^t(\xxi(\ovl 1,n-1,1),\cd,\xxi(\ovl k-1,n-1,1))=
 {}^t(1,{\ci(1,1)}^{-1},{\ci(2,1)}^{-1},\cd,{\ci(k-1,1)}^{-1}).
\]
Then, applying the elementary transformations on $W_k$ by 
($i$-th row) $-\frac{\ci(i,1)}{\ci(i-1,1)}\times$($i+1$-th row)
 for $i=1,\cd, k-1$, in the transformed matrix $\til W_k$
its $(i,j)$-entry is as follows:
\begin{lem} For $c=(\ci(k,l)|1\leq k\leq n,1\leq l< n)$ 
we set $c^{(l)}:=(\ci(1,l),\cd,\ci(n,l))$ and 
$c^{[a:b]}:=(c^{(a)},c^{(a+1)},\cd,c^{(b)})$ $(a\leq b)$.
For $i=1,\cd, k-1$  the $(i,j)$-entry $(\til W_k)_{i,j}$ is:
\begin{equation}
(\til W_k)_{i,j}=\begin{cases}
\xxi(\ovl i,n-2,k+1)(c^{[2:n-1]})&\text{ if }j=1,\\
\xxi(\ovl i,n-2,k-j+1)(c^{[2:n-1]})&\text{ if }1<j<k,\\
0&\text{ if }j=k,
\end{cases}
\end{equation}
\end{lem}
Note that for type $D_n$ by definition we have $c=c^{[1:n-1]}$.\\
{\sl Proof.}
The proof is similar to the one for Lemma \ref{lem-w-c}. Indeed, 
for type $D_n$ we also get the same formula as \eqref{xi-xi}:
\begin{equation*}
\xxi(\ovl i,n-1,j)(c)-\frac{\ci(i,1)}{\ci(i-1,1)}\xxi(\ovl{i+1},n-1,j)(c)
=\frac{\ci(i,1)}{\ci(i-1,1)} \xxi(\ovl i,n-2,j)(c^{[2,n-1]}).
\end{equation*}
Then, this proves the lemma.
\qed

Applying the above elementary transformations to the matrix $W_k$, we
find 
\begin{eqnarray*}
\det W_k
=\xxi(\ovl 1,n-k,k+1)(c^{[k:n-1]})
\end{eqnarray*}
Thus, it follows from  (\ref{d-xi})  that for $k=1,2,\cd,n-2$
\begin{eqnarray}
&&\Del_{w_0\Lm_k,s_k\Lm_k}(\Theta^-_\bfi(c))
=\xxi(\ovl 1,n-k,k+1)(c^{[k:n-1]}) \nn \\
&&=\ci(1,k)+\frac{\ci(2,k)}{\ci(1,k+1)}+\cd+
\frac{\ci(n-2,k)}{\ci(n-3,k+1)}
+\frac{\ci(n-1,k)\ci(n,k)}{\ci(n-2,k+1)}
+\frac{\ci(n-2,k+1)}{\ci(n-1,k+1)\ci(n,k+1)}+
\frac{\ci(n,k)}{\ci(n-1,k+1)}+\frac{\ci(n-1,k)}{\ci(n,k+1)}
+\frac{\ci(n-3,k+2)}{\ci(n-2,k+2)}+\cd+
\frac{\ci(k,n-1)}{\ci(k+1,n-1)}.\nn
\end{eqnarray}
%Note 15, pp60-61
The cases $k=n-1,n$ are easily obtained by 
the formula in \cite[(4.18)]{BZ2}:
\begin{equation}
\Del_{w_0\Lm_{n-1},s_k\Lm_{n-1}}(\Theta^-_\bfii0(c))=\ci(n-1,n-1),\qq
\Del_{w_0\Lm_n,s_k\Lm_n}(\Theta^-_\bfii0(c))=\ci(n,n-1).
\end{equation}
Now, the proof of Theorem \ref{thm-Dn1} has been done.
\qed
%%%%%%%%%%%%%%%%%%%%%%%%%%%%%%%%%%%%%%%%%
\subsection{Correspondence to the monomial realizations}%9.3

Except for $\Del_{w_0s_{n-1}\Lm_n,\Lm_{n-1}}(\Theta^-_\bfii0(c))$ 
and $\Del_{w_0s_n\Lm_n,\Lm_n}(\Theta^-_\bfii0(c))$
(see \ref{Dn-del-n}), we shall see
the positive answer to the conjecture for type $D_n$.

Let us see the monomial realization of $B(\Lm_1)$ associated with 
the cyclic sequence $\cd(12\cd n)(12\cd n)\cd$, which means that the sign
$p_0=(p_{i,j})$ is given by $p_{i,j}=1$ if $i<j$ and $p_{i,j}=0$ if 
$i>j$ as in the previous sections.
The crystal $B(\Lm_1)$ is described as follows:
We abuse the notation 
$B(\Lm_1):=\{v_i,v_{\ovl i}|1\leq i\leq n\}$ if 
there is no confusion.
Then the actions of $\eit$ and $\fit$ $(1\leq i\leq n)$ are defined as 
$\fit=f_i$ and $\eit=e_i$ in (\ref{d-fi}) and (\ref{d-fn}).

To see the monomial realization $B(\ci(1,k))$, 
we give the explicit forms of the monomials $A_{i,m}$:
% See Note 17 pp4-5
\begin{equation}
A_{i,m}=\begin{cases}
\ci(1,m){\ci(2,m)}^{-1}\ci(1,m+1)&\text{ for }i=1,\\
\ci(i,m){\ci(i+1,m)}^{-1}{\ci(i-1,m+1)}^{-1}\ci(i,m+1)
&\text{ for }1<i< n-2,\\
\ci(n-2,m){\ci(n-1,m)}^{-1}{\ci(n,m)}^{-1}{\ci(n-3,m+1)}^{-1}
\ci(n-2,m+1)&\text{ for }i=n-2,\\
\ci(n-1,m){\ci(n-2,m+1)}^{-1}
\ci(n-1,m+1)&\text{ for }i=n-1,\\
\ci(n,m){\ci(n-2,m+1)}^{-1}\ci(n,m+1)&\text{ for }i=n.
\end{cases}
\label{aim-d-1}
\end{equation}
Let 
where $m^{(k)}_i:B(\Lm_i)\hookrightarrow \cY(p)$ ($u_{\Lm_i}\mapsto 
\ci(i,k)$)
is the embedding of crystal as in Sect.5.
Here the monomial realization $B(\ci(1,k))=m^{(k)}_1(B(\Lm_1))=
\{m^{(k)}_1(v_j),m^{(k)}_1(v_{\ovl j})|1\leq j\leq n\}$ associated with 
$p_0$ is described explicitly:
\begin{equation}
m^{(k)}_1(v_j)=\begin{cases}
\frac{\ci(j,k)}{\ci(j-1,k+1)}&1\leq j\leq n-2,\\
\frac{\ci(n-1,k)\ci(n,k)}{\ci(n-2,k+1)}&j=n-1,\\
\frac{\ci(n,k)}{\ci(n-1,k+1)}&j=n,
\end{cases}\qq\qq
m^{(k)}_1(v_{\ovl j})=\begin{cases}
\frac{\ci(j-1,k+n-i)}{\ci(j,k+n-j)}&1\leq j\leq n-2,\\
\frac{\ci(n-2,k+1)}{\ci(n-1,k+1)\ci(n,k+1)}&j=n-1,\\
\frac{\ci(n-1,k)}{\ci(n,k+1)}&j=n,
\end{cases}
\end{equation}
where we understand $\ci(0,k)=1$.
Now, Theorem \ref{thm-Dn1} and Theorem \ref{thm-Dn2} claim the following:
\begin{thm}
We obtain 
\begin{eqnarray}
&&\Del_{w_0\Lm_k,s_k\Lm_k}(\Theta^-_\bfii0(c))
=\sum_{j=1}^nm^{(k)}_1(v_j)+\sum_{j=k+1}^nm^{(k)}_1(v_{\ovl j}),\qq
(k=1,2\cd,n-2)
\label{thm-mono-d1}\\
&&\Del_{w_0\Lm_{n-1},s_{n-1}\Lm_{n-1}}(\Theta^-_\bfii0(c))
=m^{(n-1)}_{n-1}(u_{\Lm_{n-1}}),\qq
\Del_{w_0\Lm_{n},s_k\Lm_{n}}(\Theta^-_\bfii0(c))=m^{(n-1)}_n(u_{\Lm_n}),\\
&&\Del_{w_0s_k\Lm_k,\Lm_k}(\Theta^-_\bfii0(c))
=\sum_{j=1}^km^{(k-n+1)}_1(v_{\ovl j}), \qq(k=1,2,\cd,n-2).
\label{thm-mono-d2}
\end{eqnarray}
\end{thm}
Note that the last result is derived from the fact that
$B(\ci(1,k))=B({\ci(1,n+k-1)}^{-1})$, which is the connected component
including ${\ci(1,n+k-1)}^{-1}$ as the lowest monomial.

%%%%%%%%%%%%%%%%%%%%%%%%%%%%%%%%%%%%%%%%%%
\subsection{$\Del_{w_0s_{n-1}\Lm_{n-1},\Lm_{n-1}}(\Theta^-_\bfii0(c))$ 
and $\Del_{w_0s_n\Lm_n,\Lm_n}(\Theta^-_\bfii0(c))$}\label{Dn-del-n}
%9.4
%See Note 16 pp56, 60-71
To state the results for 
$\Del_{w_0s_{n-1}\Lm_n,\Lm_{n-1}}(\Theta^-_\bfii0(c))$ and 
$\Del_{w_0s_n\Lm_n,\Lm_n}(\Theta^-_\bfii0(c))$,
we need to prepare the set of triangles $\tri'_n$ for type $D_n$
which is similar to the one for type $B_n$:
\begin{equation}
\tri'_n:=\{(\ji(k,l)|1\leq k\leq l<
n)|1\leq \ji(k,l+1)\leq \ji(k,l)<\ji(k+1,l+1)\leq n
\q(1\leq k\leq l<n-1)\}.
\end{equation}
We visualize a triangle $(\ji(k,l))$ in $\tri'_n$ as follows:
\[(\ji(k,l))=
 \begin{array}{c}
\ji(1,1)\\
\ji(2,2)\ji(1,2)\\
\ji(3,3)\ji(2,3)\ji(1,3)\\
\cd\cd\cd\cd\\
\ji(n-1,n-1)\cd\ji(2,n-1)\ji(1,n-1)
\end{array}
\]
Here we know that the set of triangles $\tri'_n$ for $D_n$
coincides with $\tri_{n-1}$ for $B_n$.
As type $B_n$ we easily obtain 
\begin{lem}
For any $k\in\{1,2,\cd,n-1\}$ there exists a unique $j$
($1\leq j\leq k+1$) such that the $k$th row of 
a triangle $(\ji(k,l))$ in $\tri_n'$ is in the following form:
\begin{equation}
k\text{-th row}\q
(\ji(k,k),\ji(k-1,k),\cd,\ji(2,k),\ji(1,k))
=(k+1, k,k-1,\cd,j+1,j-1,j-2,\cd,2,1),
\end{equation}
that is, we have $\ji(m,k)=m$ for $m<j$ and 
$\ji(m,k)=m+1$ for $m\geq j$.
\end{lem}
For a triangle $\del=(\ji(k,l))\in\tri'_n$, we list $j$'s as in the lemma:
$s(\del):=(s_1,s_2,\cd,s_{n-1})$, 
which we call the {\it label} of a triangle $\del$.
Here we have the following same as Lemma \ref{s4}:
\begin{lem}\label{s4-d}
For any $\del\in\tri_n'$ 
let $s(\del):=(s_1,s_2,\cd,s_{n-1})$ be its label.
Then 
\begin{enumerate}
\item
The label $s(\del)$  satisfies 
$1\leq s_k\leq k+1$ and $s_{k+1}=s_k$ or $s_k+1$ for $k=1,\cd,n-1$.
\item
Each  $k$-th row of a triangle $\del$ is in 
one of the following
     {\rm I, II, III, IV:}
\begin{enumerate}[\rm I.]
\item $s_{k+1}=s_k+1$ and $s_k=s_{k-1}$.
\item $s_{k+1}=s_k$ and $s_k=s_{k-1}$.
\item $s_{k+1}=s_k+1$ and $s_k=s_{k-1}+1$.
\item $s_{k+1}=s_k$ and $s_k=s_{k-1}+1$.
\end{enumerate}
Here we suppose that $s_0=1$ and $s_{n}=s_{n-2}+1$, which means that 
the 1st row must be in {\rm I,II or IV} 
and the $n-1$-th row is in {\rm I or IV}.
\end{enumerate}
\end{lem}
Now, we associate a Laurant 
monomial $m(\del)$ in variables $(\ci(i,j))_{i\in I,j\in\bbZ}$ 
with a triangle $\del=(\ji(k,l))$ by the
following way.
\begin{enumerate}
\item Let $s=(s_1\cd,s_{n-1})$  be the label of $\del\in\tri'_n$.
\item Suppose  $i$-th row is in the form I. If $1\leq i\leq n-2$, 
 then associate $\ci(i,s_i)$. For $i=n-1$, 
\begin{enumerate}
\item If $n+ s_{n-1}$ is even, then associate $\ci(n-1,s_{n-1})$.
\item If $n+ s_{n-1}$ is odd, then associate $\ci(n,s_{n-1})$.
\end{enumerate}
\item Suppose  $i$-th row is in the form IV. If $1\leq i\leq n-2$, 
 then associate ${\ci(i,s_i)}^{-1}$. For $i=n-1$, 
\begin{enumerate}
\item If $n+ s_{n-1}$ is even, then associate ${\ci(n,s_{n-1})}^{-1}$.
\item If $n+ s_{n-1}$ is odd, then associate ${\ci(n-1,s_{n-1})}^{-1}$.
\end{enumerate}
\item If  $i$-th row is in the form II or III, then associate 1.
\item  Take the product of all monomials as above for 
$1\leq i <n$, then
we obtain the monomial $m(\del)$ associated with $\del$. 
This defines the map 
$m:\tri'_n\to \cY$, where $\cY$ is the set of Laurant monomials 
in $(\ci(i,j))_{i\in I,j\in\bbZ}$.
\end{enumerate}
Here we define the involutions $\xi$ and $\ovl{\,\,}$ on $\cY$ by 
\begin{equation}
\begin{array}{l}
\xi:\ci(k,n-1)\mapsto \ci(k,n),\q \ci(k,n)\mapsto \ci(k,n-1),
\ci(k,j)\mapsto \ci(k,j) \q(j\ne n-1,n),\\
{\,}^-:\ci(i,j)\mapsto {\ci(i,n-j)}^{-1}.
\end{array}
\label{op-d}
\end{equation}
As type $B_n$, 
let us denote the special triangle such that 
$\ji(k,l)=k+1$(resp. $\ji(k,l)=k$) for any $k,l$ by $\del_h$ 
(reps. $\del_l$).
Indeed, we have
\begin{equation}
m(\del_h)=\begin{cases}
\ci(n-1,1)&\text{ if $n$ is odd,}\\
\ci(n,1)&\text{ if $n$ is even,}
\end{cases},\qq
m(\del_l)={\ci(n,n)}^{-1}.
\label{m-del-d}
\end{equation}
Now, we present $\Del_{w_0s_{n-1}\Lm_{n-1},\Lm_{n-1}}(\Theta^-_\bfii0(c))$ 
and $\Del_{w_0s_n\Lm_n,\Lm_n}(\Theta^-_\bfii0(c))$ for type $D_n$:
\begin{thm}\label{thm-Lm-n-d} 
For type $D_n$, we have the explicit forms:
\begin{eqnarray}
&&\Del_{w_0s_n\Lm_n,\Lm_n}(\Theta^-_\bfii0(c))=
\sum_{\del\in\tri_n\setminus\{\del_l\}}\ovl{m(\del)},\\
&&\Del_{w_0s_{n-1}\Lm_{n-1},\Lm_{n-1}}(\Theta^-_\bfii0(c))=
\sum_{\del\in\tri_n\setminus\{\del_l\}}\ovl{\xi\circ m(\del)}.
\end{eqnarray}
\end{thm}
The proof of this theorem will be given in \ref{proof-thm-Lm-n-d}.
\begin{ex} 
The set of triangles $\tri'_5$ is the same as $\tri_4$:
\begin{eqnarray*}
&&\begin{array}{c}
1\\21\\321\\4321
\end{array}\q
\begin{array}{c}
1\\21\\321\\5321
\end{array}\q
\begin{array}{c}
1\\21\\421\\5321
\end{array}\q
\begin{array}{c}
1\\31\\421\\5321
\end{array}\q
\begin{array}{c}
2\\31\\421\\5321
\end{array}\q
\begin{array}{c}
2\\31\\421\\5421
\end{array}\q
\begin{array}{c}
1\\31\\421\\5421
\end{array}\q
\begin{array}{c}
1\\21\\421\\5421
\end{array}\\
&&\begin{array}{c}
1\\31\\431\\5421
\end{array}\q
\begin{array}{c}
2\\31\\431\\5421
\end{array}\q
\begin{array}{c}
2\\32\\431\\5421
\end{array}\q
\begin{array}{c}
1\\31\\431\\5431
\end{array}\q
\begin{array}{c}
2\\31\\431\\5431
\end{array}\q
\begin{array}{c}
2\\32\\431\\5431
\end{array}\q
\begin{array}{c}
2\\32\\432\\5431
\end{array}\q
\begin{array}{c}
2\\32\\432\\5432
\end{array}
\end{eqnarray*}
and their labels $s(\del)$ are
\begin{eqnarray*}
&&(2,3,4,5),\q (2,3,4,4),\q(2,3,3,4),\q(2,2,3,4),\q(1,2,3,4),\q
(1,2,3,3),\q(2,2,3,3),\q(2,3,3,3),\\
&&(2,2,2,3),\q (1,2,2,3),\q(1,1,2,3),\q(2,2,2,2),\q(1,2,2,2),\q
(1,1,2,2),\q(1,1,1,2),\q(1,1,1,1).
\end{eqnarray*}
Then, we have the corresponding monomials ${m(\del)}$:
\begin{eqnarray*}
&& \frac{1}{\ci(5,5)},\q
\frac{\ci(5,4)}{\ci(3,4)},\q \frac{\ci(3,3)}{\ci(2,3)\ci(4,4)},\q
\frac{\ci(2,2)}{\ci(1,2)\ci(4,4)},\q\frac{\ci(1,1)}{\ci(4,4)},\q
\frac{\ci(1,1)\ci(4,3)}{\ci(3,3)},\q\frac{\ci(2,2)\ci(4,3)}{\ci(1,2)\ci(3,3)},\q
\frac{\ci(4,3)}{\ci(2,3)},\\
&&\frac{\ci(3,2)}{\ci(1,2)\ci(5,3)},\q
\frac{\ci(1,1)\ci(3,2)}{\ci(2,2)\ci(5,3)},\q\frac{\ci(2,1)}{\ci(5,3)},\q
\frac{\ci(5,2)}{\ci(1,2)},\q \frac{\ci(1,1)\ci(5,2)}{\ci(2,2)},\q
\frac{\ci(2,1)\ci(5,2)}{\ci(3,2)},\q \frac{\ci(3,1)}{\ci(4,2)},\q
{\ci(4,1)}.
\end{eqnarray*}
And then, we have the corresponding monomials $\ovl{m(\del)}$:
\begin{eqnarray*}
&& \ci(5,0),\q
\frac{\ci(3,1)}{\ci(5,1)},\q \frac{\ci(2,2)\ci(4,1)}{\ci(3,2)},\q
\frac{\ci(1,3)\ci(4,1)}{\ci(2,3)},\q\frac{\ci(4,1)}{\ci(1,4)},\q
\frac{\ci(3,2)}{\ci(1,4)\ci(4,2)},\q\frac{\ci(1,3)\ci(3,2)}{\ci(2,3)\ci(4,2)},\q
\frac{\ci(2,2)}{\ci(4,2)},\\
&&\frac{\ci(1,3)\ci(5,2)}{\ci(3,3)},\q
\frac{\ci(2,3)\ci(5,2)}{\ci(1,4)\ci(3,3)},\q\frac{\ci(5,2)}{\ci(2,4)},\q
\frac{\ci(1,3)}{\ci(5,3)},\q \frac{\ci(2,3)}{\ci(1,4)\ci(5,3)},\q
\frac{\ci(3,3)}{\ci(2,4)\ci(5,3)},\q \frac{\ci(4,3)}{\ci(3,4)},\q
\frac{1}{\ci(4,4)}.
\end{eqnarray*}
Thus, we have 
\begin{eqnarray*}
&&\Del_{w_0s_n\Lm_n,\Lm_n}(\Theta^-_\bfii0(c))=
\frac{\ci(3,1)}{\ci(5,1)}+ \frac{\ci(2,2)\ci(4,1)}{\ci(3,2)}+
\frac{\ci(1,3)\ci(4,1)}{\ci(2,3)}+\frac{\ci(4,1)}{\ci(1,4)}+
\frac{\ci(3,2)}{\ci(1,4)\ci(4,2)}+\frac{\ci(1,3)\ci(3,2)}{\ci(2,3)\ci(4,2)}+
\frac{\ci(2,2)}{\ci(4,2)}\\
&&\qq\qq\qq\qq+\frac{\ci(1,3)\ci(5,2)}{\ci(3,3)}+
\frac{\ci(2,3)\ci(5,2)}{\ci(1,4)\ci(3,3)}+\frac{\ci(5,2)}{\ci(2,4)}+
\frac{\ci(1,3)}{\ci(5,3)}+ \frac{\ci(2,3)}{\ci(1,4)\ci(5,3)}+
\frac{\ci(3,3)}{\ci(2,4)\ci(5,3)}+ \frac{\ci(4,3)}{\ci(3,4)}
+\frac{1}{\ci(4,4)}.
\end{eqnarray*}
\end{ex}
%%%%%%%%%%%%%%%%%%%%%%%%%%%%%%%%%%%%%%
\subsection{Crystal structure on $\tri'_n$}

We shall define certain crystal structure on $\tri'_n$ by the similar
way to type $B_n$. 
First, let us define
the actions of $\fit$ and $\eit$ as follows:
For a triangle $\del\in \tri'_n$ let
 $J_k=(\ji(k,k),\ji(k-1,k),\cd,\ji(2,k),\ji(1,k))$ be the $k$-th
row of $\del$. Thus, we denote $\del=(J_1,\cd,J_{n-1})$.
It follows from  Lemma \ref{s4-d} that 
there exists a unique $j$ such that 
$J_k=(k+1,k,\cd,j+1,j-1,\cd,2,1)$ and $J_k$ is in one of I,II,III,IV.
Set $J'_k=(k+1,\cd,j+2,j,j-1,\cd,2,1)$ and 
$J''_k=(k+1,\cd,j+1,j,j-2,\cd,2,1)$.
Then, we have
\begin{eqnarray}
&&\fit\del=\begin{cases}
(\cd,J_{i-1},J'_i,J_{i+1},\cd)&\text{ if $J_i$ is in I},\\
\qq 0&\text{otherwise,}
\end{cases}\qq(i=1,\cd,n-2),\\
&&\til f_{n-1}\del=\begin{cases}
(\cd,J_{n-2},J'_{n-1})&\text{ if $J_{n-1}$ is in I and $n+j$ is even},\\
\qq 0&\text{otherwise,}
\end{cases}\\
&&\til f_{n}\del=\begin{cases}
(\cd,J_{n-2},J'_{n-1})&\text{ if $J_{n-1}$ is in I and $n+j$ is odd},\\
\qq 0&\text{otherwise,}
\end{cases}\\
&&\eit\del=\begin{cases}
(\cd,J_{i-1},J''_i,J_{i+1},\cd)&\text{ if $J_i$ is in IV},\\
\qq 0&\text{otherwise,}
\end{cases}\qq(i=1,\cd,n-2),\\
&&\til e_{n-1}\del=\begin{cases}
(\cd,J_{n-2},J''_{n-1})&\text{ if $J_{n-1}$ is in IV and $n+j$ is odd},\\
\qq 0&\text{otherwise,}
\end{cases}\\
&&\til e_{n}\del=\begin{cases}
(\cd,J_{n-2},J''_{n-1})&\text{ if $J_{n-1}$ is in IV and $n+j$ is even},\\
\qq 0&\text{otherwise.}
\end{cases}
\end{eqnarray}
The weight of $\del=(J_1,\cd,J_{n-1})$ is defined as follows:
Let $s=(s_k)_{k=1,\cd,n-1}$ be the label of $\del$, that is, 
$J_k=(k+1,k,\cd, j+1,j-1,\cd,2,1)$ for $j=s_k$.
\begin{equation}
\wt(\del)=\begin{cases}
\displaystyle
\Lm_n-\sum_{k=1}^{n-2}(s_k-1)\al_k-\left[\frac{s_{n-1}-1}{2}\right]
\al_{n-1}
-\left[\frac{s_{n-1}}{2}\right]\al_{n}&\text{if $n$ is even,}\\
\displaystyle
\Lm_{n-1}-\sum_{k=1}^{n-2}(s_k-1)\al_k-\left[\frac{s_{n-1}}{2}\right]
\al_{n-1}
-\left[\frac{s_{n-1}-1}{2}\right]\al_{n}&\text{if $n$ is odd,}
\end{cases}
\end{equation}
where $[x]$ is the so-called Gaussian symbol, {\it i.e.,}, the maximum
integer which does not exceed $x$.
We can easily check that $\tri'_n$ is equipped with 
the crystal structure and obtain.
\begin{pro}\label{del-spin-d}
We have the following isomorphism of crystals:
\[
\tri'_n\cong \begin{cases}B(\Lm_n)&\text{if $n$ is even,}\\
B(\Lm_{n-1})&\text{if $n$ is odd.}
\end{cases}
\]
The highest weight crystal is $\del_h\in\tri'_n$.
\end{pro}
{\sl Proof.}
As was given in \ref{Dn}, we know the explicit form of the crystals
$B(\Lm_n)\cong B^{(+,n)}_{sp}$ and $B(\Lm_{n-1})\cong B^{(-,n)}_{sp}$. 
So, we shall see how 
$\tri'_n$ and $B^{(\pm,n)}_{sp}$ correspond to each other.
Let $s=(s_1,\cd,s_{n-1})$ be the label of a triangle $\del\in \tri'_n$.
Now, let us associate $\ep=(\ep_1,\cd,\ep_n)\in B^{(\pm,n)}_{sp}$ 
with $s$ by 
\begin{eqnarray}
&&\ep_1=\begin{cases}+&\text{ if }s_1=1,\\
-&\text{ if }s_1=2,
\end{cases}\qq
\ep_k=\begin{cases}+&\text{ if }s_{k}=s_{k-1},\\
-&\text{ if }s_{k}=s_{k-1}+1,
\end{cases}(2\leq k\leq n-2),
\label{ep1-d}\\
&&
(\ep_{n-1},\ep_n)=\begin{cases}
(+,+)&\text{if $s_{n-1}=s_{n-2}$ and $n+s_{n-1}$ is odd.}\\
(+,-)&\text{if $s_{n-1}=s_{n-2}$ and $n+s_{n-1}$ is even,}\\
(-,+)&\text{if $s_{n-1}=s_{n-2}+1$ and  $n+s_{n-1}$ is odd,}\\
(-,-)&\text{if $s_{n-1}=s_{n-2}+1$ and $n+s_{n-1}$ is even,}
\end{cases}
\label{ep2-d}
\end{eqnarray}
which define the map $\psi_+:\tri'_n\to B^{(+,n)}_{sp}$ if $n$ is even
and $\psi_-:\tri'_n\to B^{(-,n)}_{sp}$ if $n$ is odd.
Then, e.g., for even $n$, 
the vector $\del_h$ (resp. $\del_l$) corresponds to the 
highest (resp. lowest) weight vector $(+,+,\cd,+)$
(resp. $(-,-,\cd,-)$) since $s(\del_h)=(1,1,\cd,1)$ and 
$s(\del_l)=(2,3,\cd,n)$.
It is clear to find that the map $\psi_{\pm}$ is bijective.
The rest of the proof is almost the same as the one for 
Proposition \ref{del-spin-b}.
\qed

Next, let us show that the map $m:\tri'_n\to\cY$ gives an 
isomorphism between $\tri'_n$ and the monomial
realization of $B(\Lm_n)$ for even $n$ or $B(\Lm_{n-1})$ for odd $n$.
Consider the crystal structure on $\cY$
by taking the sign $\ovl p=(\ovl p_{i,j})$ which is the same one as (\ref{ovl-pij}).
Then,  it corresponds to the cyclic sequence of $I$ such as:
$\cd(n n-1 \cd 21)(n n-1 \cd 21)\cd$.
Here, by the map $m$ we  get
\[
m(\del_h)=\begin{cases}\ci(n,1)&\text{if $n$ is even},\\
\ci(n-1,1)&\text{if $n$ is odd.}
\end{cases}
\]
\begin{pro}\label{mono-del-d}
Let $B(\ci(n,1))$(resp. $B(\ci(n-1,1))$) 
be the connected subcrystal of $\cY(\ovl p)$ whose
highest monomial is $\ci(n,1)$ (resp. $\ci(n-1,1)$) 
and which is isomorphic to $B(\Lm_n)$ (resp.$B(\Lm_{n-1})$)
of type $D_n$. Then we have $m(\tri'_n)= B(\ci(n,1))$ if $n$ is even and 
$m(\tri'_n)= B(\ci(n-1,1))$ if $n$ is odd. 
\end{pro}
{\sl Proof.}
By the above setting for $\cY(\ovl p)$, we have 
\begin{equation}\label{aim-d}
A_{i,m}=\begin{cases}
\ci(1,m){\ci(2,m+1)}^{-1}\ci(1,m+1),&\text{ for }i=1,\\
\ci(i,m){\ci(i-1,m)}^{-1}{\ci(i+1,m+1)}^{-1}\ci(i,m+1),
&\text{ for }1<i<n-2,\\
\ci(n-2,m){\ci(n-3,m)}^{-1}{\ci(n-1,m+1)}^{-1}{\ci(n,m+1)}^{-1}
\ci(n-2,m+1),
&\text{ for }i=n-2,\\
\ci(n-1,m){\ci(n-2,m)}^{-1}\ci(n-1,m+1),&\text{ for }i=n-1,\\
\ci(n,m){\ci(n-2,m)}^{-1}\ci(n,m+1)&\text{ for }i=n.
\end{cases}
\end{equation}
For $\del\in\tri'_n$ and 
$i=1,2,\cd,n-2$ we can see that $m(\fit\del)=\fit(m(\del))$ by the same
way as in the proof of Proposition \ref{mono-del-b}.

Thus, let us see $\til f_{n-1}$ and $\til f_n$.
First, for the label $s(\del)=(s_1,\cd,s_{n-1})$, 
suppose that $n+s_{n-1}$ is odd.
Thus, the action of $\til f_{n-1}$ is trivial and 
\[
 \til f_{n}:(s_{n-3},s_{n-2},s_{n-1})=
\begin{array}{ccccccc}
 &(j-1,j,j)&\to&(j-1,j,j+1)&\q \frac{\ci(n,j)}{\ci(n-2,j)}
&\to&\frac{1}{\ci(n,j+1)}\\
 &(j,j,j)&\to&(j,j,j+1)&\q \ci(n,j)&\to&\frac{\ci(n-1,j)}{\ci(n,j+1)}.
\end{array}
\]
This means that the action of $\til f_{n}$ involves the 
multiplication of the 
monomial 
 ${A_{n,j}}^{-1}={\ci(n,j)}^{-1}\ci(n-2,j){\ci(n,j+1)}^{-1}$.
Therefore, we have $m(\til f_n\del)=\til f_nm(\del)$.

For $s=(s_1,\cd,s_{n-1})$, suppose that $n+s_{n-1}$ is even.
Thus, the action of $\til f_{n}$ is trivial and 
\[
 \til f_{n-1}:(s_{n-3},s_{n-2},s_{n-1})=
\begin{array}{ccccccc}
 &(j-1,j,j)&\to&(j-1,j,j+1)&\q \frac{\ci(n-1,j)}{\ci(n-2,j)}
&\to&\frac{1}{\ci(n-1,j+1)}\\
 &(j,j,j)&\to&(j,j,j+1)&\q \ci(n-1,j)&\to&\frac{\ci(n-1,j)}{\ci(n-1,j+1)},
\end{array}
\]
which means that the action of $\til f_{n-1}$ involves the 
multiplication of the 
monomial 
 ${A_{n-1,j}}^{-1}={\ci(n-1,j)}^{-1}\ci(n-2,j)
{\ci(n-1,j+1)}^{-1}$.
Therefore, we have 
$m(\til f_{n-1}\del)=\til f_{n-1}m(\del)$.
As for the actions of $\eit$, we can show that $m(\eit \del)=\eit
m(\del)$ by the similar way to the case of $\fit$. \qed

Then, by these results we have 
\begin{thm}
$B(\ci(n,1))=m(\tri'_n)\cong B(\Lm_n)$ if $n$ is even and 
$B(\ci(n-1,1))=m(\tri'_n)\cong B(\Lm_{n-1})$ if $n$ is odd. 
\end{thm}

%%%%%%%%%%%%%%%%%%%%%%%%%%%%%%%%%%%
\subsection{Proof of Theorem \ref{thm-Lm-n-d}}\label{proof-thm-Lm-n-d}
%9.6

By the explicit descriptions in (\ref{D-fi}) and (\ref{D-ei}),
we know that $e_i^2=f_i^2=0$ on $V^{(\pm,n)}_{sp}$.
Thus, we can write 
${\pmb x}_i(c):=\al_i^\vee(c^{-1})x_i(c)=c^{-h_i}(1+c\cdot e_i)$ 
and ${\pmb y}_i(c):=y_i(c)\al_i^\vee(c^{-1})=(1+c\cdot f_i)c^{-h_i}$ 
on 
$V^{(\pm,n)}_{sp}$
and then for $\ep=(\ep_1,\cd,\ep_n)\in B^{(\pm,n)}_{sp}$
\begin{equation}
\pmbx_i(c)\ep=\begin{cases}
c\ep+\ep'&\text{ if }
(\ep_i,\ep_{i+1})=(-,+),\\
c^{\frac{\ep_{n-1}-\ep_n}{2}}\ep&\text{ otherwise },
\end{cases}\,\,(i<n),\qq
\pmbx_n(c)\ep=\begin{cases}c\ep+\ep''
&\text{ if }(\ep_{n-1},\ep_n)=(-,-),\\
c^{-\frac{\ep_{n-1}+\ep_n}{2}}\ep&\text{ otherwise },
\end{cases}
\label{xi-Dn}
\end{equation}
%% Note 16 pp56
where $\ep'=(\cd,\ep_{i-1}, +,-,\ep_{i+2},\cd)$ for $i<n$ and 
$\ep''=(\cd,\ep_{n-2},+,+)$.

For $c^{(i)}=(\ci(1,i),\cd,\ci(n,i))\in (\bbC^\times)^n$ set
$X^{(i)}(c^{(i)}):=\pmbx_1(\ci(1,i))\cd \pmbx_n(\ci(n,i))$ and 
$X(c):=X^{(n-1)}(c^{(n-1)})\cd X^{(1)}(c^{(1)})$ where 
$c:=(c^{(n-1)},\cd,c^{(1)})$.

To obtain the explicit form of 
$\Del_{w_0s_n\Lm_n,\Lm_n}(\Theta^-_\bfii0(c))$, 
let us see 
$\Del_{w_0\Lm_n,s_n\Lm_n}(\Theta^-_{\bfii0^{-1}}(c))$
as in \ref{proof-thm-Lm-n}, 
since
for $g\in U\ovl w_0 U$ we have 
\begin{equation}\label{eta-si-d}
\Del_{w_0s_i\Lm_i,\Lm_i}(g)=\Del_{w_0\Lm_i,s_i\Lm_i}(\eta(g)),
\end{equation}
and $\eta(\Theta^-_{\bfii0}(c))=\Theta^-_{\bfii0^{-1}}(\ovl c)$
$(-:\ci(j,i)\mapsto\ci(n-j+1,i))$ as the previous case.
%See Note 16 p90-91.
Since $\omega(\Theta^-_{\bfii0^{-1}}(c))=X(c)$, 
we have
$\Del_{w_0\Lm_n,s_n\Lm_n}(\Theta^-_{\bfii0^{-1}}(c))
=\lan s_n\ep_h,X(c)\ep_l\ran$
where $s_n\ep_h=s_n(+,\cd,+)=(+,\cd,+,-,-)$ and 
\[
\ep_l=\begin{cases}(-,\cd,-,-)&\text{ if $n$ is even,}\\
(-,\cd,-,+)&\text{ if $n$ is odd.}
\end{cases}
\]
Thus, our aim is to obtain the coefficient of the vector 
$s_n\ep_h=(+,\cd,+,-,-)$ in $X(c)\ep_l$.

For the sequence $\bfii0=(12\cd n)^{n-1}$, 
let $J$ be the set of all subsequences of $\bfii0$.
The operator $X(c)$ is expanded in the form on $V_{sp}^{(\pm)}$:
\begin{equation}\label{expand-d}
X(c)=\sum_{\xi\in J}\al_\xi e_\xi,
\end{equation}
where $e_\xi=e_{i_1}\cd e_{i_m}$ for $\xi=(i_1,\cd,i_m)$ and $\al_\xi$ 
is a product of $\al_i^\vee(c)$'s
and some scalars. 
Now, we shall see which sequence $(i_1,\cd,i_k)\in J$ gives
\begin{equation}\label{non0-d}
 e_{i_1}\cd e_{i_k}\ep_l=(+,\cd,+,-,-)=s_n\ep_h.
\end{equation}
As for the weight, we have
\begin{equation}\label{wt-dif-d}
\wt(s_n\ep_h)-\wt(\ep_l)=
\sum_{i=1}^{n-2}i\al_i+\left[\frac{n-1}{2}\right]\al_{n-1}+
\left[\frac{n-2}{2}\right]\al_{n},
%\begin{cases}
%\sum_{i=1}^{n-2}i\al_i+\frac{n-2}{2}\al_{n-1}+\frac{n-2}{2}\al_{n}&
%\text{ if $n$ is even,}\\
%\sum_{i=1}^{n-2}i\al_i+\frac{n-1}{2}\al_{n-1}+\frac{n-3}{2}\al_{n}&
%\text{ if $n$ is odd.}
%\end{cases}
\end{equation}
where $[x]$ is the Gaussian symbol as before.
Thus, for $(i_1,\cd,i_k)$ in (\ref{non0-d}), we know that
if $i\ne n-1,n$ the number of $i$'s is $i$ and the number of $n-1$
(resp. $n$) is $\frac{n-2}{2}$ (resp. $\frac{n-2}{2}$) if $n$ is even or
$\frac{n-1}{2}$ (resp. $\frac{n-3}{2}$) if $n$ is odd.

Let us associate an element in $U({\mathfrak n}_+)$ 
with a triangle $\del=(\ji(k,l))\in \tri'_n$ by the following way,
which is almost same as the one for the 
type $B_n$ in \ref{proof-thm-Lm-n}:
\begin{enumerate}
\item
For $i=1,\cd,n-1$ define $r_i(\del)$ to be the set of indices as
$r_i(\del):=\{(k,l)|\ji(k,l)=i\}$ and 
$n_i(\del):=\sharp r_i(\del)$.
\item
Set $L=n_i(\del)$. 
For $r_i(\del)$, define $l_{a_1},\cd,l_{a_L}$ by setting 
$r_i(\del)=\{(k_{a_1},l_{a_i}),\cd,(k_{a_{L}},l_{a_L})\}$
with $l_{a_1}<l_{a_2}<\cd<l_{a_L}$.
\item For $i=1,2,\cd, n-1$ set 
\[
 E_i(\del):=e_{l_{a_1}}e_{l_{a_2}}\cd e_{l_{a_{L-1}}}\times
\begin{cases}e_{l_{a_L}}&\text{ if $l_{a_L}\ne n-1$},\\
e_{n-1}&\text{if $l_{a_L}=n-1$ and $k_{a_L}+l_{a_L}$ is odd,}\\
e_{n}&\text{if $l_{a_L}=n-1$ and $k_{a_L}+l_{a_L}$ is even.}
\end{cases}
\]
and $E(\del):=E_{n-1}(\del)\cdot E_{n-2}(\del)\cd E_1(\del)$.
\end{enumerate}
\begin{ex}
For $n=5$ and a triangle 
\[
\del=\begin{array}{c}
2\\32\\431\\5421
\end{array}\qq
\text{we obtain }\q
\begin{array}{ll}
r_1(\del)=\{(1,3),(1,4)\},&E_1(\del)=e_3e_4,\\
r_2(\del)=\{(1,1),(1,2),(2,4)\},&E_2(\del)=e_1e_2e_5\\
r_3(\del)=\{(2,2),(2,3)\},&E_3(\del)=e_2e_3\\
r_4(\del)=\{(3,3),(3,4)\},&E_4(\del)=e_3e_4,
\end{array}
\]
and then $E(\del)=(e_3e_4)(e_2e_3)(e_1e_2e_5)(e_3e_4)$. 
This example is quite similar to Example \ref{EX-B}, but $r_2(\del)$
 differs from the one in Example \ref{EX-B}.
\end{ex}

\begin{lem}\label{tri-ep-d}
Let $(j_1,\cd,j_m)\in J$ be a subsequence of $\bfii0$ 
and $E=e_{j_1}\cd e_{j_m}$ be the associated
monomial of $e_i$'s. Then we have 
\begin{equation}
E\ep_l=s_n\ep_h=(+,\cd,+,-,-)\text{ if and only if there exists } 
\del\in \tri_n\setminus\{\del_l\}\text{ such that }
E=E(\del).
\end{equation}
\end{lem}
Note that $\del=(\ji(k,l))\ne\del_l$ if and only if $\ji(n-1,n-1)=n$.

\nd
{\sl Proof.}
The proof is similar to the one for Lemma \ref{tri-ep}.
Indeed, we can easily see that $E(\del_l)\ep_l\ne s_n\ep_h$.
Suppose that there exists $\del\in\tri'_n\setminus\{\del_l\}$ such that 
$E=E(\del)$.
Writing $E(\del)=e_{l_m}\cd e_{l_1}$ we may 
show that $e_{l_k}\cd e_{l_1}\ep_l\ne0$ for any $k=1,2,\cd,m$ by the
induction as in the case $B_n$.
For this $\del$, suppose that $e_{l_k}$ appears in $E_j(\del)$. 
The cases $i=l_k=1,2,\cd,n-2$ are done by the same way as the proof 
of Lemma \ref{tri-ep}. 
So, let us see the case $i=l_k=n-1$ or $n$. 
The $n-2$-th row and the $n-1$-th row of $\del$ around $j$ are as
follows:
\[
 (i)
\begin{array}{cc}
\syl n-2\text{-th row}&\syl *\q j-1\q\\
\syl n-1\text{-th row}&\q \textcircled{$\scriptstyle j$}\q \syl j-1
\end{array}\qq
(ii)
\begin{array}{cc}
\syl n-2\text{-th row}&\syl *\q j-1\q\\
\syl n-1\text{-th row}&\q \textcircled{$\scriptstyle j$}\q \syl j-2
\end{array}\qq
(iii)
\begin{array}{cc}
\syl n-2\text{-th row}&\syl *\q j-2\q\\
\syl n-1\text{-th row}&\q \textcircled{$\scriptstyle j$}\q \syl j-2
\end{array}
\]
where $*$ means $j$ or $j+1$ and 
$\textcircled{$\syl j$}$ is the entry which gives $e_{l_k}=e_{n-1}$
or $e_n$.

As for (i), we have the following cases:
\begin{enumerate}
\item[(a)]
$E_j(\del)=\cd e_n\cd$ and 
$E_{j-1}(\del)=\cd e_{n-2}e_{n-1}\cd$.
\item[(b)]
$E_j(\del)=\cd e_{n-1}\cd$ and 
$E_{j-1}(\del)=\cd e_{n-2}e_{n}\cd$.
\end{enumerate} 
Now, we consider the case (a).
Suppose that $j-1$ in the $n$-th row gives $e_{n-1}=e_{l_p}$ 
in $E_{j-1}(\del)$.
Denote the vector before being applied 
$e_{n-1}=e_{l_p}$ by $\ep'=(\ep'_i)$, that is, 
$\ep'=e_{l_{p-1}}\cd e_{l_1}\ep_l$, 
which satisfies $\ep'_{n-2}=\ep'_{n-1}=-,\ep'_{n}=+$ and then 
\[
 (\ep'_{n-2}, \ep'_{n-1},\ep'_n)=(-,-,+)\mapright{e_{n-1}}
(-,+,-)\mapright{e_{n-2}}(+,-,-).
\]
This shows that the vector before being applied $e_n=e_{l_k}$, 
say $\ep''=(\ep''_i)$, has the form $\ep''_{n-1}=\ep''_{n}=-$ and then 
we have $e_n\ep''=e_{l_k}\cd e_{l_1}\ep_l\ne 0$.
The case (b) is also shown similarly and the cases (ii) and (iii)
are also shown similarly.
Therefore, we obtain $E(\del)\ep_l\ne0$ for any $\del\in\tri_n'\setminus\{\del_l\}$.

Assuming $E\ep_l=s_n\ep_h=(+,\cd,+,-,-)$ we see that each 
$e_i$ ($i\ne n$) should appear in $E$ $i$-times $(1\leq i\leq n-2)$,
$e_{n-1}$ should appear $\left[\frac{n-1}{2}\right]$-times
and $e_n$ should appear $\left[\frac{n-2}{2}\right]$-times 
by (\ref{wt-dif-d}). Moreover, $e_{n-1}$ and $e_n$ appear alternatively
since only $e_n$ can change $\ep_n=-$ to $+$ and only $e_{n-1}$ can
change $\ep_n=+$ to $-$.
In the sequence $\bfii0$ there are $n-1$-cycles $12\cd n$ just as
\[
\begin{array}{ccccc}
\bfii0=&  \underbrace{12\cd n}&
\cd  &\underbrace{12\cd n}&\underbrace{12\cd n}\\
&(n-1)\text{th cycle}&&2\text{nd cycle}&1
\text{st cycle}
\end{array}
\]
Since the indices $n-1$  and $n$ appear alternatively 
$\left[\frac{n-1}{2}\right]+\left[\frac{n-2}{2}\right]=n-2$-times
in $\bfii0$, 
we should choose $n-1$ or $n$ from one cycle alternatively. 
Listing the number of 
such cycles from the right, we obtain 
$(n-1,\cd,j+1,j-1,\cd,2,1)
=(\ji(n-2,n-1),\ji(n-3,n-1),\cd,\ji(2,n-1),\ji(1,n-1))$ for 
some $j$, which is 
the bottom row($=n-1$-th row) of an element in $\tri'_n$ except the
first entry $\ji(n-1,n-1)=n$.
Next, we define the $n-2$-th row.
We can apply $e_{n-2}$ after applying $e_{n-1}$ or $e_n$ and can apply
$e_{n-1}$ or $e_n$ after applying $e_{n-2}$. Thus, this means 
that $n-2$ must be in between $n-1$ and $n$ chosen in the above process.
Then, we obtain the $n-2$-th row 
$(\ji(n-2,n-2),\ji(n-3,n-2),\cd,\ji(2,n-2),\ji(1,n-2))$, which 
should satisfy the condition $\ji(k,n-1)\leq \ji(k,n-2)<\ji(k+1,n-1)$
$(1\leq k\leq n-2)$ since this $n-2$ can be in the same cycle as 
the previous $n-1$ or $n$.
This condition just coincides with the ones for $\tri'_n$.
Repeating this process, we obtain a triangle $\del$ 
and find that this $\del\in\tri'_n$ satisfies $E=E(\del)$.
\qed

Next, let us find the coefficient of the vector 
generated from $\ep_l$ by 
applying $a_\xi e_{i_1}\cd e_{i_k}$ 
in (\ref{expand}).
Set ${\bf E}_i(c)=c\al_i^\vee(c^{-1})e_i$. Then we can write
$\pmbx_i(c)=\al_i^{\vee}(c^{-1})+{\bf E}_i(c)$.
We consider the cases that $\del$ is in I$\sim$IV.
In the case that $\del$ is in I,
there exists $j$ such that the $i-1$, $i$, $i+1$th rows of $\del$
are in the form:
\[
\begin{array}{lcccc}
i-1\text{-th row}&\cd j+3\q j+2\q j+1\q j-1\q j-2\\
i\text{-th row}&\cd \q\q j+3\q j+2\q j+1\q j-1\q j-2\\
i+1\text{-th row}&\cd\q\q\q\q j+3\q j+2\q j\q\q j-1\q j-2
\end{array}
\]
This means that in the expansion (\ref{expand-d}) we know that 
this $\del$ gives the following monomial in which $\al_i^{\vee}$
appears once:
\begin{equation}
\cd \al_i^\vee{({\ci(i,j)}^{-1})}{\bf E}_{i+1}(\ci(i+1,j))\cd
{\bf E}_{i-1}(\ci(i-1,j-1)){\bf E}_{i}(\ci(i,j-1))
{\bf E}_{i+1}(\ci(i+1,j-1))\cd 
\end{equation}
Hence, by the action of $\al_i^\vee{({\ci(i,j)}^{-1})}$
we obtain the coefficient $\ci(i,j)$.
In  the other cases II,III,IV, we can discuss similarly and 
obtain that for II and III, we have the coefficient 1 and 
for IV we have ${\ci(i,j)}^{-1}$, which coincides with the recipe
after Lemma \ref{s4-d} and then we know that 
the desired coefficient is the same as $m(\del)$.

As has been seen in Proposition \ref{del-spin-d}, the set of monomials 
$m(\tri'_n)$ has the crystal structure 
isomorphic to $B(\Lm_{n-1})$ or $B(\Lm_n)$.
The following lemma is immediate from Theorem \ref{thm-dual}.
\begin{lem}
The set of monomials $\ovl{m(\tri'_n)}$ (resp. $\ovl{\xi\circ m(\tri_n')}$)
is a crystal isomorphic to
$B(\ci(n,0))\cong B(\Lm_n)$ (resp. $B(\Lm_{n-1})$), where ${\,}^-\,$ and 
$\xi$ are the involutions as in \eqref{op-d}.
\end{lem}
{\sl Proof.}
We find by the definition of the map ${}^-$ 
that this is the case $a=n$ in Theorem \ref{thm-dual}.
By Proposition \ref{del-spin-d}, we know that $\tri'_n\cong 
B(\ci(n,1))$ if $n$ is even and 
$\tri'_n\cong  B(\ci(n-1,1))$ if $n$ is odd.
Thus, by Theorem \ref{thm-dual} we have that 
for an even $n$, $\ovl{B(\ci(n,1))}=B({\ci(n,n-1)}^{-1})=B(\ci(n,0))$
and 
for an odd $n$,
$\ovl{B(\ci(n-1,1))}=B({\ci(n-1,n-1)}^{-1})=B(\ci(n,0))$,
which means the desired result.\qed

Thus, by the formula (\ref{eta-si-d}), we obtain 
\[
 \Del_{w_0s_n\Lm_n,\Lm_n}(\Theta^-_\bfii0(c))=
\Del_{w_0\Lm_n,s_n\Lm_n}(\Theta^-_{\bfii0^{-1}}(\ovl c)),
\]
where for $c=(c^{(1)},\cd,c^{(n)})$ ($c^{(i)}=(\ci(l,i))_{1\leq l\leq n}$)
we define $\ovl c=(\ovl{c^{(n-1)}},\cd,\ovl{c^{(1)}})$ 
$(\ovl{c^{(i)}}=({\ci(n-l+1,i)}^{-1})_{1\leq l\leq n})$.
The case $k=n-1$
\[
 \Del_{w_0s_{n-1}\Lm_{n-1},\Lm_{n-1}}(\Theta^-_\bfii0(c))=
\Del_{w_0\Lm_{n-1},s_{n-1}\Lm_{n-1}}(\Theta^-_{\bfii0^{-1}}(\ovl c)),
\]
is also obtained by considering the map $\xi$, that is, by flipping
$n-1\leftrightarrow n$.
Therefore, we have completed the proof of 
Theorem \ref{thm-Lm-n-d}.
\qed 

Hence, we get the positive answer to our conjecture for type $D_n$.
\bibliographystyle{amsalpha}

\end{document}